\PassOptionsToPackage{dvipsnames}{xcolor}
\documentclass[11pt]{article}
\usepackage[margin=1in]{geometry}
\usepackage{authblk}

\usepackage{url} 
\usepackage{lineno}
\usepackage{soul}
\usepackage{amsmath}
\usepackage{array}
\usepackage{amsthm}

\usepackage{graphicx}
\usepackage{float}
\usepackage{microtype}
\usepackage{amsbsy}
\usepackage{unicode-math}
\usepackage{comment}
\usepackage{xcolor}
\usepackage{algorithm}
\usepackage{algpseudocode}
\usepackage{hyperref}
\usepackage{apacite}
\let\cite\shortcite
\usepackage{cleveref}
\hypersetup{breaklinks = true,
	        colorlinks = true,
	        linkcolor  = blue,
            urlcolor = RoyalPurple,
	        citecolor  = Blue}

\usepackage[inline]{trackchanges} 
\renewcommand{\boldsymbol}[1]{\symbfit{#1}}
\newcommand{\vect}[1]{\symbfit{#1}}
\newcommand{\bu}{{\symbf{u}}}
\newcommand{\mO}{\mathcal{O}}
\newcommand{\e}{\varepsilon}

\newcommand{\bx}{{\symbf{x}}}
\newcommand{\bz}{{\symbf{z}}}
\newcommand{\by}{{\symbf{y}}}
\newcommand{\bU}{{\symbf{U}}}
\newcommand{\bG}{{\symbf{G}}}
\newcommand{\bnabla}{\symbf{\nabla}}
\newcommand{\pd}[2]{\frac{\partial #1}{\partial #2}}

\newcommand{\defn} {\ensuremath{\stackrel{\mathrm{def}}{=}}}

\providecommand{\citep}{\cite}

\date{}
\begin{document}

\title{A conservative Runge--Kutta framework for\\nonlinear free-surface ocean models}

\author[1,2]{Simone Silvestri}
\author[2]{Jean-Michel Campin}
\author[2,3]{Gregory Wagner}
\author[4,5]{Navid~C.~Constantinou}
\author[7]{Siddhartha Bishnu}
\author[2,6]{Xin Kai Lee}
\author[2,6]{Raffaele Ferrari}

\affil[1]{DIATI, Politecnico di Torino, Torino, Italy}
\affil[2]{Earth, Atmospheric, and Planetary Sciences, Massachusetts Institute of Technology, Cambridge, MA, USA}
\affil[3]{Aeolus Labs, San Francisco, CA, USA}
\affil[4]{University of Melbourne, Parkville, VIC, Australia}
\affil[5]{Australian Research Council Center of Excellence for the Weather of the 21st Century, Australia}
\affil[6]{Center for Computational Science and Engineering, Massachusetts Institute of Technology}
\affil[7]{Department of Earth Sciences, University of Cambridge}

\newcommand{\correspondence}{Correspondence to S. Silvestri (\href{mailto:silvestri.simone0@gmail.com}{silvestri.simone0@gmail.com})}

\maketitle
\begin{center}\small\correspondence\end{center}
\begin{abstract}
We present a Runge--Kutta framework for nonlinear free-surface ocean models implemented in Oceananigans.
We consider two formulations that differ in the time discretization of the coupled barotropic mode: \emph{(i)}~an implicit free-surface discretization and \emph{(ii)}~a split--explicit free-surface discretization.
We test two three-stage, low-storage Runge--Kutta formulations, one of which is strong stability preserving, and a four-stage formulation.
We analyze the stability of the resulting schemes in a vertical normal-mode framework by separating the coupling effects from the uncoupled evolution, and show how stability depends on the barotropic time discretization, the averaging filter, and the reconstruction of the slow forcing.
We identify configurations that are stable and second-order accurate in time in the linearized analysis.
By adapting exact variance-budget diagnostics to multi-stage integrators, we directly quantify the contribution of temporal discretization to numerical mixing.
In a nearly linear internal tide and in an eddying re-entrant channel, the Runge--Kutta formulations produce less numerical mixing than the Adams--Bashforth time discretization with the same spatial discretization, while using timesteps 3.5 to 5.5 times larger.
This reduction does not occur in a quarter-degree global simulation, where we suspect that the spatial discretization dominates the numerical mixing, although the split--explicit Runge--Kutta schemes still increase throughput by 14--25\%.

\end{abstract}

\section*{Plain Language Summary}
Ocean models simulate how currents, temperature, and salinity evolve over time.
Advancing these simulations is challenging because fast surface waves must be treated together with the slower motions of the ocean interior.
We develop and test a family of time-stepping methods that combine these motions, and examine the choices needed to keep the simulations stable and accurate.
In two idealized ocean configurations, these methods also reduce artificial mixing introduced by the numerical calculations by up to 50\% compared with a commonly used method.
In a global simulation at coarser resolution, the mixing is similar across the time-stepping schemes, possibly because spatial discretization errors dominate, but the simulation runs faster.
The framework is available in the Oceananigans ocean model, providing a practical approach to improving the accuracy of ocean simulations.
\section*{Key Points}
\begin{itemize}
\item Low-storage and strong-stability-preserving Runge--Kutta formulations preserve tracer conservation and constancy at every coupled stage.
\item Barotropic averaging with low-dissipation filters and a reconstruction of the slow forcing stabilize split-explicit Runge--Kutta schemes.
\item Runge–Kutta time-stepping reduces numerical mixing relative to Adams–Bashforth in idealized tests while taking larger timesteps.
\end{itemize}

\section{Introduction}
\label{sec:1}
Ocean general circulation models must resolve processes spanning vastly different timescales, from external gravity waves propagating at hundreds of meters per second to eddy advection and internal waves occurring over days to weeks. This separation, typically of two to three orders of magnitude, has driven the development of split--explicit and implicit-explicit (IMEX) approaches that decouple the fast barotropic dynamics from the slower three-dimensional evolution \cite{berntsen1981, blumberg1987, bleck1990, killworth1991, hallberg1997, higdon1997, campin2004, shchepetkin}.
The baroclinic component, governing currents, tracers, and stratification, advances on timesteps constrained by advection and internal gravity waves, while the barotropic component employs either sub-cycling or an implicit solver to handle the fast external gravity waves efficiently (e.g., \citeA{griffies2001}).

Baroclinic time-stepping in ocean models has long relied on low-order methods, often with filtering or modified coefficients to improve stability.
Leapfrog schemes, typically used in atmospheric models \cite{klower2024} and in ocean models such as LICOM \cite{licom3}, offer nominally non-dissipative, second-order accurate integration, but require filtering to control the computational mode, which introduces numerical dissipation \cite{asselin1972frequency}. Before version~5.0, NEMO used a modified leapfrog formulation, incorporating filtering and a modified coupling between tracer and momentum updates to improve stability \cite{madec2017}.
The ROMS models avoid plain leapfrog by using third-order Adams--Bashforth and leapfrog--Adams--Moulton predictor--corrector schemes, with the choice depending on the model branch and prognostic variable \cite{shchepetkin, shchepetkin2009commentary}.
MOM6 \cite{adcroft2019} instead adopts a two-stage predictor-corrector approach based on a modified version of the second-order Runge--Kutta method, which adds the damping that plain RK2, unstable for purely oscillatory problems, lacks.
The Adams--Bashforth method, which avoids time-staggering between state variables, remains the most prevalent choice, used in MITgcm \cite{marshall1997}, ICON-O \cite{korn2017}, Veros \cite{hafner2018}, the current MPAS-Ocean \cite{mpaso2026code} and Oceananigans \cite{oceananigans}.
These ocean models mostly use over-extrapolated, formally first-order variants of the second-order Adams--Bashforth scheme; an exception is the third-order Adams--Bashforth scheme \cite{durran1990} implemented in the MITgcm, which achieves higher accuracy at the same computational cost \cite{eccov4}. However, the latter is rarely used since \emph{(i)}~it requires substantially more memory and \emph{(ii)}~it reduces the stability limit of diffusive processes, pairing quite poorly with high-order diffusive spatial discretizations such as WENO schemes \cite{simoneWENO}.

Multi-stage schemes such as Runge--Kutta (RK) methods offer a promising alternative, providing higher-order accuracy and larger stable Courant numbers while also reducing implicit dissipation without requiring state variable staggering.
Low-storage RK formulations, which require the same memory as a second-order Adams--Bashforth method, have been successfully implemented in atmospheric models \cite{wicker, skamarock2008, skamarock2012, baldauf2011, clima_atmos} and are increasingly being explored for ocean applications \cite{ishimwe2023, lanrk3, ducousso2026}.
Avoiding the storage of multiple time levels of the whole state vector is particularly important for GPU-based ocean models, such as Oceananigans, where memory leanness is a strict requirement \cite{silvestri2025}.

However, coupling multi-stage schemes with the evolution of the free surface requires satisfying three constraints simultaneously. 
First, \textit{free surface reconciliation} requires consistency between the explicitly evolved free surface and the grid stretching diagnosed from the continuity equation in response to the free surface displacements; violations lead to spurious mass sources and unphysical vertical velocities \cite{campin2004}.
Second, \textit{tracer conservation} must hold both globally (preserving the total tracer content) and locally (keeping constant tracer fields constant), which becomes particularly challenging with moving vertical coordinates, where tracer and continuity discretizations must match exactly; \citeA{griffies2001} exposed this tension between global and local conservation without fully resolving it.
Third, \textit{numerical stability} of the coupled barotropic--baroclinic system requires careful analysis of mode interactions, as improper coupling can destabilize either fast or slow modes. Foundational analyses of this stability problem were developed by \citeA{higdon1996}, \citeA{higdon1997}, and \citeA{hallberg1997}.

Several recent studies have advanced RK implementations for ocean models.
\citeA{ishimwe2023} developed a robust split--explicit RK3 scheme paired with a Discontinuous Galerkin (DG) spatial discretization, using a three-stage RK scheme for the primitive equations coupled with a two-stage RK scheme for the barotropic mode. They show that their method is stable and tracer-conserving. In contrast, \citeA{lanrk3} implemented the classical three-stage strong-stability-preserving (SSP) Runge--Kutta formulation for both the three-dimensional evolution and the free-surface evolution in the MPAS-Ocean framework.
Most recently, NEMO has transitioned from its modified leapfrog scheme to a third-order Runge--Kutta method in version 5.0, citing improved stability and computational efficiency \cite{nemo5}. This Runge--Kutta method uses a novel ``first-stage'' coupling strategy where the barotropic mode is evolved only once per timestep \cite{ducousso2026}.
Runge--Kutta integration has also been combined with a semi-implicit free-surface treatment in WAVETRISK-OCEAN \cite{kevlahan2022}.

We focus on comparing this new framework with the quasi-second-order Adams--Bashforth (AB2) method for two reasons: AB2 is widely used across ocean models, and it is the only current time discretization for the Oceananigans hydrostatic model \cite{silvestri2025}, the ocean model we aim to improve.
AB2 offers significant computational advantages: it requires only one tendency evaluation per timestep and minimal memory overhead by storing just two time levels.
However, AB2 has a rather small stability range and introduces implicit numerical dissipation.
While the impact of spatial discretization on numerical mixing, the spurious diffusion of tracers introduced by discretization errors, has been extensively studied, and spurious diapycnal mixing from advection schemes has been shown to rival or exceed explicitly parameterized physical mixing \cite{klingbeil2014, hill2012}, the contribution of time discretization to numerical mixing has received comparatively less attention.
The dissipative nature of AB2 suggests that time-stepping schemes may also contribute to spurious mixing in realistic ocean configurations and degrade model fidelity by affecting ocean stratification over long integrations~\cite{ilicak2012}.

To quantify numerical mixing, several diagnostic tools have been developed.
\citeA{klingbeil2014} introduced methods for quantifying the relative contributions of physical and numerical mixing from the tracer variance decay. \citeA{hill2012} used the spreading of virtual tracer release experiments.
\citeA{ilicak2012} introduced the reference potential energy (RPE) framework, an exact, globally integrated measure of numerical mixing, and subsequently extended this approach to quantify the spatial distribution of numerical mixing~\cite{ilicak2016}.
More recently, \citeA{banerjee2024} extended the variance budget method to provide diagnostics that are exact locally in both space and time, a property essential for isolating the contribution of time discretization to numerical mixing.

In this study, we develop low-storage RK time-stepping schemes for the baroclinic evolution and couple them with both the implicit (RK-IM) and the split--explicit (RK-SE) free-surface time discretizations, advancing the barotropic mode at every stage as in the time-split scheme of \citeA{wicker} and in \citeA{hallberg1997}. 
The contribution of this work is threefold. First, we formulate the coupled stage so that any barotropic solver returns a free-surface increment and a transport that satisfy the discrete continuity equation, which guarantees free-surface reconciliation, tracer conservation and constancy preservation at every stage of the Runge--Kutta discretization.
Second, we provide a detailed linear stability analysis of the coupled barotropic--baroclinic system, extending the modal decomposition approach of~\citeA{demange} by separating the coupling contributions from the uncoupled evolution, which turns the choice of barotropic integrator, averaging filter and slow-forcing reconstruction into explicit stability conditions. Third, we adapt the exact variance budget of \citeA{banerjee2024} to multi-stage RK schemes and to the multi-level AB2 time discretization, isolating the contribution of time discretization to numerical mixing.
We quantify this contribution in two idealized test cases, the topographic generation of an internal tide and an equilibrated channel flow, and we assess robustness and computational cost in a quarter-degree global ocean simulation.

The remainder of this paper is organized as follows. In section~\ref{sec:2} we introduce two three-stage Runge–Kutta formulations and a four-stage extension, and compare their stability and damping with AB2 for a generic ODE.
In section~\ref{sec:3} we detail the implementation of the primitive equations, addressing free-surface consistency and tracer conservation for both implicit and split--explicit treatments, followed in section~\ref{sec:4} by the linear stability analysis of the coupled barotropic--baroclinic system, which Appendix~\ref{sec:rotation-results} extends to rotation.
We finally evaluate the proposed schemes in section~\ref{sec:5} and conclude in section~\ref{sec:conclusions}.

\section{Runge--Kutta time discretizations}
\label{sec:2}

We seek to replace the quasi-second-order Adams--Bashforth scheme used in Oceananigans' hydrostatic model \cite{silvestri2025, oceananigans} with low-storage Runge–Kutta schemes for the slow baroclinic component, coupled to a separate solver for the fast barotropic component.
Two three-stage formulations are good candidates: the low-storage formulation of \citeA{wicker}, adapted to the ocean by \citeA{ishimwe2023} and adopted by NEMO 5.0 \cite{nemo5} with the ``first-stage'' coupling strategy described in \citeA{ducousso2026}, and the strong-stability-preserving formulation of \citeA{shuosher}, which underlies the split--explicit scheme that \citeA{lanrk3} implemented in MPAS-Ocean. 

Before describing the complete scheme, we compare the outer discretizations without barotropic coupling. These stability and damping properties provide the reference for the coupled analysis in~Section~\ref{sec:4}. For a generic state $\varphi$ and its tendency $G(\varphi)$,
\begin{equation}
    \frac{\mathrm{d}\varphi}{\mathrm{d}t} = G(\varphi) \ , 
    \label{eq:ODE}
\end{equation}
the quasi-second-order Adams--Bashforth (AB2) scheme used in Oceananigans reads
\begin{equation}
    \varphi^{n+1} = \varphi^n +\Delta t \left[ \left( \frac{3}{2} + \xi\right) G(\varphi^n) - \left( \frac{1}{2} + \xi\right) G(\varphi^{n-1}) \right] \ ,
\end{equation}
where $\xi$ is the over-extrapolation parameter. For $\xi>0$, the scheme is formally first order, with a leading truncation error proportional to $\xi$. A positive $\xi$ is needed to stabilize oscillatory modes, but introduces numerical diffusion in the solution. In the rest of the manuscript we use Oceananigans' default, $\xi=0.1$, and retain the ``AB2'' name (even though the scheme is first order).

A Runge--Kutta scheme advances the solution from $t^n$ to $t^{n+1} = t^n + \Delta t$ through a sequence of stages. We write $\varphi^{n,m}$ for the output of stage $m$, with $\varphi^{n,0}\equiv\varphi^n$ and $\varphi^{n+1}$ the output of the last stage. The low-storage scheme of~\citeA{wicker}, referred to as WRK3 from here on, reads
\begin{equation}
    \varphi^{n,m} = \varphi^{n} + c_m\,\Delta t\,G\left(\varphi^{n,m-1}\right) \ , \qquad m = 1, 2, 3 \ , \qquad
    c = \left(\frac13,\ \frac12,\ 1 \right) \ ,
    \label{eq:ws-RK}
\end{equation}
where the stage fraction $c_m$ specifies the stage end-time $t^n+c_m\Delta t$. Each stage restarts from $\varphi^n$ and uses the tendency evaluated at the preceding stage. Therefore, the final increment uses only the last stage tendency. WRK3 is third order for linear problems, but only second order for nonlinear ones~\cite{wicker}. 

The strong-stability-preserving scheme of~\citeA{shuosher}, SRK3 from here on, is third order for nonlinear problems and can be written in the ``Shu--Osher'' form, as a convex combination of forward-Euler steps,
\begin{equation}
    \begin{gathered}
    \varphi^{n,m} = \left(1 - b_m\right)\varphi^{n} + b_m\left[\varphi^{n,m-1} + \Delta t\,G\left(\varphi^{n,m-1}\right)\right] \ , \qquad m = 1, 2, 3 \ , \qquad
    b = \left(1,\ \frac14,\ \frac23\right) \ .
    \end{gathered}
    \label{eq:ssp-shuosher}
\end{equation}
The SRK3 stage states, corresponding to intermediate approximations that are combined in the final state, evolve the solution to times $t^n + \Delta t$, $t^n + \Delta t/2$ and $t^{n} +\Delta t$. Unrolling the stage recursion, it is possible to express the final update as a weighted sum of the stage tendencies,
\begin{equation}
    \varphi^{n+1} = \varphi^{n} + \Delta t\sum_{m=1}^{3}\beta_m\,G\left(\varphi^{n,m-1}\right) \ , \qquad
    \beta_m = \prod_{l \ge m} b_l \ , \qquad \beta = \left(\frac16,\ \frac16,\ \frac23\right) \ .
    \label{eq:ssp-weights}
\end{equation}
Since every stage is a convex combination of the beginning-of-step state and a forward-Euler step, any bound that a forward-Euler step preserves at a given timestep is preserved by SRK3 at the same timestep (the “strong-stability-preserving” property), a guarantee that WRK3 does not provide. This is a useful property when this scheme is paired with non-oscillatory spatial reconstructions, such as those used in Oceananigans \cite{simoneWENO}.

We add a third formulation which has the potential to increase the efficiency of the discretization, a \emph{modified} four-stage extension of WRK3,
\begin{equation}
    \varphi^{n,m} = \varphi^{n} + c_m\,\Delta t\,G\left(\varphi^{n,m-1}\right) \ , \qquad m = 1, \dots, 4 \ , \qquad
    c = \left(a,\ \frac13,\ \frac12,\ 1\right) \ ,
    \label{eq:mrk4-RK}
\end{equation}
named MRK4 from here on. With $a = 1/4$ its stability polynomial is the fourth-order Taylor polynomial, making the scheme fourth order for linear problems, while for nonlinear problems it remains second order, like WRK3. We retain a free parameter $a$, which controls damping when the scheme is coupled to the barotropic mode. With a slight abuse of notation, in the remainder of the manuscript, we drop the leading $n$ and write $\varphi^m$ for $\varphi^{n,m}$, retaining $\varphi^n$ and $\varphi^{n+1}$ for the solutions at the beginning and at the end of the timestep.

\label{sec:2-amplification}
\begin{figure}[H]
    \centering
    \includegraphics[width=\linewidth,height=0.78\textheight,keepaspectratio]{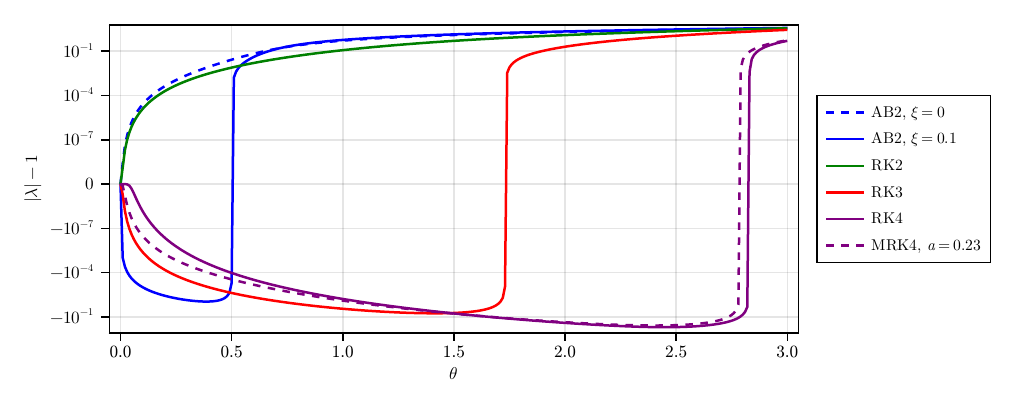}
    \caption{Departure of the amplification-factor magnitude from unity, $|\lambda| - 1$, against the phase change per step $\theta$, for AB2 and for the stability polynomials of Runge--Kutta schemes of two, three and four stages. WRK3 and SRK3 share the three-stage polynomial; RK4 denotes the stability polynomial of the four-stage Runge--Kutta for reference, of which MRK4 here is a slightly modified version with $a = 0.23$.}
    \label{fig:amplification-factor}
\end{figure}

We examine the numerical damping using the unit-speed advection equation, $\partial_t\varphi+\partial_x\varphi=0$. Each Fourier mode satisfies $\partial_t\hat\varphi=-ik\hat\varphi$, with exact amplification factor $e^{-i\theta}$, where $\theta=k\Delta t$ is the phase change over one timestep. A numerical amplification factor with magnitude below or above one indicates damping or growth while an error in the phase indicates dispersion.
An $M$-stage Runge--Kutta scheme whose order for linear problems is $M$ advances the mode as
\begin{equation}
    \hat{\varphi}^{n+1} = R_M(-i \theta)\hat{\varphi}^n,
    \qquad
    R_M(z) = \sum_{r=0}^{M}\frac{z^r}{r!} \ ,
\end{equation}
where the stability polynomial $R_M(z)$ is the degree-$M$ Taylor approximation of the exact amplification factor $e^z$.
WRK3 and SRK3 differ in nonlinear order and the SSP property; however, they share $R_3$ and therefore produce identical updates for simple linear systems.  
The stability polynomial of MRK4 is
\begin{equation}
    R(z) = 1 + z + \frac{z^2}{2} + \frac{z^3}{6} + \frac{a}{6}\,z^4 = R_3(z) + \frac{a}{6}\,z^4 \ ,
    \label{eq:mrk4-R}
\end{equation}
with an adjustable coefficient of $z^4$, which reduces to $R_3$ at $a = 0$ and to $R_4$ at $a = 1/4$.
The amplification factor of the Runge--Kutta schemes can be computed as $|\lambda_\text{RK}| = |R(z)|$.
For AB2, substituting $G(\hat{\varphi})=-ik\hat{\varphi}$ into the time-discrete equation gives
\begin{equation}
    \hat{\varphi}^{n+1}
    =
    \left[
        1-i \theta\left(\frac{3}{2}+\xi\right)
    \right]\hat{\varphi}^n
    +
    i \theta\left(\frac{1}{2}+\xi\right)
    \hat{\varphi}^{n-1}.
\end{equation}
Because AB2 is a two-step method, we seek solutions of the form
$\hat{\varphi}^n = \lambda \hat{\varphi}^{n-1} = \ldots = \lambda^n\hat{\varphi}^0$. Substitution yields the
quadratic characteristic equation
\begin{equation}
    \lambda^2
    -
    \left[
        1-i \theta \left(\frac{3}{2}+\xi\right)
    \right]\lambda
    -
    i \theta \left(\frac{1}{2}+\xi\right)
    =0.
\end{equation}
The physical root is the one that tends to one (neither damping nor growth) as $\theta$ tends to zero. The second root is the parasitic mode associated with the two-step formulation and approaches zero as $\theta$ tends to zero. The magnitude of the physical root's amplification factor is
\begin{align}
    \left|\lambda_{\mathrm{AB2}}\right|
    =
    \frac{1}{2}
    \left|
        1
        -
        \left(\frac{3}{2}+\xi\right)i\theta
        +
        \sqrt{
            1
            -
            i(1-2\xi)\theta
            -
            \left(\frac{3}{2}+\xi\right)^2
            \theta^2
        }
    \right|.
    \label{eq:amplificationAB2}
\end{align}
Figure~\ref{fig:amplification-factor} compares the magnitude of the amplification factors for different values of $\xi$ and $M$. Both AB2 with $\xi = 0$ and RK2 are linearly unstable for purely oscillatory systems ($|\lambda| > 1$) and require dissipation (either explicit or implicit through the spatial discretization) for stability. The three- and four-stage RK discretizations have wider stability intervals than modified AB2, and, more importantly, for small $\theta$ values, have much lower implicit damping. This weak damping at small $\theta$ is particularly relevant in ocean models, where the timestep is constrained by the most restrictive grid cells, while much of the resolved flow evolves at smaller local Courant numbers. For AB2, reducing $\xi$ lowers the damping, which scales as $\xi\theta^2$ at small $\theta$, but also the imaginary-axis limit, which is $0.2$, $0.4$ and $0.5$ for $\xi = 0.01, 0.05, 0.1$, respectively. Lowering $a$ from $1/4$ to $0.23$ slightly reduces the imaginary-axis limit of MRK4 and increases its damping. Table~\ref{tab:outer-schemes} collects the order and the imaginary-axis limit of the four schemes as well as the memory requirements in terms of three-dimensional fields per variable.

\begin{table}[t]
    \centering
    \begin{tabular}{lcccc}
        \hline
         & AB2, $\xi = 0.1$ & WRK3 & SRK3 & MRK4, $a = 0.23$ \\
        \hline
        tendency evaluations per step & $1$ & $3$ & $3$ & $4$ \\
        order, linear problems        & $1$ & $3$ & $3$ & $3$ \\
        order, nonlinear problems     & $1$ & $2$ & $3$ & $2$ \\
        memory requirements (3D fields) & $3$ & $3$ & $3$ & $3$ \\
        imaginary-axis limit $\theta_{\max}$ & $0.50$ & $\sqrt{3}$ & $\sqrt{3}$ & $2.79$ \\
        \hline
    \end{tabular}
    \caption{Properties of the \textit{uncoupled} outer time discretizations. AB2 is second order at $\xi = 0$, and MRK4 is fourth order for linear
    problems at $a = 1/4$, where $\theta_{\max} = 2\sqrt{2}$.}
    \label{tab:outer-schemes}
\end{table}

\section{Time discretization of the primitive equations}
\label{sec:3}

We apply the three schemes described in section~\ref{sec:2} to the primitive equations, taking care of the coupling between the three-dimensional evolution and the barotropic variables.  Asterisks denote predictor quantities; daggers denote transport quantities corrected to maintain tracer consistency.
We use a $z$-star vertical coordinate formulation \cite{zstar}, with the vertical coordinate $z$ given by
\begin{equation}
    z = r \frac{H+\eta}{H} + \eta \ .
    \label{eq:zstar}
\end{equation}
Here, $H$ is the water column height at rest, $\eta$ is the free surface elevation as measured from the rest state of the ocean ($z = 0$), and $r$ is a reference (fixed) geopotential vertical coordinate with $r = 0$ at $z = \eta$ and $r = -H$ at $z = -H$. The reference vertical coordinate is denoted $z^*$ by \citeA{zstar};
we follow the notation used in the Oceananigans source code, denoting the reference vertical coordinate by $r$ and the specific thickness by $\sigma=\partial_r z$, to facilitate comparison with the model documentation and implementation.
The specific thickness (or dimensionless vertical stretching) is
\begin{equation}
    \sigma = \pd{z}{r} = \frac{H + \eta}{H} \ .
\end{equation}

The primitive equations in $r$ coordinates \cite{adcroft2006}, \citeA<derived in>{oceananigans-docs-gvc}, read
\begin{align}
    \pd{\bu}{t} &= \bG_u - g \bnabla \eta \ , \label{eq:momentumeq}\\ 
    \pd{p}{r} &= \sigma b \ , \\
    \pd{\sigma}{t} &= - \bnabla \cdot (\sigma \bu) - \pd{\omega}{r} \ ,  \label{eq:continuity}\\
    \pd{\sigma C}{t} &= - \bnabla \cdot (\sigma \bu C) - \pd{(\omega C)}{r} \ , \label{eq:tracereq}
\end{align}
where $g$ is the gravitational acceleration, $\bnabla$ is the horizontal gradient operator ($\hat{\bx}\,\partial_x + \hat{\by} \,\partial_y$) at constant $r$, $\bu$ is the horizontal velocity field $u \hat{\bx} + v \hat{\by}$, $\omega$ is the vertical velocity across $r$ surfaces, $C$ indicates a (passive or active) tracer, $p$ is the baroclinic hydrostatic pressure divided by the reference density, and $b$ is the buoyancy (which depends on the temperature $T$, salinity $S$, and depth $z$).
Finally, $\bG_u$ includes all the ``slow'' baroclinic tendency terms (advection, Coriolis, hydrostatic pressure gradient, vertical mixing, and boundary conditions).
The above equations are complemented by the evolution equation for the free surface, which, in the absence of precipitation or evaporation, reads
\begin{equation}
    \pd{\eta}{t} = - \bnabla \cdot \underbrace{\int_{-H}^{0} \sigma \bu \, \mathrm{d}r}_{\defn \bU}  \ ,
\end{equation}
where $\bU$ is the barotropic transport. In this manuscript we do not consider surface volume fluxes (precipitation, evaporation and runoff) or open boundaries. The free surface must also be consistent with equation~\eqref{eq:zstar}, which implies,
\begin{equation}
    \eta = \int_{-H}^0 \sigma \, \mathrm{d}r - H \ , 
    \label{eq:free-surface-consistency}
\end{equation}
where the dynamic column thickness is $H + \eta$.

We require our time discretization to satisfy three properties: (1) consistency between free surface and column thickness, (2) tracer conservation and constancy preservation, and (3) stability of the coupled barotropic-baroclinic system at large timesteps.

The first is satisfied by construction. The grid is not evolved explicitly, as in models that use more sophisticated vertical coordinates \cite<e.g.,>{adcroft2019}, but is diagnosed from equation~\eqref{eq:zstar}, so that~\eqref{eq:free-surface-consistency} always holds and the continuity equation~\eqref{eq:continuity} is used only to diagnose the vertical velocity $\omega$.
The $\partial_t \sigma$ term in the continuity equation follows from the solve for the free surface $\eta$.

The second requires compatibility between the discrete continuity equation and the tracer fluxes to ensure both global tracer conservation and constancy preservation \cite{shchepetkin}. In the finite-volume discretization of Oceananigans, the tracer flux leaving one cell enters its neighbor, so internal fluxes cancel and the total tracer content should change only through boundary fluxes. Additionally, since the evolution equation for a uniform tracer reduces to the continuity equation~\eqref{eq:continuity}, constancy preservation (a uniform tracer remains uniform) requires the time discretization of the grid evolution to match that of the tracer equation. Since we use a $z$-star formulation, the grid evolves with the free surface,
\begin{equation}
    \pd{\sigma}{t} = \frac{1}{H} \pd{\eta}{t} = - \frac{1}{H} \bnabla \cdot \bU \ ,
\end{equation}
so an explicit free-surface update using the same transport discretization as the tracers would preserve a uniform tracer by construction. Such a pairing, however, restricts the timestep through the fast surface gravity waves, making it prohibitively expensive for practical hydrostatic ocean models.

Two approaches relax the timestep restriction imposed by surface gravity waves: an implicit (or semi-implicit) discretization \cite{dukowicz1994, kang2021}, or a split--explicit sub-cycling approach, where the barotropic mode is advanced with a much smaller timestep within each baroclinic timestep \cite{killworth1991, shchepetkin}. We implement both, coupling the free surface and the
three-dimensional evolution at each Runge--Kutta stage. With the implicit solver this corresponds to an implicit--explicit (IMEX) RK scheme
with the same coefficients for the explicit and the implicit components, which is first order in the implicit component~\cite{kar}. Since the implicit scheme is used here as a benchmark, IMEX schemes of higher implicit order~\cite{conde, atmosphereSSP}, as well as semi-implicit formulations~\cite{kevlahan2022} are not considered.

\subsection{The RK stage design}
\label{sec:3-stage}

The schemes share four ingredients: stage fractions, slow forcing, a stage-composition rule, and a barotropic solver. WRK3 and MRK4 restart each stage from the beginning-of-step state; SRK3 blends that state with an update from the preceding stage. We refer to these as \textbf{restarted} and \textbf{blended} stages, respectively.
To satisfy grid and free surface consistency, the grid follows the free surface at every stage,
\begin{equation}
    \sigma^{m} = \frac{H + \eta^{m}}{H} \ .
    \label{eq:sigma-stage}
\end{equation}
This implies that, since $\eta^m$ is advanced with a different time-discretization than the tracers, the tracer transport must be corrected to ensure constancy preservation.
Each barotropic solve spans $\gamma_m$, with
\begin{equation}
    \gamma_m \defn c_m\,\Delta t \ , \qquad
    c = \begin{cases}
        \left(\dfrac13,\ \dfrac12,\ 1\right) & \text{WRK3} \ , \\[4pt]
        \left(a,\ \dfrac13,\ \dfrac12,\ 1\right) & \text{MRK4} \ , \\[4pt]
        \left(1,\ 1,\ 1\right) & \text{SRK3} \ .
    \end{cases}
    \label{eq:defgamma}
\end{equation}
We call $\mathcal{B}$ the \emph{barotropic operator} that advances the increment of the free surface $\Delta \eta^m$ and the barotropic transport $\Delta\bU^m$ from the beginning of the step over the stage timestep $\gamma_m$, driven by a depth-integrated slow forcing $\widehat{\bG}_u^{\,m}$.
Besides the two increments, $\mathcal{B}$ returns the thickness-weighted transport $\bU^\dagger$ consistent with the free-surface increment,
\begin{equation}
    \Delta\eta^{m},\ \Delta\bU^{m},\ \bU^\dagger = \mathcal{B}\left(\eta^n,\ \bU^n,\ \widehat{\bG}_u^{\,m},\ \gamma_m\right) \ , \qquad
    \Delta\eta^{m} = - \gamma_m\,\bnabla \cdot \bU^\dagger \ .
    \label{eq:transport-identity}
\end{equation}
The specific computation of $\bU^\dagger$ depends on the barotropic operator $\mathcal B$. We construct corrected velocities $\bu^\dagger$ (with thickness-weighted vertical integral equal to $\bU^\dagger$), and $\omega^\dagger$ (from continuity) to advect tracers such that, for a uniform tracer, the tracer update reduces to the thickness update implied by $\Delta\eta^m$ and the scheme preserves constant tracers.
$\widehat{\bG}_u^{\,m}$ is the slow forcing collecting the depth-integrated ``slow'' tendencies $\bG_u$; the default choice is to simply use the value computed from the previous stage
\begin{equation}
    \widehat{\bG}_u^{\,m} = \overline{\bG}{}^{\,m-1} \ , \qquad
    \overline{\bG}{}^{\,j} \defn \int_{-H}^{0} \sigma^{j}\, \bG_u^{j} \, \mathrm{d}r \ ,
    \label{eq:Ghat-frozen}
\end{equation}
and keep it fixed during the barotropic solve. We refer to this implementation as \emph{frozen} forcing.

\subsubsection{The restarted stage (WRK3 and MRK4)}
The stage of the restarted formulations reads
\begin{align}
    \Delta \eta^{m}, \Delta \bU^{m}, \bU^\dagger &= \mathcal{B}\left(\eta^n, \bU^n, \widehat{\bG}_u^{\,m}, \gamma_m\right) \ , \\
    \eta^{m},\ \bU^{m} &= \eta^{n} + \Delta\eta^{m},\ \bU^{n} + \Delta\bU^{m} \ , \\
    \bu^\dagger &= \bu^{m-1} - \frac{1}{H + \eta^{m-1}}\left(\int_{-H}^0 \sigma^{m-1} \bu^{m-1} \, \mathrm{d}r - \bU^\dagger\right) \label{eq:transport-evolution}  \ , \\
    \omega^\dagger &= - \int_{-H}^r \left[\bnabla \cdot (\sigma^{m-1} \bu^{\dagger}) - \frac{1}{H}\bnabla \cdot \bU^{\dagger} \right] \, \mathrm{d}r^\prime \ , \\
    \sigma^{m} &= \frac{H + \eta^{m}}{H} \label{eq:sigma-m+1} \ , \\
    C^{m} &= \frac{\sigma^n}{\sigma^{m}} C^n - \frac{\gamma_m}{\sigma^{m}}\left[\bnabla \cdot (\sigma^{m-1} \bu^\dagger C^{m-1})  + \pd{(\omega^\dagger C^{m-1})}{r} \right] \ , \\
    \bu^* &= \bu^n + \gamma_m \bG_u^{m-1}  \label{eq:stage-velocity} \ , \\ 
    \bu^{m} &= \bu^* - \frac{1}{H + \eta^{m}}\left(\int_{-H}^0 \sigma^{m} \bu^* \, \mathrm{d}r  - \bU^{m}\right) \label{eq:stage-reconcile} \ , \\
    \omega^{m} &= - \int_{-H}^r \left[ \bnabla \cdot (\sigma^{m} \bu^{m}) - \frac{1}{H}\bnabla \cdot \bU^{m} \right] \, \mathrm{d}r^\prime \ .
\end{align}
The velocity advance in $\eqref{eq:stage-velocity}$ is complemented by the free-surface reconciliation in~\eqref{eq:stage-reconcile}, which, for the implicit free surface, takes the form of a surface-pressure correction. In a closed lateral domain, in the absence of sources and sinks, the global tracer conservation is ensured by the flux-form finite volume discretization together with the fact that at the surface
\begin{equation}
\omega^\dagger(0) = -\bnabla\cdot \underbrace{\int_{-H}^{0}\sigma^{m-1}\bu^\dagger\,\mathrm{d}r}_{\bU^\dagger} +\bnabla\cdot\bU^\dagger =0 \ .
\end{equation}
Constancy preservation also follows by construction. Assuming the tracer is initially uniform ($C^n=C^{m-1}=C_0$), from the definition of $\omega^\dagger$,
\begin{equation}
\bnabla\cdot\left(\sigma^{m-1}\bu^\dagger\right)
+\frac{\partial\omega^\dagger}{\partial r}
=\frac{1}{H}\bnabla\cdot\bU^\dagger \ ,    
\end{equation}
substituting the constant tracer into the stage update gives
\begin{equation}
\sigma^m C^m
=C_0\left(\sigma^n-\frac{\gamma_m}{H}
\bnabla\cdot\bU^\dagger\right).    
\end{equation}
But the free-surface update implies
\begin{equation}
\sigma^m
=\sigma^n+\frac{\Delta\eta^m}{H}
=\sigma^n-\frac{\gamma_m}{H}\bnabla\cdot\bU^\dagger \ .
\end{equation}
Thus $\sigma^m C^m=C_0\sigma^m$ giving $C^m=C_0$.

\subsubsection{The blended stage (SRK3)}
\label{sec:3-ssp}
For SRK3, we must specify how the Shu--Osher blend acts on the barotropic variables. In this case, each barotropic solve spans the full $\Delta t$, although the blended states correspond to a nominal time of ($t^n + \Delta t$, $t^n + \Delta t / 2$, $t^n + \Delta t$).
For this reason, applying the blend to the full barotropic state output of the barotropic solve would produce an over-integration of $\eta$ and $\bU$ with respect to the slow variables. Therefore, we restart each barotropic solve from $(\eta^n, \bU^n)$ exactly as in the restarted stage, and we blend the \emph{increment} each solve returns rather than the state itself. The slow forcing that drives $\mathcal{B}$ is the frozen forcing of~\eqref{eq:Ghat-frozen}.
Accounting for the Shu--Osher blend, the stage reads
\begin{align}
    \Delta\eta^{m},\ \Delta\bU^{m},\ \bU^\dagger & = \mathcal{B}\left(\eta^{n}, \bU^{n}, \widehat{\bG}_u^{\,m}, \gamma_m\right) \ , \\ 
    \eta^{m} &= \left(1 - b_m\right)\eta^{n} + b_m\left(\eta^{m-1} + \Delta\eta^{m}\right) \ ,
      \label{eq:ssp-etablend} \\
    \bU^{m} &= \left(1 - b_m\right)\bU^{n} + b_m\left(\bU^{m-1} + \Delta\bU^{m}\right) \ ,
      \label{eq:ssp-Ublend} \\
    \sigma^{m} &= \frac{H + \eta^{m}}{H} \ , \label{eq:ssp-sigma} \\
    \bu^\dagger &= \bu^{m-1} - \frac{1}{H +\eta^{m-1}}
        \left(\int_{-H}^0 \sigma^{m-1}\bu^{m-1}\,\mathrm{d}r - \bU^\dagger\right) \ , \label{eq:ssp-udag} \\
    \omega^\dagger &= -\int_{-H}^{r}\left[\bnabla\cdot\left(\sigma^{m-1}\bu^\dagger\right)
        - \frac{1}{H}\bnabla\cdot\bU^\dagger\right]\mathrm{d}r^\prime \ , \label{eq:ssp-omegadag} \\
    \sigma^{m} C^{m} &= \left(1 - b_m\right)\sigma^{n} C^{n}
        + b_m\Big\{\sigma^{m-1} C^{m-1}
        - \gamma_m\Big[\bnabla\cdot\left(\sigma^{m-1}\bu^\dagger C^{m-1}\right)
        + \pd{\left(\omega^\dagger C^{m-1}\right)}{r}\Big]\Big\} \ , \label{eq:ssp-tracer} \\
    \bu^{\,\ast} &= \left(1 - b_m\right)\bu^{n} + b_m\left(\bu^{m-1} + \gamma_m\,\bG_u^{m-1}\right) \ , \label{eq:ssp-velocity} \\
    \bu^{m} &= \bu^{\,\ast} - \frac{1}{H+\eta^{m}}
        \left(\int_{-H}^0 \sigma^{m}\bu^{\,\ast}\,\mathrm{d}r - \bU^{m}\right) \ , \label{eq:ssp-reconcile} \\
    \omega^{m} &= -\int_{-H}^{r}\left[\bnabla\cdot\left(\sigma^{m}\bu^{m}\right)
        - \frac{1}{H}\bnabla\cdot\bU^{m}\right]\mathrm{d}r^\prime \ . \label{eq:ssp-omega}
\end{align}
As in the restarted stage, the barotropic reconciliation~\eqref{eq:ssp-reconcile} runs at every stage, so that the last stage returns $\bu^{n+1}$ and $\omega^{n+1}$, and global tracer conservation follows from the same flux cancellation as in the restarted stage. 
The increment blend weights the three increments by $\beta$, so that the depth-integrated slow forcing receives the same weights as the baroclinic tendencies. 
This coupling differs from the scheme of \citeA{lanrk3}, which instead drives the final barotropic solve with the weighted forcing $\overline G{}^0/6 +\overline G{}^1/6+2\overline G{}^2/3$.
Our formulation also accommodates an arbitrary barotropic operator $\mathcal B$, restarts each solve from the beginning-of-step state, and blends the barotropic increments rather than the state. 
The increment blend is required to ensure constancy preservation: setting $C^n = C^{m-1} = C_0$ in the tracer update yields 
\begin{equation}
        \sigma^{m} C^{m} = C_0 \Big(\left(1 - b_m\right)\sigma^{n} 
        + b_m\Big[\sigma^{m-1}
        + \frac{\Delta \eta^m}{H} \Big]\Big)\ , 
\end{equation}
where we have used the continuity relation~\eqref{eq:ssp-omegadag} and $\Delta \eta^m = -\gamma_m \bnabla\cdot \bU^\dagger$ to substitute the free surface increment. Substituting the definition of $\sigma$, multiplying by $H$ and dividing by $C_0$ throughout gives
\begin{equation}
    (H +\eta^{m}) \frac{C^m}{C_0} = \left(1 - b_m\right)(H +\eta^{n}) + b_m\left(H +\eta^{m-1} + \Delta \eta^m\right) \ ,
\end{equation}
so that $C^m = C_0$ if and only if the free surface follows the increment blend~\eqref{eq:ssp-etablend}. For this reason, equation~\eqref{eq:ssp-etablend} is the only free-surface update consistent with constancy preservation. Equation~\eqref{eq:ssp-Ublend} is then necessary to keep the barotropic variables synchronized.

This blended stage also preserves the strong-stability property for the tracers. We define the forward-Euler step of stage $m$, which advances the tracer content of the preceding stage with the corrected transport,
\begin{equation}
    \widetilde{\sigma}^{m}\widetilde{C}^{m} = \sigma^{m-1} C^{m-1}
        - \gamma_m\Big[\bnabla\cdot\left(\sigma^{m-1}\bu^\dagger C^{m-1}\right)
        + \pd{\left(\omega^\dagger C^{m-1}\right)}{r}\Big] \ , \qquad
    \widetilde{\sigma}^{m} = \sigma^{m-1} + \frac{\Delta\eta^m}{H} \ ,
    \label{eq:ssp-euler}
\end{equation}
where, by virtue of~\eqref{eq:ssp-omegadag} and of the transport identity~\eqref{eq:transport-identity}, $\widetilde{\sigma}^{m}$ is the thickness consistent with the transport of the step. Substituting~\eqref{eq:ssp-euler} in the tracer update~\eqref{eq:ssp-tracer}, and the free-surface blend~\eqref{eq:ssp-etablend} in~\eqref{eq:ssp-sigma}, gives
\begin{equation}
    \sigma^{m} C^{m} = \left(1 - b_m\right)\sigma^{n} C^{n} + b_m\,\widetilde{\sigma}^{m}\widetilde{C}^{m} \ , \qquad
    \sigma^{m} = \left(1 - b_m\right)\sigma^{n} + b_m\,\widetilde{\sigma}^{m} \ ,
\end{equation}
so that the stage concentration is a pointwise convex combination of the initial concentration and of the forward-Euler update,
\begin{equation}
    C^{m} = \vartheta_m\, C^{n} + \left(1 - \vartheta_m\right)\widetilde{C}^{m} \ , \qquad
    \vartheta_m \defn \frac{\left(1 - b_m\right)\sigma^{n}}{\sigma^{m}} \ ,
    \label{eq:ssp-convex}
\end{equation}
where $\vartheta_m$ lies between zero and one since $\sigma^{n}$ and $\widetilde{\sigma}^{m}$ are positive quantities. As a consequence, if the forward-Euler step~\eqref{eq:ssp-euler} keeps $\widetilde{C}^{m}$ within the bounds of $C^{m-1}$ for timesteps below a limit $\Delta t_{\mathrm{FE}}$, as for a monotone or positivity-preserving advection scheme, the SSP property ensures that $C^{n+1}$ remains within the bounds of $C^{n}$ for the same timestep. This property does not extend to the restarted stage, where the tendency of $C^{m-1}$ is added to $C^{n}$ rather than to $C^{m-1}$.

\subsection{The implicit solver}
\label{sec:3-barotropic}
The implicit free-surface solver provides unconditional stability for the evolution of the barotropic mode. It advances the barotropic variables with one elliptic solve per stage, extending the approach of~\citeA{dukowicz1994} and~\citeA{campin2004} to each Runge--Kutta stage.
We use it here as a benchmark, rather than as a production option, since efficient elliptic solves require suitable preconditioning, whose parallel implementation can be challenging on the GPUs targeted by our model.
Also, for simplicity, we apply it only to the restarted stages (WRK3 and MRK4).
The free surface is evolved with the updated velocity at stage $m$. Theoretically, using the evolved $\eta^m$ everywhere on the right-hand side, the free surface obeys the fully implicit nonlinear evolution given by
\begin{equation}
    \frac{\eta^{m} - \eta^{n}}{\gamma_m} = -\bnabla\cdot\left[\int_{-H}^0 \sigma^{m}\left( \bu^n + \gamma_m \bG_u^{m-1} \right) dr - \gamma_m\, g\, (H + \eta^{m})\bnabla\eta^{m}\right] \ .
\end{equation}
However, the above equation needs an expensive nonlinear solve; therefore, the typical approach is to linearize the above using information at the previous stage for the geometry terms,
\begin{equation}
    \frac{\eta^{m} - \eta^{n}}{\gamma_m} = -\bnabla\cdot\left[\int_{-H}^0 \sigma^{m-1}\left( \bu^n + \gamma_m \bG_u^{m-1} \right) dr - \gamma_m\, g\, (H + \eta^{m-1})\bnabla\eta^{m}\right] \ ,
    \label{eq:im-evolution}
\end{equation}
which rearranges into the linear system
\begin{equation}
    \eta^{m} - \gamma_m^2 g \, \bnabla \cdot \left[ (H+\eta^{m-1}) \bnabla \eta^{m} \right]
    = \eta^n - \gamma_m\bnabla \cdot \int_{-H}^0 \sigma^{m-1}\left( \bu^n + \gamma_m \bG_u^{m-1} \right) dr \ .
    \label{eq:im-stage-solve}
\end{equation}
We solve this system with a standard conjugate-gradient method and a relative tolerance of $10^{-7}$. Finally, the velocity is corrected by the surface-pressure gradient, which replaces the reconciliation~\eqref{eq:stage-reconcile},
\begin{equation}
    \begin{gathered}
    \bu^{m} = \bu^n + \gamma_m \bG_u^{m-1} - \gamma_m\, g\,\bnabla\eta^{m} \ , \qquad
    \bU^{m} \defn \int_{-H}^0 \sigma^{m} \bu^{m} \, \mathrm{d}r \ ,
    \end{gathered}
    \label{eq:im-transport}
\end{equation}
and, consistently with the linearization of~\eqref{eq:im-stage-solve}, we compute $\bU^\dagger$ using the evolved velocity $\bu^m$ integrated with geometry terms of the previous stage $\sigma^{m-1}$, 
\begin{equation}
    \bU^\dagger \defn \int_{-H}^0 \sigma^{m-1}\bu^m dr
\end{equation}
so that the transport identity~\eqref{eq:transport-identity} holds and constancy preservation is ensured by construction (up to the tolerance of the solver).

\subsection{The split--explicit solver}
The split--explicit barotropic solver resolves the fast dynamics explicitly, with $N_S$ substeps per stage of size
$\gamma_B$, a fraction $\Delta\tau$ of the stage timestep, 
\begin{equation}
    \gamma_B \defn \Delta\tau\,\gamma_m \ , \qquad \Delta\tau \simeq \frac{2}{N_S} \ ,
    \label{eq:defgammaB}
\end{equation}
where $\Delta\tau$ differs from $2/N_S$ by the rescaling that enforces the centering of the averaging filter
(Appendix~\ref{sec:appendix-filters}). The factor two extends the barotropic integration beyond the stage endpoint, allowing the time filter to be always centered at $t^n+\gamma_m$. 
We investigate two barotropic discretizations. The first is a forward-backward (FB) discretization, implemented as:
\begin{alignat}{3}
    \eta^{0}, \ \bU^{0} &= \eta^n, \ \bU^n \ , && \text{at the beginning of each barotropic solve, } \\
    \bU^{s+1} &= \bU^{s} + \gamma_B \left(\widehat{\bG}_u^{\,m} - g (H + \eta^{s}) \bnabla \eta^{s}\right) \ \ && \text{for} \ \ s = 0,\dots, N_S - 1 \ , \\
    \eta^{s+1} &= \eta^{s} - \gamma_B \bnabla \cdot \bU^{s+1} \ \ && \text{for} \ \ s = 0,\dots, N_S - 1 \ ,
\end{alignat}
where $s$ indicates the sub-cycling index. The second is the same WRK3 applied to the barotropic mode, which, with an inner stage index $j$ restarting from the beginning of each substep, reads
\begin{alignat}{3}
    \eta^{0}, \ \bU^{0} &= \eta^n, \ \bU^n \ ,  && \text{at the beginning of each barotropic solve, } \\
    \eta^{s,0}, \ \bU^{s,0} &= \eta^{s}, \ \bU^{s} \ , && \text{at the beginning of each substep, } \\
    \bU^{s,j} &= \bU^{s} + c_j \gamma_B \left(\widehat{\bG}_u^{\,m} - g (H + \eta^{s,j-1}) \bnabla \eta^{s,j-1}\right) \ \ && \text{for} \ \ j = 1, 2, 3 \ , \\
    \eta^{s,j} &= \eta^{s} - c_j \gamma_B \bnabla \cdot \bU^{s,j-1} \ \ && \text{for} \ \ j = 1, 2, 3 \ , \\
    \eta^{s+1}, \ \bU^{s+1} &= \eta^{s,3}, \ \bU^{s,3} \ \ && \text{for} \ \ s = 0,\dots, N_S - 1 \ ,
\end{alignat}
where $c = \left(1/3, 1/2, 1\right)$ are the WRK3 fractions of~\eqref{eq:ws-RK}, now applied to the substep rather than to the baroclinic step. 

For both barotropic time discretizations, we follow the filtering approach of~\citeA{shchepetkin}, where the sub-cycled sequence is averaged with a time filter with output weights $w_s$ and transport weights $w_s^\dagger$, such that the output of the barotropic solve $\mathcal{B}$ (in terms of increments) is
\begin{equation}
    \Delta \eta^{m} = \sum_{s=1}^{N_S} w_s\,\eta^{s}  - \eta^n\ , \qquad \text{and} \qquad
    \Delta \bU^{m} = \sum_{s=1}^{N_S} w_s\,\bU^{s} - \bU^n\ ,
    \label{eq:se-filtered-outputs}
\end{equation}
while the corrected transport accumulates, with $w^\dagger_s$, the transport that advances the free surface from $\eta^{s-1}$ to $\eta^{s}$, which is $\bU^{s}$ for FB and $\bU^{s-1,2}$ for RK3,
\begin{equation}
    \bU^\dagger = \sum_{s=1}^{N_S} w^\dagger_s\,\begin{cases} \bU^s & \text{for FB} \\ \bU^{s-1, 2}  &  \text{for RK3}\end{cases}\ , \qquad \text{where} \qquad
    w^\dagger_s = \Delta\tau \sum_{s^\prime=s}^{N_S} w_{s^\prime} \ .
    \label{eq:se-transport}
\end{equation}
The relation between (secondary) transport weights $w_s^\dagger$ and the (primary) output weights $w_s$ follows from summing the continuity updates of the substeps, so that $(\Delta \eta^{m}, \bU^\dagger)$ satisfies the transport identity~\eqref{eq:transport-identity} exactly. A detailed derivation can be found in \citeA{shchepetkin}.
The filtering weights are constrained by consistency and by the centering on the stage timestep $\gamma_m$,
\begin{equation}
    \sum_{s} w_s = 1 \ , \qquad \sum_{s} w_s\,\tau_s = 1 \ , \qquad \tau_s \defn \frac{s\,\gamma_B}{\gamma_m} \ ,
    \label{eq:filter-moments}
\end{equation}
where $\tau_s$ is the fractional sub-cycle time. For a fixed $N_S$ and filter shape, the same filter weights are used for all the outer stages and for every outer time-discretization (WRK3, SRK3 or MRK4). The higher moments of the filter
\begin{equation}
    M_p \defn \sum_{s} w_s\,(\tau_s - 1)^p
    \label{eq:filter-highmoments}
\end{equation}
set the order and the dissipation of the filtering procedure.

\subsection{Damping of barotropic waves}
\label{sec:3-barotropic-damping}
\begin{figure}[H]
    \centering
    \includegraphics[width=\linewidth,height=0.78\textheight,keepaspectratio]{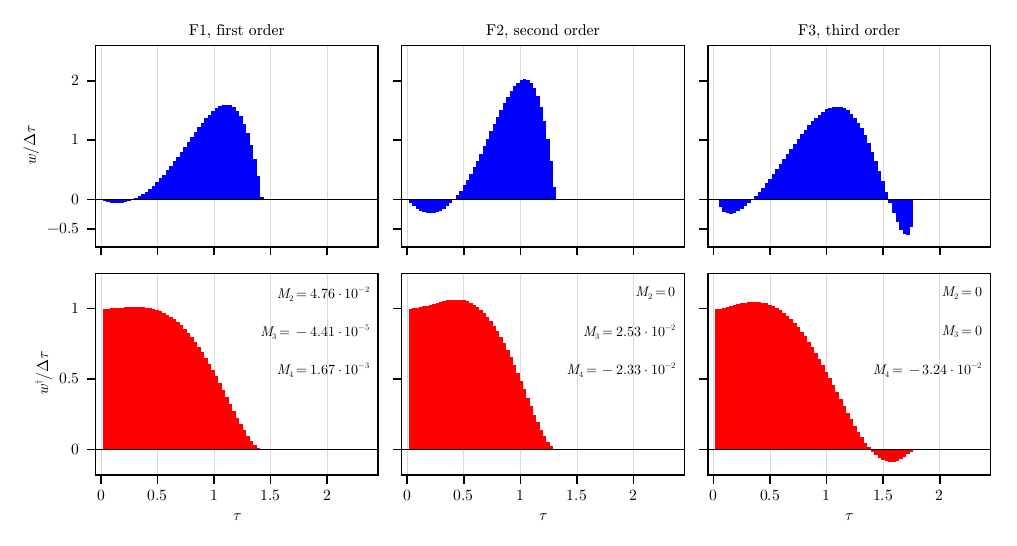}
    \caption{The first- and second-order averaging filters of~\protect\citeA{shchepetkin}, F1 and F2, and the third-order filter F3, for $64$ substeps: averaging weights $w$ (top) and transport weights $w^\dagger$ (bottom), divided by the spacing $\Delta\tau$ of the fractional sub-cycle times, with the moments $M_2$, $M_3$ and $M_4$.}
    \label{fig:averaging-weights}
\end{figure}
Following the amplification factor analysis of section~\ref{sec:2}, we compare the damping introduced by the barotropic solvers and their time filters. We consider a horizontal Fourier mode of wavenumber $k$ of the unforced barotropic equations linearized about rest at constant depth $H$ ($\sigma = 1$). We use $\bar{u}$ and $\eta$ to denote the depth-averaged velocity and free surface normalized by $\sqrt{gH}$ and $H$, respectively. Then the barotropic operator $\mathcal{B}$ is represented by the amplification matrix $\mathbf{W}$ over the stage interval $\gamma_m$, 
\begin{equation}
    \begin{bmatrix}\bar u^m \\ \eta^m\end{bmatrix}= \mathbf{W}\begin{bmatrix}\bar u^n \\ \eta^n\end{bmatrix} \ .
\end{equation}
The exact evolution is a pure rotation with
\begin{equation}
    \mathbf{W} = e^{-i\mu_m\mathbf{J}} = \cos\mu_m\,\mathbf{I} - i\sin\mu_m\,\mathbf{J} \ , \qquad \text{where} \qquad 
    \mu_m \defn \sqrt{gH}\,k\,\gamma_m \ , \qquad \text{and} \qquad
    \mathbf{J} \defn \begin{bmatrix} 0 & 1 \\ 1 & 0 \end{bmatrix} \ ,
    \label{eq:windows}
\end{equation}
with $\mu_m$ the barotropic Courant number of the stage. The two eigenvalues of $\mathbf{W}$ are the amplification factors for the two wave propagation directions, so their magnitudes determine damping or growth (as in section~\ref{sec:2}); we call one of them $\lambda_b$. The implicit updates of~\eqref{eq:im-evolution} and~\eqref{eq:im-transport} reduce to backward Euler for the barotropic system, giving the simple amplification matrix
\begin{equation}
    \mathbf{W} = \big(\mathbf{I} + i\mu_m\,\mathbf{J}\big)^{-1} \ ,
    \label{eq:W-im}
\end{equation}
whose eigenvalues have magnitude equal to
\begin{equation}
    |\lambda_b| = 1 / \sqrt{1+ \mu_m^2} < 1 \ , 
\end{equation} 
so that the implicit solve damps barotropic waves more strongly as $\mu_m$ increases. The filtered barotropic solve gives the amplification matrix
\begin{equation}
    \mathbf{W} = \sum_{s} w_s\,\mathbf{W}_{B}^{\,s} \ ,
    \label{eq:W-se}
\end{equation}
\begin{figure}[H]
    \centering
    \includegraphics[width=\linewidth,height=0.78\textheight,keepaspectratio]{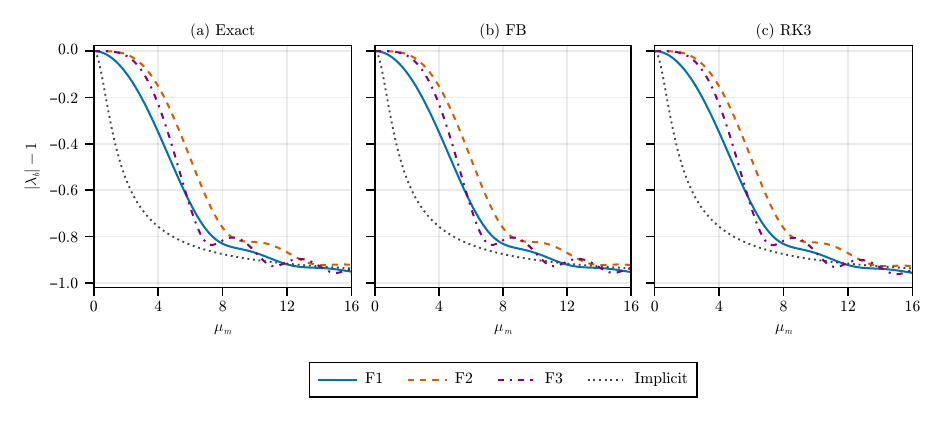}
    \setlength{\abovecaptionskip}{3pt}
    \caption{Departure of the barotropic amplification-factor magnitude from unity, $|\lambda_b|-1$, against the stage Courant number $\mu_m$, for filters F1, F2 and F3 with (a) exact, (b) forward--backward (FB), and (c) RK3 substeps using $64$ substeps. The implicit solver is shown in grey.}
    \label{fig:barotropic-filter-damping}
\end{figure}
with $\mathbf{W}_B$ the evolution matrix for one substep. To isolate the filter response, we first take each substep to be exact. As in \citeA{shchepetkin}, the amplification factor for one direction can be computed analytically with
\begin{equation}
    \begin{aligned}
    \lambda_b = \sum_{s} w_s\,e^{-i\mu_m\tau_s}
    &= e^{-i\mu_m}\sum_{p}\frac{(-i\mu_m)^{p}}{p!}\sum_{s} w_s\,(\tau_s - 1)^{p} \\
    &= e^{-i\mu_m}\Big[1 - \frac{M_2}{2}\mu_m^{2} + i\,\frac{M_3}{6}\mu_m^{3} + \frac{M_4}{24}\mu_m^{4} + \mathcal{O}(\mu_m^{5})\Big] \ ,
    \end{aligned}
    \label{eq:filter-transfer}
\end{equation}
where the normalization and the centering of~\eqref{eq:filter-moments} set $M_0 = 1$ and $M_1 = 0$. A filter is of order $p$ if $M_1,\ldots,M_p$ vanish.
The magnitude of this amplification factor is
\begin{equation}
    \big|\lambda_b| = 1 - \frac{M_2}{2}\mu_m^{2} + \frac{M_4}{24}\mu_m^{4} + \mathcal{O}(\mu_m^{6}) \ .
    \label{eq:filter-gain}
\end{equation}
The leading damping is second order in $\mu_m$ when $M_2\neq0$, and fourth order in $\mu_m$ (for negative $M_4$) otherwise.
Increasing the number of substeps does not eliminate this dissipation, because the moments are determined by the filter shape.

Filters of different order are shown in Figure~\ref{fig:averaging-weights}: F1 (first order) and F2 (second order) belong to the low-dispersion family of~\citeA{shchepetkin}, F1 being the default in Oceananigans and F2 the member of the same family with $M_2 = 0$, while F3, used in most of the experiments of section~\ref{sec:5}, is constructed from a polynomial multiplied by a taper, satisfies $M_2=M_3=0$, and is designed to be more dissipative than F2 at higher Courant numbers, while dissipating slow waves less than F1. Their definitions and their discretization are given in Appendix~\ref{sec:appendix-filters}.
Figure~\ref{fig:barotropic-filter-damping} compares the resulting damping with that of the implicit solver. With $64$ nominal substeps, the exact, FB and RK3 results are nearly identical over the range shown, indicating that the damping is dominated by the time filter.

\subsection{The slow forcing reconstruction}
\label{sec:3-mrk4}
The implicit solver and SRK3 use a frozen forcing. This is similar to the approach of \citeA{wensch2009} and \citeA{knoth2014} where the slow forcing is updated between stages and held constant during fast integration. For split--explicit WRK3 and MRK4, the tendencies computed at the previous stages allow a \emph{progressive} reconstruction of the vertically integrated forcing $\widehat{\bG}_u^{\,m}$ during the barotropic substeps. Here we take inspiration from the MRI-GARK formulation of \citeA{sandu2019}, in which the weights of the slow tendencies can vary polynomially in time during fast integration.
Denoting by $t'$ the time elapsed since $t^n$, we write the reconstructed forcing as
\begin{equation}
    \widehat{\bG}_u^{\,m}(t')
    = \sum_{p=0} \bG_p^m\,t'^p,
    \label{eq:progressive-reconstruction}
\end{equation}
where $\bG_p^m$ denotes the coefficient of $t'^p$ at stage $m$. The three-dimensional tendencies used to compute the coefficients lie at $t'=c_j\Delta t$, with $(c_0,c_1,c_2)=(0,1/3,1/2)$ for WRK3 and $(c_0,c_1,c_2,c_3)=(0,a,1/3,1/2)$ for MRK4. 
At the beginning of the timestep, only $\overline{\mathbf G}{}^{0}$ is available. The forcing used during the first stage is therefore constant:
\begin{equation}
\widehat{\mathbf G}_u^{\,1}(t')
=\mathbf G_0^1
=\overline{\mathbf G}{}^{\,0}.
\end{equation}
After the first stage, we also have $\overline{\mathbf G}{}^{\,1}$, associated with $t'=c_1\Delta t$. We can therefore construct a linear forcing for the second stage;
requiring that $\widehat{\mathbf G}_u^{\,2}(0)=\overline{\mathbf G}{}^{\,0}$ and $\widehat{\mathbf G}_u^{\,2}(c_1\Delta t)=\overline{\mathbf G}{}^{\,1}$ gives
\begin{equation}
\mathbf G_0^2=\overline{\mathbf G}{}^{\,0},
\qquad
\mathbf G_1^2=
\frac{\overline{\mathbf G}{}^{\,1}-\overline{\mathbf G}{}^{\,0}}
     {c_1\Delta t} \ .
\end{equation}
For the third stage, a quadratic reconstruction that adds information from
$\overline{\mathbf G}{}^{\,2}$, computed at $t'=c_2\Delta t$, retains
$\mathbf G_0^3=\mathbf G_0^2$ and introduces
\begin{equation}
    \mathbf G_2^3
    =\frac{1}{(c_2-c_1)\Delta t^2}
    \left[
        \frac{\overline{\mathbf G}{}^{\,2}-\overline{\mathbf G}{}^{\,0}}{c_2}
        -\frac{\overline{\mathbf G}{}^{\,1}-\overline{\mathbf G}{}^{\,0}}{c_1}
    \right] \ ,
\end{equation}
while the linear coefficient becomes $\mathbf G_1^3=\mathbf G_1^2-c_1\Delta t\,\mathbf G_2^3$.
For the fourth stage of MRK4, we retain a quadratic reconstruction and fit it to all four available stage values by least squares. This uses the additional stage information without increasing the polynomial degree, which we found damaging in our stability tests. We therefore determine
$\mathbf G_0^4$, $\mathbf G_1^4$, and $\mathbf G_2^4$ by minimizing
\begin{equation}
    \sum_{j=0}^{3}
    \left\|
        \mathbf G_0^4
        +\mathbf G_1^4 c_j\Delta t
        +\mathbf G_2^4(c_j\Delta t)^2
        -\overline{\mathbf G}{}^{\,j}
    \right\|^2
\end{equation}
at each horizontal location. The reconstructed forcing is sampled at the midpoint of each barotropic substep. Since $c_j$ are fixed, the least-squares weights depend only on $a$ and are computed once at startup. The procedure to compute the weights can be found in the repository supporting the manuscript~\cite{timesteppertestcases}. At each stage, the coefficients are obtained by applying these weights to the available stage values.

\section{Eigenvalue perturbation in a linear framework}
\label{sec:4}

Following the linear analysis framework of~\citeA{demange}, we examine how coupling affects the stability and the accuracy of the proposed Runge--Kutta schemes. 
Since the vertical modes of the continuous linear problem evolve independently, we can identify the errors introduced when the numerical discretization couples them together.
Computing the amplification factors of the evolution matrices is the same approach used by~\citeA{ducousso2026} to assess stability of their Runge--Kutta configurations, of which the ``Each-Stage strategy'' couples the barotropic solver with the baroclinic equations at every stage as in \citeA{hallberg1997}, \citeA{lanrk3}, \citeA{ishimwe2023} and the present work.
We expand this modal decomposition with an eigenvalue perturbation analysis to understand which choices of barotropic solver, filter, and slow forcing preserve stability, and how much damping an outer scheme needs to add to compensate for the coupling errors.

We linearize the non-rotating, inviscid hydrostatic equations in an $x$--$z$ plane about rest, with constant depth $H$, constant stratification $N^2$, keeping the vertical coordinate fixed $(\sigma = 1)$ since the motion of the levels does not change this analysis at linear order.
Momentum advection vanishes at linear order, while vertical advection of the background buoyancy remains. The perturbations satisfy
\begin{equation}
    \pd{u}{t} = -\pd{}{x}(g\eta + p) \ , \qquad
    \pd{p}{z} = b \ , \qquad
    \pd{b}{t} = -N^2 w \ , \qquad
    \pd{u}{x} + \pd{w}{z} = 0 \ ,
\end{equation}
where $p$ is the baroclinic pressure perturbation divided by the reference density, and $b$ is the buoyancy perturbation. The boundary conditions are
\begin{equation}
    w(-H)=0 \ , \qquad
    w(0)=\pd{\eta}{t} \ , \qquad
    p(0)=0 \ .
\end{equation}
We expand in horizontal Fourier modes and vertical modes $\phi_q(z)$, normalized to unit mean square over the depth. For a fixed horizontal wavenumber $k$,
\begin{equation}
    \frac{\widehat{u}(z,t)}{\sqrt{gH}}
    = \sum_{q=0}^{\infty} u_q(t)\phi_q(z) \ ,
    \qquad
    \frac{g\widehat{\eta}(t)+\widehat{p}(z,t)}{gH}
    = \sum_{q=0}^{\infty} h_q(t)\phi_q(z) \ ,
\end{equation}
where hats denote horizontal Fourier amplitudes. $u_q$ and $h_q$ are the dimensionless velocity and total-pressure amplitude of the vertical mode $q$ at the selected $k$.
The associated vertical eigenvalue problem and its derivation are given in Supporting Information, section~\ref{si-sec:appendix-modes}. Its solutions and their surface and depth-averaged weights are
\begin{equation}
    \begin{gathered}
    \phi_q = A_q\cos\Big[x_q\Big(1 + \frac{z}{H}\Big)\Big] \ , \qquad
    x_q\tan x_q = \varepsilon \defn \frac{N^{2}H}{g} \ , \qquad
    c_q = \frac{NH}{x_q} \ , \qquad
    \alpha_q \defn \frac{c_q}{\sqrt{gH}} = \frac{\sqrt{\varepsilon}}{x_q} \ , \\
    \psi_q \defn \phi_q(0) \ , \qquad
    \chi_q \defn \frac{1}{H}\int_{-H}^{0}\phi_q\,\mathrm{d}z = \psi_q\,\alpha_q^{2} \ , \qquad \Pi_q \defn \chi_q - \psi_q
    \end{gathered}
    \label{eq:modes}
\end{equation}
where $c_q$ is the horizontal wave speed, $x_q$ is the dimensionless vertical wavenumber, $\psi_q$ is the surface value of mode $\phi_q$, $\chi_q$ its depth average, and $A_q$ fixes the normalization. 
The barotropic mode is the root $x_0$ in $(0, \pi/2)$ and the baroclinic modes are the roots $x_q$ in
$(q\pi, q\pi + \pi/2)$. Expressions of $\alpha_q$, $\psi_q$ and $\chi_q$ for small $\varepsilon$ are given in Supporting Information, section~\ref{si-sec:appendix-modes}. 
Each vertical mode evolves independently according to
\begin{equation}
    \frac{\mathrm{d}u_q}{\mathrm{d}t}
    = -ik\sqrt{gH}\,h_q \ ,
    \qquad
    \frac{\mathrm{d}h_q}{\mathrm{d}t}
    = -ik\sqrt{gH}\,\alpha_q^2 u_q \ .
\end{equation}
The two propagation directions have exact amplification factor $e^{\pm i \theta_q}$ over $\Delta t$, where
\begin{equation}
    \theta_q \defn c_q\,k\,\Delta t = \alpha_q\,\mu_0 \ , \qquad \mu_0 \defn \sqrt{gH}\,k\,\Delta t \ .
    \label{eq:modal-phase}
\end{equation}
Here $\theta_q$ is the phase change of mode $q$ (which generalizes the $\theta$ for unit-speed advection of section~\ref{sec:2}), while $\mu_0$ is the barotropic Courant number that uses the reference speed $\sqrt{gH}$ (not $c_0$). 
Since the exact amplification factors have unit magnitude, any damping or coupling, in this framework, comes from the time discretization.

We now apply the numerical stages described in section~\ref{sec:3} to these modes. For this analysis, we keep the barotropic mode ($q = 0$) and $Q$ baroclinic modes. We collect the velocity amplitudes in $\boldsymbol{u}=(u_0,\ldots,u_Q)^T$ and the total-pressure amplitudes in $\vect{h}=(h_0,\ldots,h_Q)^T$ and use the same vector notation for the modal weights ($\vect{\chi}$ and $\vect{\psi}$). Applying $\vect{\chi}^T$ to a modal vector returns its scalar depth-average, while $\vect{\chi}$ applied to a scalar projects it onto the vertical modes.
As in section~\ref{sec:3-barotropic-damping}, $\bar{u}$ and $\eta$ denote the depth-averaged velocity and free surface Fourier amplitudes scaled by $\sqrt{gH}$ and $H$, respectively. The expansion of $\bar{u}$ and $\eta$ gives
\begin{equation}
    \bar{u} = \vect{\chi}^T \vect{u} \ , \qquad  \eta = \vect{\psi}^T \vect{h} \ ,  \qquad F = (\vect{\chi} - \vect{\psi})^T \vect{h} = \vect{\Pi}^T \vect{h} \ . 
\end{equation}
In this linearized configuration, the only slow forcing in the barotropic momentum equation is the baroclinic pressure gradient; this forcing is represented by the horizontal gradient of $F$, the depth-averaged baroclinic pressure divided by $gH$. Since $F$ in the barotropic momentum equation is treated in the same way as $\eta$, the barotropic equations integrated by the solver are
\begin{equation}
    \frac{d}{dt} \begin{bmatrix}
        \bar{u} \\ \eta
    \end{bmatrix} = - ik\sqrt{gH} \mathbf{J} \begin{bmatrix}
        \bar{u} \\ \eta + F
    \end{bmatrix} \ , 
\end{equation}
with \(\mathbf J\) defined in~\eqref{eq:windows}. When $F$ is held constant during the solve (the \emph{frozen} forcing), $\bar{u}$ and $\eta + F$ obey the unforced barotropic equations. Therefore, we can use $\mathbf{W}$ of section~\ref{sec:3-barotropic-damping} by adding $F$ to the initial free surface. Subtracting the initial state gives the increment returned by the solver
\begin{equation}
    \begin{bmatrix}
        \Delta \bar{u}^m \\ \Delta\eta^m
    \end{bmatrix} = (\mathbf{W} - \mathbf{I}) \begin{bmatrix}
        \bar{u}^n \\ \eta^n + F
    \end{bmatrix} \ .
    \label{eq:barotropcic-system}
\end{equation}
The solver also returns the transport used to advect the tracers. We denote the depth-average of this transport, normalized by $\sqrt{gH}$, by $\bar{u}^\dagger$. For wavenumber $k$, the continuity identity then reads $\Delta \eta^m = -i\mu_m \bar{u}^\dagger$. Therefore, in this setting, we can explicitly compute $\bar{u}^\dagger$ as
\begin{equation}
    \bar{u}^\dagger = \widetilde{W}_{11} \bar{u}^n + \widetilde{W}_{12}(\eta^n + F) \ , \qquad \text{where} \qquad [\widetilde{W}_{11} \quad \widetilde{W}_{12}] = \frac{[W_{21} \quad W_{22}-1]}{-i\mu_m} \ ,
\end{equation}
and $W_{ij}$ are the entries of the barotropic matrix $\mathbf{W}$. These equations only describe the evolution under a frozen forcing, the extension for reconstructed forcing is derived in Supporting Information, section~\ref{si-sec:app-stage-matrices}.

Combining the barotropic solve with the corrected tracer transport and velocity reconciliation gives the complete numerical timestep
\begin{equation}
    \vect{y}^{n+1}  = \begin{bmatrix}
        \vect{u}^{n+1} \\ \vect{h}^{n+1} 
    \end{bmatrix} = \mathbf{M} \begin{bmatrix}
        \vect{u}^{n} \\ \vect{h}^{n}  
    \end{bmatrix}  = \mathbf{M} \vect{y}^n\ .
\end{equation}
The matrices for each stage and their combination are given in Supporting Information, section~\ref{si-sec:app-stage-matrices}. The eigenvalues of $\mathbf{M}$ are the numerical amplification factors. 

\subsection{Coupled stability}
\label{sec:4-coupled}
The outer scheme and the barotropic solver each introduce damping, as shown in sections~\ref{sec:2} and~\ref{sec:3-barotropic-damping}. Coupling then can either add to this damping or oppose it, rendering the coupled scheme unstable through spurious growth. To understand when it produces growth, we write the complete step matrix as an uncoupled evolution plus a correction,
\begin{equation}
    \mathbf{M} = \mathbf{M}_\mathrm{u} + \mathbf{V} \ , 
\end{equation}
where $\mathbf{M}_\mathrm{u}$ advances the barotropic mode with $\mathbf{W}(\mu_0)$ and each baroclinic mode with the outer stability polynomial $R$. 
The difference $\mathbf{V}$ contains both the transfers between the modes and the changes to their individual evolution.
We take the amplification factors of $\mathbf{M}_\mathrm{u}$, which are stable within the uncoupled stability limits, as a reference and compute how $\mathbf{V}$ changes them. 
Each baroclinic mode has two uncoupled amplification factors $R( \pm i \theta_q)$, one for each propagation direction. We denote them by $\lambda_j$, where $j$ indicates both the mode and its propagation direction. We then search for an expansion that gives $\delta \lambda_j$ and $\delta\vect{r}_j$ for which
\begin{equation}
    \mathbf{M}(\vect{r}_j + \delta\vect{r}_j) = (\lambda_j + \delta\lambda_j)(\vect{r}_j + \delta\vect{r}_j) \ , \ \ \text{with} \ \ \mathbf{M}_\mathrm{u}\vect{r}_j = \lambda_j\vect{r}_j
\end{equation}
with $\vect{r}_j$ the right eigenvector of the uncoupled matrix. For isolated eigenvalues ($\lambda_j \ne \lambda_l$), the amplification factor correction due to the presence of $\mathbf{V}$, expressed as $\delta \lambda_j$, is
\begin{equation}
    \delta\lambda_j = V_{jj} + \sum_{l\neq j} \frac{V_{jl}V_{lj}}{\lambda_j - \lambda_l} + \mathcal{O}\left(\lVert \mathbf{V} \rVert^3\right) \  \ \ \text{where} \ \ V_{jl} = \vect{\ell}_j\mathbf{V}\vect{r}_l \ \ \text{and} \ \ \vect{\ell}_j\vect{r}_l = \delta_{jl} \ ,
    \label{eq:perturbation}
\end{equation}
with $\vect{\ell}_j$ the left eigenvector of $\mathbf{M}_\mathrm{u}$ associated with mode and propagation direction $j$. This expansion is derived in Supporting Information, section~\ref{si-sec:appendix-perturbation}.
To interpret these coefficients, we take one uncoupled eigenvector as the input to a timestep, $\vect{y}^n = \vect{r}_j$. Then, the uncoupled evolution matrix just multiplies this initial condition by $\lambda_j$. Contrarily, coupling both changes the coefficient and generates components along the other eigenvectors
\begin{equation}
    \mathbf{M} \vect{r}_j = (\lambda_j + V_{jj}) \vect{r}_j + \sum_{l \ne j} V_{lj} \vect{r}_l
\end{equation}
Therefore, $V_{jj}$ changes the initial condition \textit{directly}, while $V_{lj}$ measures the spurious transfer from $j$ to $l$ during this timestep.
We then call $V_{jj}$ the \textit{direct correction}.
The indirect correction instead involves a transfer from $j$ to another mode $l$ and back, through the product $V_{jl}V_{lj}$, weighted by their eigenvalue separation $\lambda_j - \lambda_l$.
This denominator makes the indirect correction particularly sensitive to modes with nearly coincident eigenvalues.
The contribution of these eigenvalue corrections to growth and damping can be approximated with
\begin{equation}
    |\lambda_j + \delta \lambda_j| - 1 \approx (|\lambda_j| -  1) + \frac{\operatorname{Re}(\overline{\lambda_j}\delta \lambda_j)}{|\lambda_j|} \ ,
    \label{eq:linearized-eig-equation}
\end{equation}
where $|\lambda_j| - 1$ is the uncoupled growth, where a negative value indicates damping. The second term is the leading change in magnitude caused by coupling: a positive value reduces the damping while a negative value adds to it. A positive correction causes the overall scheme to become unstable only if it exceeds the available uncoupled damping.
\begin{figure}[H]
    \centering
    \includegraphics[width=\linewidth,height=0.78\textheight,keepaspectratio]{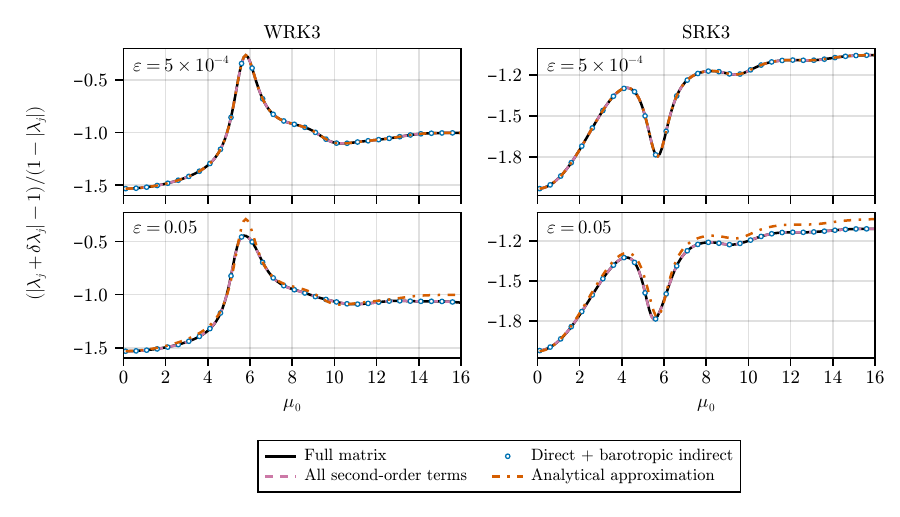}
    \caption{Retained coupling terms for the first baroclinic eigenvalue of WRK3 (left) and SRK3 (right), using RK3 substeps and F3, at $\varepsilon=5\cdot10^{-4}$ (top) and $\varepsilon=0.05$ (bottom). The full-matrix growth is compared with the complete second-order perturbation sum and with the sum retaining only direct and barotropic indirect terms, both evaluated with~\eqref{eq:linearized-eig-equation}. All values are divided by the damping $1-|\lambda_j|$. The calculation retains eight baroclinic modes and uses $48$ barotropic substeps. The additional dash-dotted lines plot the predictions from the analytical approximations.}
    \label{fig:c4-balance}
\end{figure}

Figure~\ref{fig:c4-balance} compares the growth computed from $\mathbf{M}$ with the perturbation estimate using \cref{eq:perturbation,eq:linearized-eig-equation}, retaining either all indirect terms ($l$) or only those involving the two barotropic eigenvalues ($l=b \in \ $barotropic). The two estimates nearly overlap and reproduce the growth computed with $\mathbf{M}$, showing that the exchanges with other baroclinic modes contribute little to the coupling. We consequently assume that, at leading order, it is possible to approximate the change in amplification factor with
\begin{equation}
    \delta \lambda_j \approx V_{jj} + \sum_{b\,\in\,\text{barotropic}} \frac{V_{jb}V_{bj}}{\lambda_j - \lambda_b} \ , 
\end{equation}
where $V_{bj}$, the reverse transfer, measures the barotropic growth generated from an initially baroclinic only state $\vect{r}_j$ and $V_{jb}$, the forward transfer, measures the baroclinic response to an initially barotropic only state $\vect{r}_b$. 

In Supporting Information, section~\ref{si-sec:app-transfers}, we derive the leading-order non-zero terms of $V_{jj}$, $V_{jb}$, and $V_{bj}$ for $\varepsilon \rightarrow 0$ at fixed $\mu_0$ with $\lambda_j \ne \lambda_b$
\begin{equation}
\begin{aligned}
    V_{jb} & \simeq \frac{\chi_q}{\sqrt2}\left[(\lambda_b-1)r_{b,1} + i\mu_0\left\{\mathbf{A}\vect{r}_b\right\}_2\right] \ , \qquad
    \mathbf{A} =
    \begin{cases}
        \mathbf{W}\left(\dfrac{\mu_0}{2}\right) & \text{(restarted)} \ , \\[8pt]
        \dfrac{1}{2}\left[\mathbf{I}+\mathbf{W}(\mu_0)\right] & \text{(blended)} \ ,
    \end{cases} \\[4pt]
    V_{bj} & \simeq \vect{\ell}_b\,\delta\vect{y}_0 \ , \qquad
    \delta\vect{y}_0 \defn
    \begin{bmatrix}
        \bar u^M-\bar u^n \\
        -i\mu_0(\bar u^{\dagger,M}-\bar u^n)
    \end{bmatrix} \ , \\[4pt]
    V_{jj} & \simeq \frac{\chi_q}{\sqrt2}
    \begin{cases}
        \bar u^M-\bar u^n+\dfrac{\mu_0^2}{2}\left(\bar u^{\dagger,M-1}-\bar u^n\right) & \text{(restarted)} \ , \\[8pt]
        \bar u^M-\bar u^n+\dfrac{\mu_0^2}{6}\left(2\bar u^{\dagger,1}+\bar u^{\dagger,2}-3\bar u^n\right) & \text{(blended)} \ ,
    \end{cases}
\end{aligned}
    \label{eq:vbj}
\end{equation}
where $\{\cdot \}_2$ indicates the second component of the vector, $\lambda_b = \lambda_b(\mu_0)$, and the barotropic eigenvectors use the normalization specified in the Supporting Information, section~\ref{si-sec:appendix-perturbation}. $M$ is the final stage and $M-1$ spans half a timestep, $\bar u^{\dagger,m}$ is the transport returned by the barotropic solve of stage $m$, and $\bar u^{\dagger,M}$ is the transport of the final stage. The restarted expressions are derived in~\eqref{si-eq:general-forward-transfer},~\eqref{si-eq:general-reverse-transfer} and~\eqref{si-eq:general-direct-projection}, and the blended ones in~\eqref{si-eq:ssp-forward-transfer},~\eqref{si-eq:ssp-general-reverse-transfer} and~\eqref{si-eq:ssp-general-direct-projection} of the Supporting Information. These expressions are valid for all the Runge--Kutta configurations considered here, independently of the slow forcing treatment.

The forward transfer $V_{jb}$ starts from an initially barotropic-only state, which produces only an $\mO(\e)$ slow forcing, so the forcing does not appear in its leading-order expression, and it affects the baroclinic mode only through the depth average $\chi_q$. The stage formulation differs only in $\mathbf{A}$, which represents the barotropic response that drives the last stage: the restarted stage uses the half-step state $\mathbf{W}(\mu_0/2)$, while SRK3 averages the initial and final states.

The reverse transfer and the direct correction ($V_{bj}$ and $V_{jj}$) start from a baroclinic-only state, whose depth-averaged velocity $\bar u^n = \chi_q / \sqrt{2}$ does not change at leading order in the uncoupled step, such that any change in the depth-averaged velocity, transport or free surface during the coupled step is spurious. 
Therefore, the reverse transfer projects this spurious barotropic state onto the barotropic eigenvectors, while the direct correction feeds the spurious velocity and the free surface of the spurious transports back onto mode $q$ within the same step.

Summing the indirect terms over the two barotropic directions, with $\lambda_j \approx 1$ in the denominators, gives (Supporting Information,~\eqref{si-eq:general-leading-coupling} and~\eqref{si-eq:ssp-general-leading-coupling})
\begin{equation}
\begin{gathered}
    \delta\lambda_j \simeq \frac{\chi_q}{\sqrt{2}}
    \begin{cases}
        \dfrac{\mu_0^2}{2}\left(\bar u^{\dagger,M-1}-\bar u^n\right)
        + i\mu_0\left\{\mathbf{W}\left(\dfrac{\mu_0}{2}\right)\left[\mathbf{I}-\mathbf{W}(\mu_0)\right]^{-1}\delta\vect{y}_0\right\}_2 & \text{(restarted)} \ , \\[10pt]
        \dfrac{\mu_0^2}{6}\left(2\bar u^{\dagger,1}+\bar u^{\dagger,2}-3\bar u^n\right)
        + \dfrac{i\mu_0}{2}\left\{\left[\mathbf{I}+\mathbf{W}(\mu_0)\right]\left[\mathbf{I}-\mathbf{W}(\mu_0)\right]^{-1}\delta\vect{y}_0\right\}_2 & \text{(blended)} \ ,
    \end{cases} \\[6pt]
    \left[\mathbf{I}-\mathbf{W}(\mu_0)\right]^{-1} = \sum_b \frac{\vect{r}_b\vect{\ell}_b}{1-\lambda_b} \ .
\end{gathered}
    \label{eq:generic-dlambda}
\end{equation}
The spurious velocity $\bar u^M-\bar u^n$ has opposite signs in the direct and the indirect corrections and cancels, so the coupling is the effect of a spurious free surface on the last stage.

With frozen forcing, the barotropic outputs in~\eqref{eq:generic-dlambda} are explicit functions of $\mathbf{W}$, and $\delta\lambda_j$ reduces to a closed form in $\mu_0$ and $\mathbf{W}$, given in~\eqref{si-eq:generic-total-coupling} for the restarted stage and in~\eqref{si-eq:ssp-total-coupling} for the blended stage. We refer to these frozen-forcing closed forms as the analytical approximation. They depend only on $\chi_q(\varepsilon)$, $\mu_0$ and $\mathbf{W}$, so they provide an easy method to assess the effect of the coupling with a generic barotropic solver without requiring to build the modal matrix. For the implicit free surface both are real and strictly negative~\eqref{si-eq:implicit-coupling} proving that coupling with an implicit free surface always adds to the damping.

\begin{figure}[H]
    \centering
    \includegraphics[width=\linewidth,height=0.78\textheight,keepaspectratio]{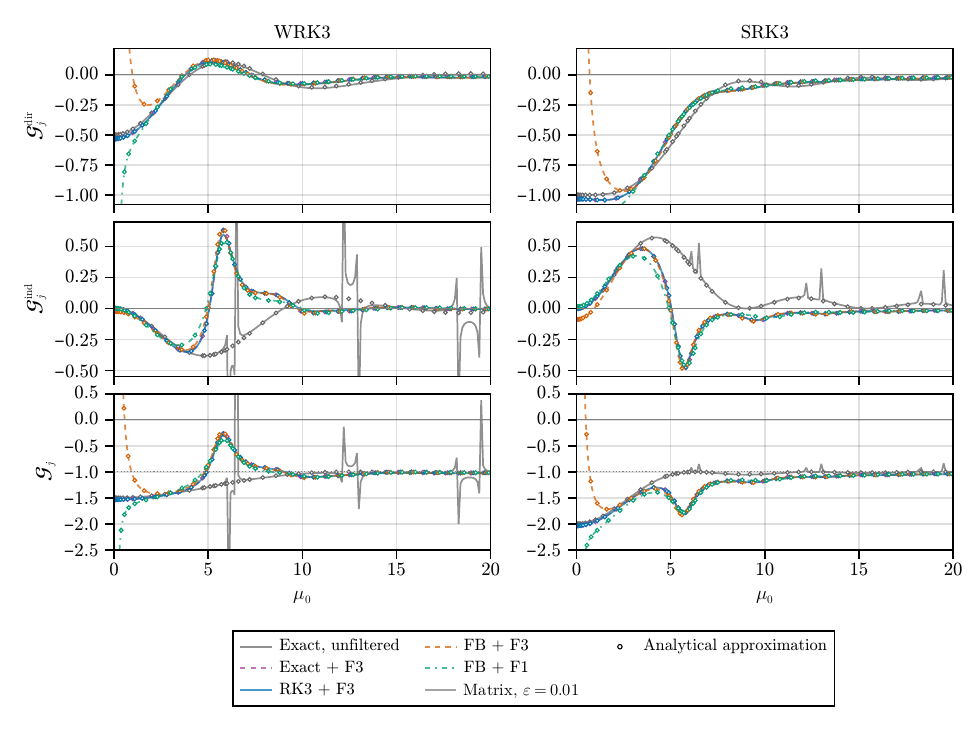}
    \caption{Damping balance for the first baroclinic eigenvalue of WRK3 (left) and SRK3 (right) at $\varepsilon=0.01$ using $48$ substeps. Rows show the magnitude corrections of the direct and barotropic indirect terms and the total growth, normalized by the uncoupled damping: $\mathcal{G}_j^{\mathrm{dir}}=\mathrm{Re}(\bar\lambda_j V_{jj})/[|\lambda_j|(1-|\lambda_j|)]$, $\mathcal{G}_j^{\mathrm{ind}}=\mathrm{Re}(\bar\lambda_j\,\delta\lambda_{j,\mathrm{ind}})/[|\lambda_j|(1-|\lambda_j|)]$ with $\delta\lambda_{j,\mathrm{ind}}=\sum_b V_{jb}V_{bj}/(\lambda_j-\lambda_b)$, and $\mathcal{G}_j=(|\lambda_j+\delta\lambda_j|-1)/(1-|\lambda_j|)$. Circles show the analytical approximation, while curves show exact terms computed from the full matrix with with $8$ retained modes. Values above zero in the last row indicate instability.}
    \label{fig:total-growth-validation}
\end{figure}

The split--explicit correction can change sign depending on the barotropic solver and filter.
Figure~\ref{fig:total-growth-validation} separates the direct and indirect contributions for the first baroclinic mode at $\varepsilon = 0.01$, computed from $\mathbf V$ with 8 retained baroclinic modes, and compares them with the analytical approximation showing very good agreement between the two. The comparison at $\varepsilon = 5\cdot 10^{-4}$ and $\varepsilon = 0.05$ in Figure~\ref{fig:c4-balance} shows how this agreement changes with the stratification: the analytical expressions capture the magnitude at small $\varepsilon$ and retain the main trends at larger $\varepsilon$. Figure~\ref{fig:total-growth-validation} also highlights four features of the coupling corrections. 

\textit{(i), Resonances.} With exact, unfiltered barotropic evolution, the indirect correction has sharp peaks, which filtering removes. The peaks come from the denominator $\lambda_j - \lambda_b$ reaching zero when the eigenvalues of the separate, uncoupled solvers approach each other. uncoupled barotropic evolution is undamped, so $|\lambda_b| = 1$ nearly coincides with the weakly damped baroclinic eigenvalues whenever the phases align. This is rather simple to verify since it depends only on the uncoupled evolution, writing $\lambda_j \approx e^{i\theta_j}$ and $\lambda_b = \rho_b e^{i\theta_b}$ gives
\begin{equation}
    |\lambda_j-\lambda_b|^2
    \approx
    (1-\rho_b)^2
    +4\rho_b
    \sin^2\!\left(\frac{\theta_j-\theta_b}{2}\right) \ge (1 - \rho_b)^2\ , 
\end{equation}
where $\rho_b=1$ for the unfiltered solve, so the denominator can become arbitrarily small as $\theta_b$ approaches $\theta_j+2\pi n$. 
Therefore, with a quasi-neutral baroclinic mode $|\lambda_j| \approx 1$, filtering supplies a minimum separation of approximately $1 - \rho_b$ even when the phases align, explaining why the filtered curves lose the sharp resonance peaks. The peak heights themselves are not present in the analytical approximation (the circles) because the eigenvalue perturbation formula assumes separated eigenvalues. This substantiates that the cause of the instability is the eigenvalues approaching each other rather than the numerator ($V_{jb}V_{bj}$) growing out of bound.

\textit{(ii), Effect of filtering on the sign of the correction.} Beyond $\mu_0 \simeq 4.5$, the filtered and unfiltered indirect corrections have opposite signs, so filtering can make a stabilizing coupling destabilizing and viceversa. Here the eigenvalues are well separated, so the change comes from the transfers rather than from the denominator. A simple example is the exact evolution averaged uniformly over $0 \le \tau\le 2$, for which
\begin{equation}
    \lambda_b^{\text{filt.}}(\mu_0)
    =\frac12\int_0^2 e^{-i\mu_0\tau}\,\mathrm{d}\tau
    =e^{-i\mu_0}\frac{\sin\mu_0}{\mu_0} \ .
\end{equation}
Averaging multiplies the exact eigenvalue by $\sin \mu_0 / \mu_0$ which always damps the eigenvalue but also reverses the sign of the eigenvalue whenever $sin \mu_0 < 0$. Since the transfer coefficients depend on the barotropic eigenvalue, filtering can reverse the sign of the indirect correction (as shown in the figure).

\textit{(iii), Difference between WRK3 and SRK3.} The indirect corrections have opposite sign over much most of the range with $\mu_0 \ge 5$, both with and without filtering. The difference comes from the barotropic state used in the forward transfer $V_{jb}$ which differs between the two formulations. WRK3 uses the half-step matrix $\mathbf{W}(\mu_0 / 2)$ while SRK3 averages the initial state and the full-step response $[\mathbf I + \mathbf W(\mu_0)] /2$. A simple example is the exact evolution at $\mu_0 = 2\pi$ for which $\lambda_b(2\pi) = 1$ and $\lambda_b(\pi) = -1$, then the forward response gives
\begin{equation}
V_{jb}^{\text{WRK3}} = - V_{jb}^{\text{SRK3}} = - \sqrt{2} i\pi \chi_q r_{b, 2} \ , 
\end{equation}
so that the forward transfers of the two schemes are equal and opposite.

\textit{(iv), Direct correction and FB substeps.} The direct correction is smooth and depends little on the solver over most of the range because it involves only entries of the full- and half-step barotropic matrices, and no eigenvalue separation. At small timesteps, however, FB substeps produce a sharp positive contribution with FB that can overwhelm the outer damping.
Deriving the direct correction for $\mu_0 \rightarrow 0$ (Supporting Information, section~\ref{si-sec:appendix-small-step}) gives
\begin{equation}
    V_{jj}
    \simeq
    \frac{\chi_q^2\mu_0^2}{4}
    \begin{cases}
        (1-\alpha_q^2)\dfrac{\gamma_B}{\gamma_m}
        -(1+\alpha_q^2)M_2, & \text{FB},\\[5pt]
        -(1+\alpha_q^2)M_2, & \text{RK3},
    \end{cases} \ .
\end{equation}
For RK3, this contribution is damping when $M_2>0$ and vanishes for F3, which has $M_2=0$. 
FB, on the other hand, adds a positive term proportional to $\mu_0^2$. 
Since the outer scheme's damping is quartic, this positive quadratic term dominates as $\mu_0 \rightarrow 0$.
A filter needs to compensate at the same quadratic order, requiring
\begin{equation}
    M_2
    \geq
    \frac{\gamma_B}{\gamma_m}
    \frac{1-\alpha_q^2}{1+\alpha_q^2} \ .
    \label{eq:fb-requirement}
\end{equation}
This explains the opposite small-step behavior for FB with F1 and F3. F3 has $M_2 = 0$ and leaves the positive quadratic term uncompensated, while F1 supplies enough damping to offset the positive contribution.
Since $\alpha_q$ is small for baroclinic modes, we can rewrite~\eqref{eq:fb-requirement} as
\begin{equation}
    M_2 \gtrsim \frac{2}{N_S}
\end{equation}
where $N_S$ is the number of substeps. Therefore, the required quadratic filter damping decreases as the number of substeps increases. 

\begin{figure}[H]
    \centering
    \includegraphics[width=\linewidth,height=0.78\textheight,keepaspectratio]{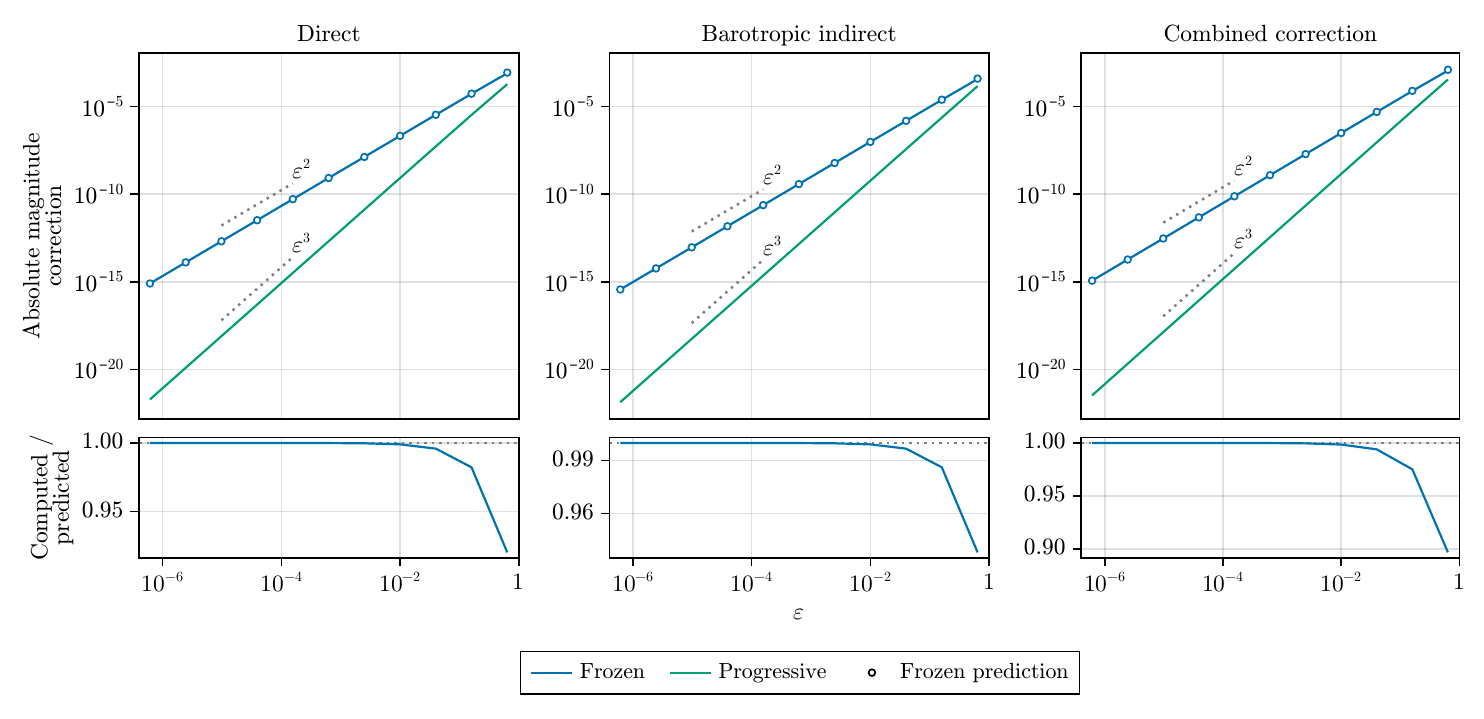}
    \caption{Absolute direct, barotropic indirect, and combined magnitude corrections for WRK3 at $\mu_0=2$, using exact, unfiltered barotropic integration and continuously integrated forcing. Curves compare frozen and progressive forcing. The circles show the leading frozen-forcing predictions, with four retained baroclinic modes. The lower panel shows the ratio of the computed and predicted correction.}
    \label{fig:forcing-reconstruction}
\end{figure}

We next use the stage information to reduce the coupling error itself. 
Figure~\ref{fig:forcing-reconstruction} compares frozen and progressive forcing for WRK3 using exact, unfiltered barotropic integration.
The computed frozen corrections scale as $\varepsilon^2$, while the progressive corrections scale as $\varepsilon^3$. The difference comes from how the forcing follows the free surface during each barotropic solve.

As shown by~\eqref{eq:vbj}, the coupling correction is driven by the spurious change in the depth-averaged velocity and transport that the coupled evolution generates from a baroclinic-only initial state $\vect{y}^n = \vect{r}_j$. During each barotropic solve, the evolution of the barotropic variables approximately follow 
\begin{equation}
\frac{d\bar u}{dt'}
\approx -i\frac{\mu_0}{\Delta t}\bigl(\eta+F^m(t')\bigr),
\qquad
\frac{d\eta}{dt'}
\approx -i\frac{\mu_0}{\Delta t}\bar u \ .
\end{equation}
For this initial state, the free surface should change linearly in time following an almost constant $\bar u$. This requires that $\eta(t^\prime) \approx - F^m(t^\prime)$. This can be achieved with progressive forcing reconstruction which (beyond the first stage) gives
\begin{equation}
    \eta(t^\prime) - \eta^n \simeq - i \frac{\mu_0}{\Delta t}\bar u^n t^\prime \ , \quad
    F^m(t^\prime) - F_0^1 \simeq + i \frac{\mu_0}{\Delta t}\bar u^n t^\prime \ ,
\end{equation}
so that the two changes cancel in $\eta + F^m$. Therefore, at leading order $d \bar u / dt \approx 0$ as in the uncoupled evolution, leaving $\bar u^m \approx \bar u^{\dagger, m} \approx \bar u^n$, so the spurious difference in the transfer coefficients disappears. The same cannot be achieved with frozen forcing for which $dF^m/dt'=0$. This cancellation occurs when the slow forcing is only the the pressure gradient since $F^m$ behaves like $\eta$ in the barotropic evolution. The same reconstruction performed for other slow terms does not necessarily lead to the same benefits, as shown in Appendix~\ref{sec:rotation-results} for the Coriolis force.

Reducing the coupling error allows us to use MRK4, which has a lower outer damping than the three-stage schemes. Its parameter $a$ controls the damping and the uncoupled stability through
\begin{equation}
    |R(-i\theta, a)|^2 = 1 - \frac{1 - 4a}{12}\theta^4 + \frac{1 -6a}{36}\theta^6 + \frac{a^2}{36}\theta^8 \ .
\end{equation}
At $a = 1/4$ the quartic damping vanishes and the uncoupled scheme is fourth order for linear problems.

\begin{figure}[H]
    \centering
    \includegraphics[width=\linewidth,height=0.78\textheight,keepaspectratio]{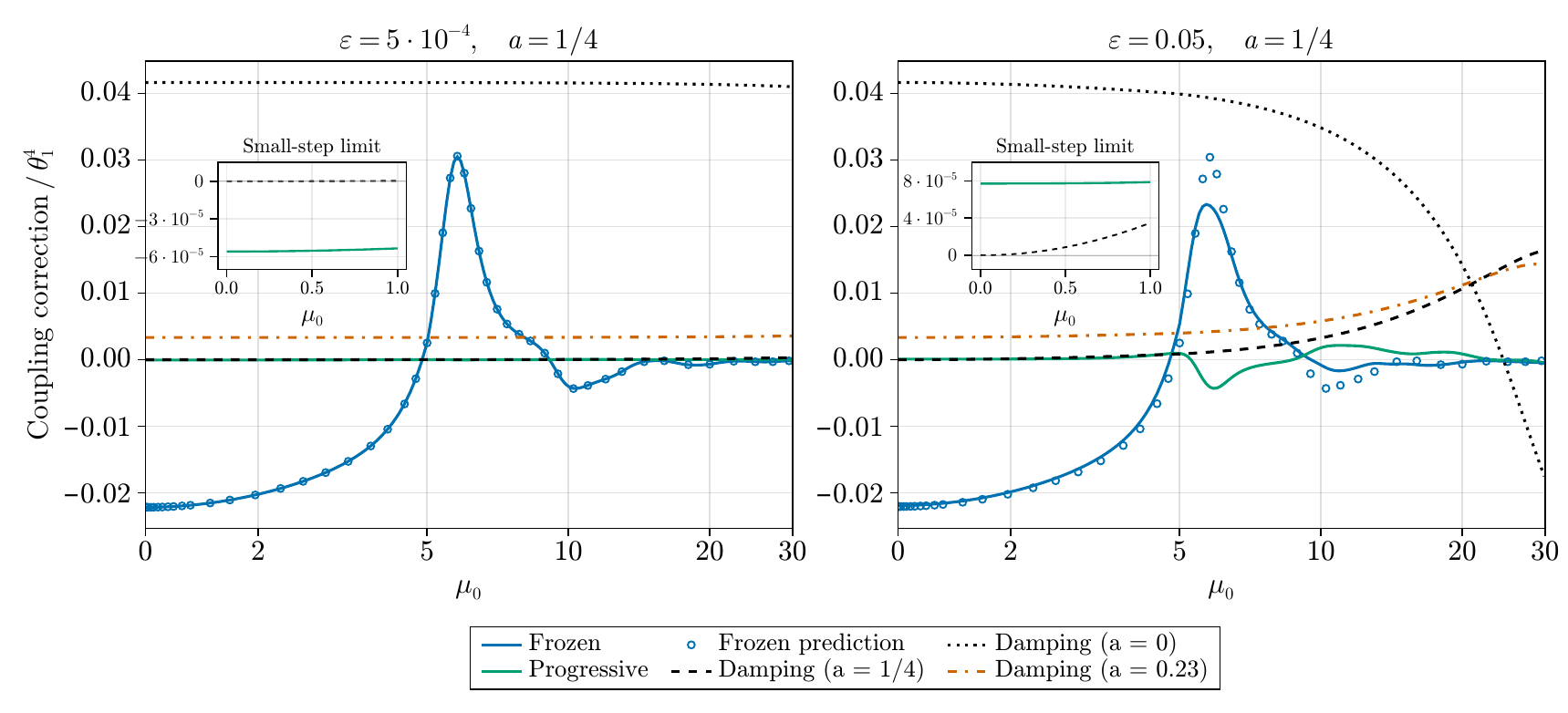}
    \caption{Normalized coupling corrections $(|\lambda_j + \delta \lambda_j| - |\lambda_j|) / \theta_1^4$ for the first baroclinic eigenvalue of MRK4 at $a = 1/4$ for $\varepsilon = 5\cdot 10^{-4}$ and $\varepsilon = 0.05$ using 48 substeps and 24 retained baroclinic modes. Solid curves compute the terms from the full coupled matrix, while circles use the frozen-forcing analytical approximation. The dashed and dotted curves show the uncoupled damping for different values of the parameter $a$.}
    \label{fig:mrk4-forcing-decomposition}
\end{figure}

\begin{figure}[H]
    \centering
    \includegraphics[width=\linewidth,height=0.78\textheight,keepaspectratio]{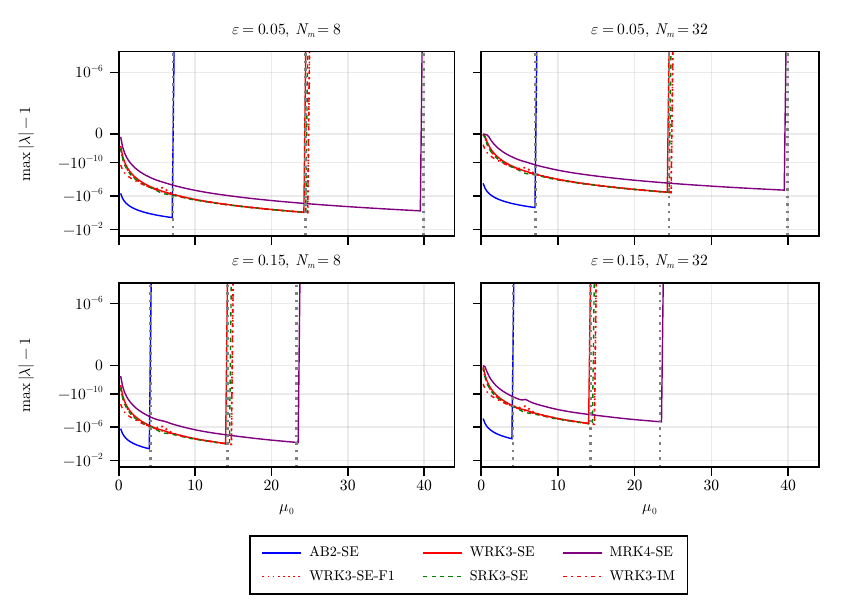}
    \caption{Coupled growth per step of the configurations of Table~\ref{tab:discretization} against the barotropic Courant number, at two stratifications and with $8$ and $32$ retained baroclinic modes. Dotted vertical lines mark $0.29\sqrt{3}/\alpha_1$, $\sqrt{3}/\alpha_1$ and $2.79/\alpha_1$.}
    \label{fig:schemes}
\end{figure}

Figure~\ref{fig:mrk4-forcing-decomposition} shows a positive frozen-forcing peak near $\mu_0 \approx 6$ that exceeds the uncoupled damping at $a = 1/4$. Using progressive forcing, we can reduce the peak below the uncoupled damping by suppressing the leading order corrections in $V_{bj}$ and $V_{jj}$. The improvement is a consequence of reducing the barotropic motion excited by mode $q$ and its feedback on $q$ itself (the $\lambda_b - \lambda_j$ mediated $V_{bj} \rightarrow V_{jb}$ feedback).

The insets in Figure~\ref{fig:mrk4-forcing-decomposition} show the small-step limit of the progressive forcing, which we need to treat differently. We define
\begin{equation}
    \Gamma_j = \lim_{\mu_0 \rightarrow 0} \frac{\operatorname{Re}(\overline{\lambda_j}\delta\lambda_j)}{|\lambda_j|\theta_q^4} \ , 
\end{equation}
where $\Gamma_j$ is the small-step coefficient of the net correction. 
A positive (destabilizing) $\Gamma_j$ is offset by quartic outer damping when $(1 - 4a)/24 > \Gamma_j$. For exact, unfiltered barotropic evolution with frozen forcing, the small-step expansion of~\eqref{si-eq:generic-total-coupling} gives $\delta\lambda_j\simeq-\chi_q^2\mu_0^4/96$.
Since $\lambda_j \approx 1$ and $\chi_q^2 \approx 2\alpha_q^4$, the leading magnitude correction gives
\begin{equation}
    \Gamma_j \approx - \frac{1}{48}
\end{equation}
Therefore, coupling with a frozen forcing adds damping in this small-step limit.
However, as seen in the insets of Figure~\ref{fig:mrk4-forcing-decomposition}, progressive forcing cancels the leading-order $\varepsilon^2$ term only for an exact, continuously forced evolution. $\Gamma_j$ with a practical scheme such as the one shown in the figure does not cancel but keeps a small residual which might be positive or negative. For this reason, $a = 1/4$, which removes the quartic outer damping at small timesteps, is not a suitable choice. We choose $a = 0.23$ to retain a quartic damping with coefficient $1/300$ as a safe choice.

The stability of the proposed schemes is shown in Figure~\ref{fig:schemes}, where we retain $8$ and $32$ modes at $\varepsilon = 0.05$ and at $\varepsilon=0.15$.
All split--explicit configurations use 48 nominal barotropic substeps. WRK3-SE and MRK4-SE use RK3 substeps, the F3 filter, and progressive forcing. SRK3-SE uses the same substeps and filter with frozen forcing. AB2-SE and WRK3-SE-F1 use forward--backward substeps, the F1 filter, and frozen forcing. WRK3-IM uses the implicit free-surface solver.
The stability limit of each scheme corresponds to the zero crossing. MRK4 retains the widest stability limit, while increasing the number of retained modes changes the stability limits very little.

\begin{figure}[H]
    \centering
    \includegraphics[width=\linewidth,height=0.78\textheight,keepaspectratio]{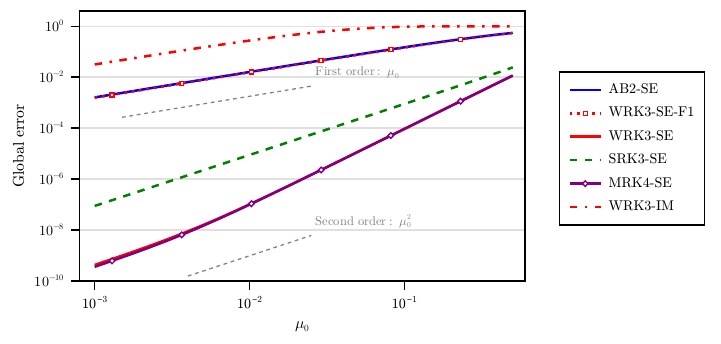}
    \caption{Global error against the barotropic Courant number for the configurations of Table~\ref{tab:discretization}, with $\varepsilon=0.05$ and two retained baroclinic modes. 
    Global error in the $2$-norm at $n\Delta t=64/(kc_0)$, with modal velocities normalized by $c_0$ and modal thicknesses by $H\alpha_q$. The split--explicit configurations use $48$ nominal barotropic substeps: FB with F1 for AB2-SE and WRK3-SE-F1, and RK3 with F3 for WRK3-SE, SRK3-SE, and MRK4-SE.
    WRK3-IM uses the implicit free-surface solver.}
    \label{fig:order-ceiling}
\end{figure}

The order of accuracy of these schemes is shown in Figure~\ref{fig:order-ceiling}. A first-order filter such as F1 departs from the exact advance by $\mathcal{O}(\mu_0^2)$ per step, as~\eqref{eq:filter-transfer} shows, and limits the barotropic mode to first-order accuracy. Consequently, replacing AB2-SE with WRK3-SE-F1 does not improve the convergence order in this linearized test, although their behavior differs in the practical configurations considered in the next section. Progressive forcing reduces the error without increasing the convergence order: errors in the stage values used for reconstruction leave a second-order global error. SRK3-SE, which uses frozen forcing, is also second order, but its error is about two orders of magnitude larger than that of WRK3 and MRK4, which use progressive forcing. The implicit solve produces a first-order global error in the barotropic mode, consistent with the IMEX analysis of~\citeA{kar}. The configurations using RK3 substeps and F3 are second order.

The analysis above neglects rotation. In a rotating ocean, the Coriolis force acting on the depth-integrated flow is part of the slow forcing of the barotropic solve, and it couples the inertial oscillations with the barotropic gravity waves.
Appendix~\ref{sec:rotation-results} extends the modal analysis to an $f$-plane, with the derivation in Supporting Information, section~\ref{si-sec:appendix-rotation}. The admissible timestep depends on the gravity-wave Courant number and $|f|\Delta t$. With frozen Coriolis forcing, WRK3, MRK4 and SRK3 remain stable well above the values of the simulations in section~\ref{sec:5}. Reconstructing the pressure gradient progressively raises this limit, reconstructing the Coriolis force together can lower it, and SRK3 has the largest stability region.

\section{Results and validation}
\label{sec:5}
We now test the discretizations of section~\ref{sec:3} with a suite of test cases, designed to highlight differences in performance across different dynamical regimes. The discretizations, collected in Table~\ref{tab:discretization}, are organized so as to isolate individual aspects of the discretizations: Oceananigans' quasi-AB2 (AB2-SE, at $0.29$ times the third-order Runge--Kutta timestep); the same barotropic scheme with a different outer time discretization (WRK3-SE-F1); the improved barotropic treatment (WRK3-SE, with RK3 substeps and the filter F3 of section~\ref{sec:3-barotropic-damping}); the strong-stability-preserving alternative with identical settings (SRK3-SE); and the four-stage scheme (MRK4-SE), run at $1.6\,\Delta t$ (the ratio of its imaginary-axis limit to that of the three-stage schemes). The implicit branch (WRK3-IM) is run as a benchmark in the idealized cases to test the results of an unconditionally stable barotropic solver, and a low-order tracer-advection variant (WRK3-UP) isolates the interplay between temporal and spatial dissipation.
We assess the numerical proficiency of the discretizations based on the resolved kinetic energy and the numerical mixing introduced.
To estimate the numerical mixing, we implement the method described by \citeA{banerjee2024} that provides an exact value (local in time and space) for the numerically induced dissipation of tracer variance in the $x$, $y$, and $z$ directions ($P_x$, $P_y$, and $P_z$). 
Appendix~\ref{sec:appendix-B} extends it to AB2 and the Runge--Kutta schemes.
With the numerical variance production, we estimate the numerical diffusivity as
\begin{equation}
    \kappa_{\rm num} = - \frac{1}{2}\left(\overline{P}_x + \overline{P}_y + \overline{P}_z\right) \left(\overline{ | \bnabla b|^2 + (\partial_z b)^2}\right)^{-1} \ ,
    \label{eq:numerical-diffusivity}
\end{equation}
where the bar indicates the time or space averaging used for each test case.
Averaging is introduced to avoid divisions by zero when the flow is locally unstratified.
\begin{table}[t]
    \centering
    \begin{tabular}{ccccccc}
         \hline
         case & outer scheme & free surface & substep & filter & timestep & tracer advection \\
         \hline
         \textbf{AB2-SE}       & AB2  & split--explicit & FB  & F1  & $0.29\,\Delta t$  & WENO, 7th order \\
         \textbf{WRK3-SE-F1} & WRK3 & split--explicit & FB  & F1 & $\Delta t$    & WENO, 7th order \\
         \textbf{WRK3-SE}      & WRK3 & split--explicit & RK3 & F3   & $\Delta t$    & WENO, 7th order \\
         \textbf{SRK3-SE}      & SRK3 & split--explicit & RK3 & F3      & $\Delta t$    & WENO, 7th order \\
         \textbf{MRK4-SE}      & MRK4 & split--explicit & RK3 & F3      & $1.6\,\Delta t$ & WENO, 7th order \\
         \textbf{WRK3-IM}      & WRK3 & implicit        & --- & ---     & $\Delta t$    & WENO, 7th order \\
         \textbf{WRK3-UP}      & WRK3 & split--explicit & RK3 & F3      & $\Delta t$    & Upwind, 3rd order \\
         \hline
    \end{tabular}
    \caption{Every case is run at the same fraction of the stability limit of its own scheme on the first baroclinic mode. The idealized cases use seventh-order WENO for the tracers; the global configuration uses ninth-order WENO and runs AB2-SE at $0.29\,\Delta t$.}
    \label{tab:discretization}
\end{table}
In the idealized cases, we use a 7th-order WENO reconstruction for tracer advection and, where not explicitly stated, the momentum advection uses the WENO vector-invariant scheme, as described by \citeA{simoneWENO}.
WRK3-UP uses the same progressive forcing and barotropic solver as WRK3-SE. Only the tracer advection differs, allowing us to distinguish the mixing due to spatial discretization from that due to time discretization.
\begin{figure}[H]
    \centering
    \includegraphics[width=\linewidth,height=0.78\textheight,keepaspectratio]{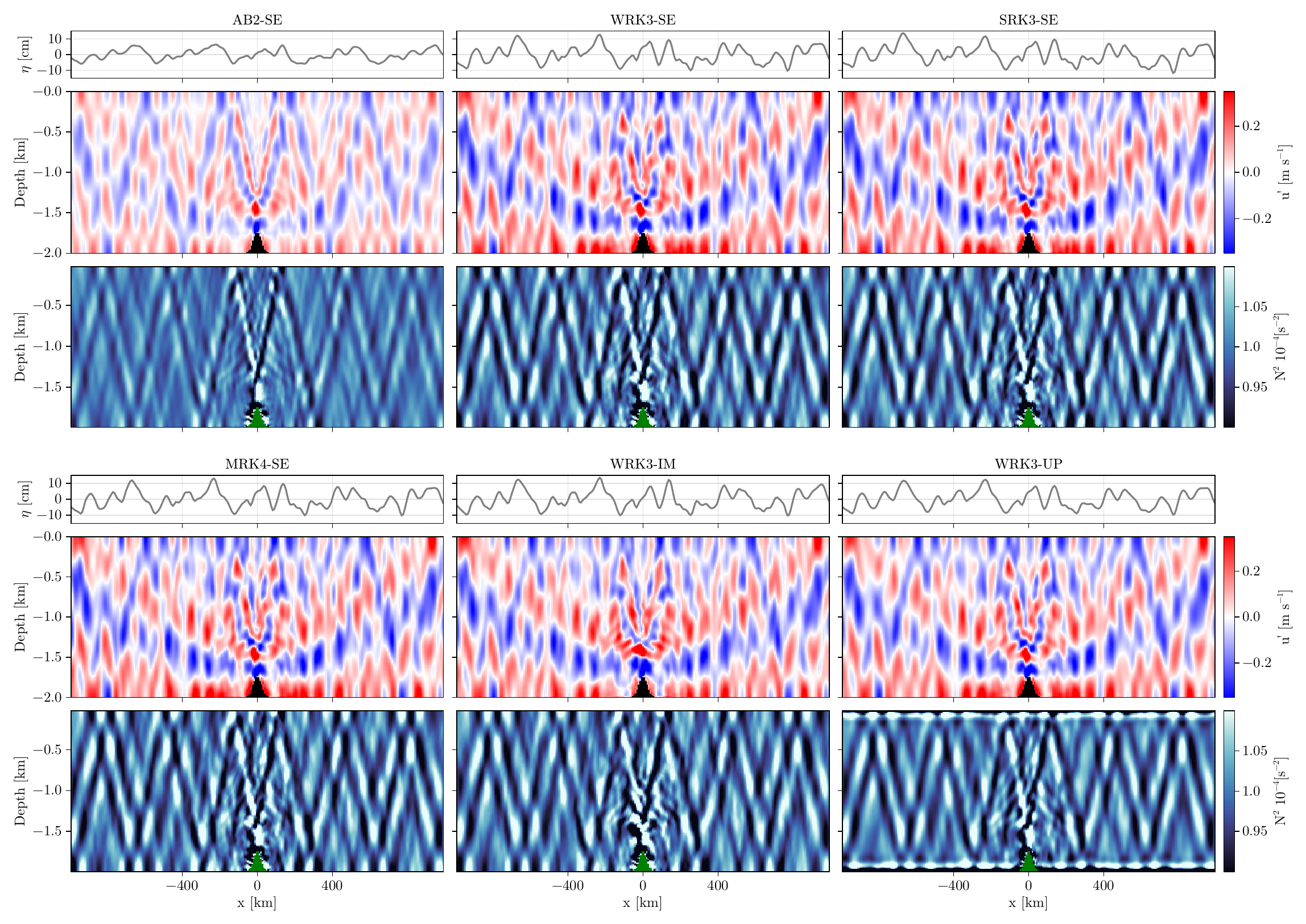}
    \caption{The internal tide configuration after 40 days of evolution, for six of the seven cases of
    Table~\ref{tab:discretization}. The top panels show the free surface elevation, while the center and bottom panels show the instantaneous velocity fluctuations $u^\prime = u - V^{-1}\int_V u \, \mathrm{d}V$ and the stratification $N^2 = \partial_z b$.}
    \label{fig:internal-tide-snapshot}
\end{figure}

\subsection{Internal tide}
The internal tide configuration mirrors the Oceananigans documentation example of the same name \cite{oceananigans-docs}.
It consists of a two-dimensional domain, 2000 kilometers long and 2 kilometers deep, initialized with a constant stratification $N^2 = 10^{-4}$ s$^{-2}$, and forced by a zonal oscillatory tide applied as a body forcing:
\begin{equation}
    F_U = \frac{\omega^2 - f^2}{\omega} U \sin{\left(\omega t\right)} \ ,  \ \ \text{with} \ \ \ \begin{cases}
        \ U \!\! &= 0.1 \omega D \ ,\\
        \ D \!\! &= 20\text{ km} \ ,\\
        \ \omega \!\! &= 2\pi / T  \ ,\\
        \ T \!\! &= 12.421\text{ hr} \ , \\
        \ f \!\! & = -1.03 \cdot 10^{-4} \text{s}^{-1} \ .
    \end{cases}
\end{equation}
The forced tide flows over a 250\,meter-tall Gaussian seamount with a standard deviation of 20~km, represented with a staircase bathymetry; the tidal excursion (2~km) and Froude number ($\approx 0.11$) are small, and the topography is subcritical, with a maximum slope of about 0.8 times that of the tidal rays.
The horizontal and vertical resolutions are 7.8125~km and 15.625~m (128 layers spanning 2~km), respectively.
This test case isolates the role of time discretization: the spatial advection operator, which is usually dominant in causing numerical mixing, plays a secondary role because this setup is mostly linear.
Therefore, we expect to see results that closely match the analytical predictions derived in section~\ref{sec:4}.
The timestep is set to the stability limit of each scheme on the first baroclinic mode. This gives $\Delta t = 11.3$~minutes for the three-stage schemes, $3.3$~minutes for AB2-SE and $18.1$~minutes for MRK4-SE. For the SE cases, we use $N_S = 64, \ 96$ for the RK3 substep for the three-stage RK and MRK4, respectively, and $N_S = 96, \ 32$ for the FB substeps for WRK3-SE-F1 and AB2-SE, respectively.
For momentum advection, in this 2D-setting, we use a 5th order WENO scheme.
Simulations are initialized with $u = U$ and run for 40 days.
No drag or explicit mixing is applied, so any tracer or momentum diffusion stems solely from the numerical mixing introduced by the time and space discretizations.

Figure~\ref{fig:internal-tide-snapshot} shows snapshots of free-surface elevation $\eta$, velocity fluctuations $u^\prime$, and stratification $N^2$ after 40~days.
The Runge--Kutta solutions are more energetic than the AB2-SE one: both $u^\prime$ and $N^2$ exhibit sharper structures than in AB2-SE, indicating the weaker dissipation of the Runge--Kutta schemes. 
Away from the bottom the Runge--Kutta configurations radiate beams that are visually indistinguishable at this resolution; only the WRK3-UP case shows a somewhat smoother stratification, owing to the larger dissipation introduced by the advection scheme.

\begin{figure}[H]
    \centering
    \includegraphics[width=\linewidth,height=0.78\textheight,keepaspectratio]{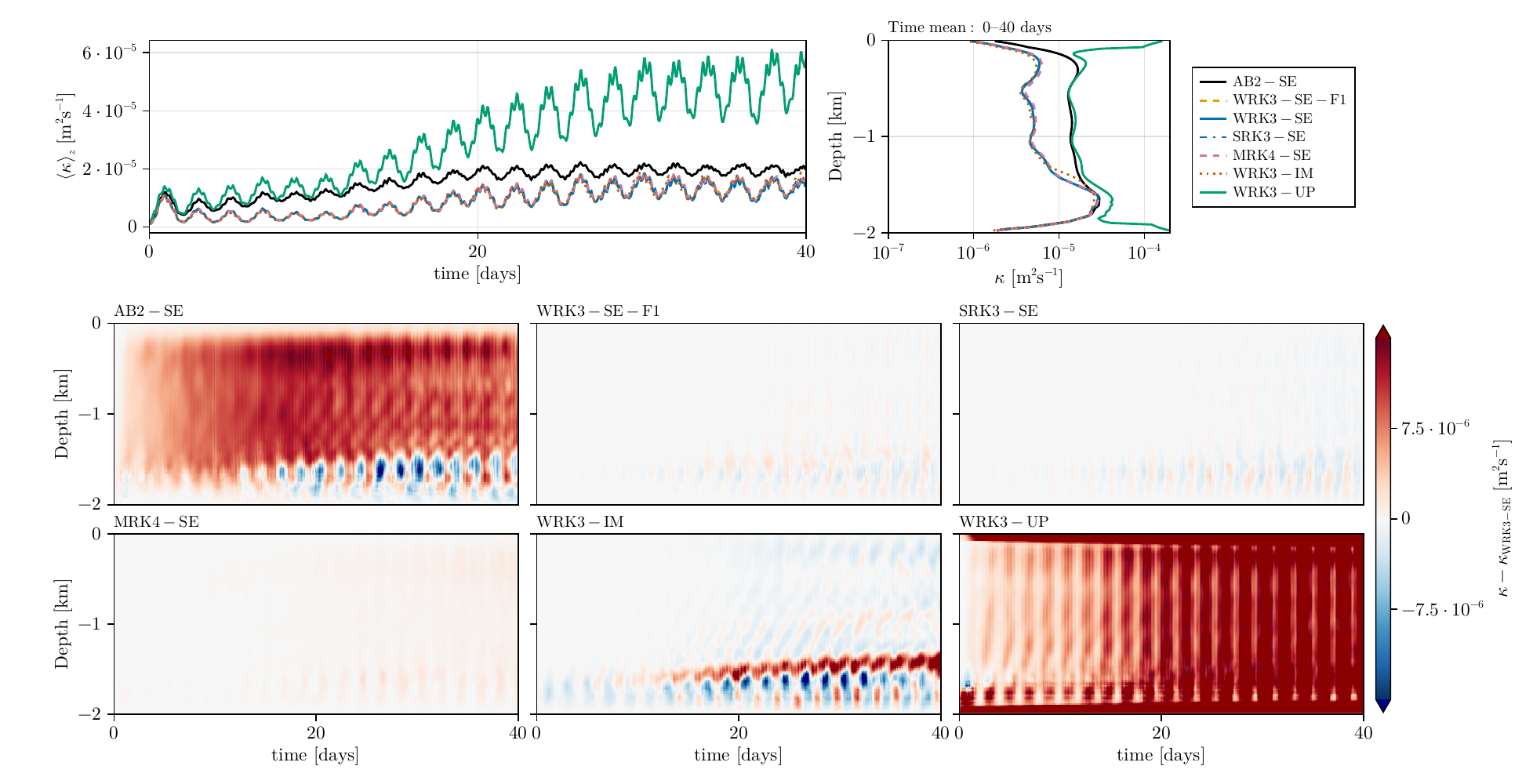}
    \caption{Numerical diffusivity of the internal tide configuration, computed via~\eqref{eq:numerical-diffusivity} with the bar indicating an average over the zonal $x$-direction.
    The upper left panel shows the space-averaged diffusivity of each case and the upper right the time-averaged diffusivity profile. The lower panels show $\kappa$ as a function of depth and time, as a difference against WRK3-SE, so that red marks a case mixing
    more and blue less.}
    \label{fig:diffusivity-internal-tide}
\end{figure}

Figure~\ref{fig:diffusivity-internal-tide} shows the numerical diffusivity of all the configurations.
Over the whole depth, AB2-SE is about twice as diffusive as WRK3-SE, with the difference concentrated in the interior of the domain. Here AB2-SE reaches a diffusivity comparable to the low-order spatial discretization of WRK3-UP. 
Near the bottom, where the diffusivity peaks due to interaction with a staircase bathymetry (which reduces the order of the advective operators), AB2-SE is as diffusive as WRK3-SE. It is also locally less diffusive near the seamount, which we attribute to the lower kinetic energy it retains, since the time discretization, contrarily to WRK3-UP, damps the momentum in addition to the tracer field.
Notably, AB2-SE is more diffusive than the Runge--Kutta configurations while also taking a timestep $0.29$ times as long.
WRK3-SE-F1 is indistinguishable from WRK3-SE, so that in this configuration the barotropic treatment, which section~\ref{sec:4} shows to govern the stability of the coupling, leaves no measurable signature on the numerical mixing.
MRK4-SE is likewise indistinguishable while advancing at $1.6$ times the timestep. 
SRK3-SE does not deviate from WRK3-SE either, indicating that the Runge--Kutta arrangement (once stable) does not affect mixing in this configuration.

\begin{figure}[H]
    \centering
    \includegraphics[width=\linewidth,height=0.78\textheight,keepaspectratio]{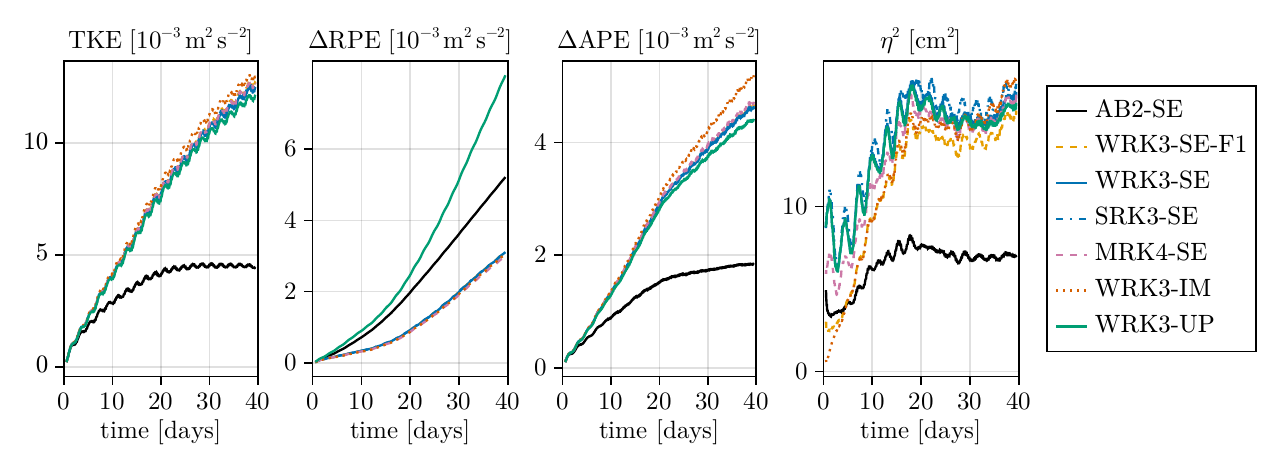}
    \caption{Energetics of the internal tide configuration for the different numerical discretizations. All quantities are volume averaged and shown as running averages over about two tidal periods to remove the oscillation at the tidal frequency.}
    \label{fig:energetics-internal-tide}
\end{figure}

Figure~\ref{fig:energetics-internal-tide} compares the volume-averaged kinetic energy of velocity fluctuations (TKE), reference potential energy (RPE), available potential energy (APE), and free-surface variance, computed as
\begin{align}
    \text{TKE} = \frac{1}{2V} \int_V \left({u^\prime}^2 + w^2 \right)\, \mathrm{d}V \ , \ \ 
    \text{RPE} = - \frac{1}{V}\int_V b z_{\text{res}}\, \mathrm{d}V  \ , \ \ 
    \text{APE} = \frac{1}{V}\int_V b (z_{\text{res}} - z)\, \mathrm{d}V \ , 
\end{align}
where $b$ is the buoyancy, $u^\prime \defn u - V^{-1} \int_V u \, \mathrm{d}V$, and $z_{\text{res}}$ is the re-sorted height computed as in the work of~\citeA{ilicak2016} and \citeA{hill2012}.
The RPE metric tracks integrated mixing in the domain: the larger the RPE increase, the greater the spurious diffusion introduced by the numerical discretization.
The largest difference between AB2-SE and the Runge--Kutta schemes lies in the energy of the flow: the TKE of AB2-SE saturates after about 20 days of evolution.

All the other discretizations show similar trajectories in TKE, APE and $\eta^2$. The RPE of AB2-SE grows twice as much as in the Runge--Kutta schemes, although less than in WRK3-UP, where RPE growth accelerates at the end of the 40 days as the higher TKE induces more mixing through the advection term. Normalized by TKE, AB2-SE shows the largest RPE growth of all the configurations, with irreversible mixing per unit kinetic energy about five times larger than the Runge--Kutta counterparts.

\begin{figure}[H]
    \centering
    \includegraphics[width=\linewidth,height=0.78\textheight,keepaspectratio]{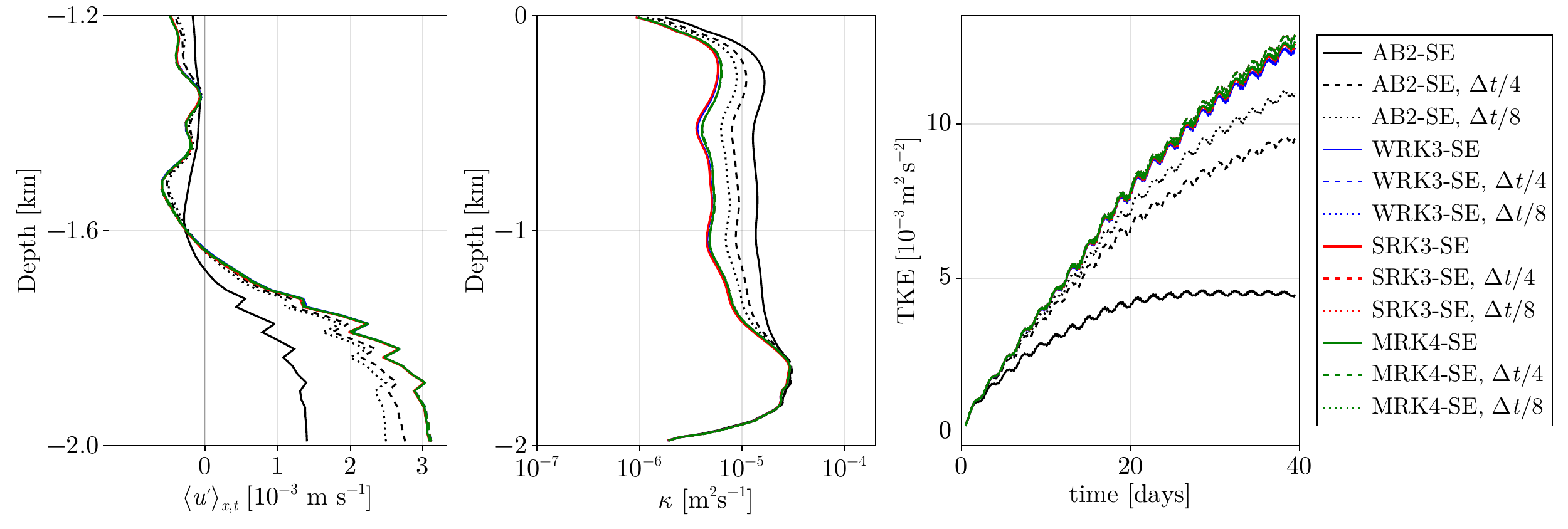}
    \caption{Temporal-resolution test of AB2-SE, WRK3-SE, SRK3-SE and MRK4-SE at the nominal timestep (solid) and at one fourth (dashed) and one eighth (dotted) of it. From left to right: velocity fluctuation averaged in the zonal direction and over the last ten days, time-averaged numerical diffusivity, and TKE.}
   \label{fig:rectified-internal-tide}
\end{figure} 

As a test of convergence we refined the timestep by factors of four and eight for AB2-SE, WRK3-SE, SRK3-SE and MRK4-SE. Figure~\ref{fig:rectified-internal-tide} shows that the mean velocity fluctuations, the numerical diffusivity and the TKE of the Runge--Kutta schemes do not change under time refinement, while AB2-SE approaches the Runge--Kutta solution as its timestep is reduced, confirming that both its lower energy and its larger mixing stem from the time discretization. Notably, AB2-SE at $\Delta t / 8 = 0.41$~min is still not as accurate as the Runge--Kutta configurations at $\Delta t = 11.3$~min.

\begin{figure}[H]
    \centering
    \includegraphics[width=\linewidth,height=0.78\textheight,keepaspectratio]{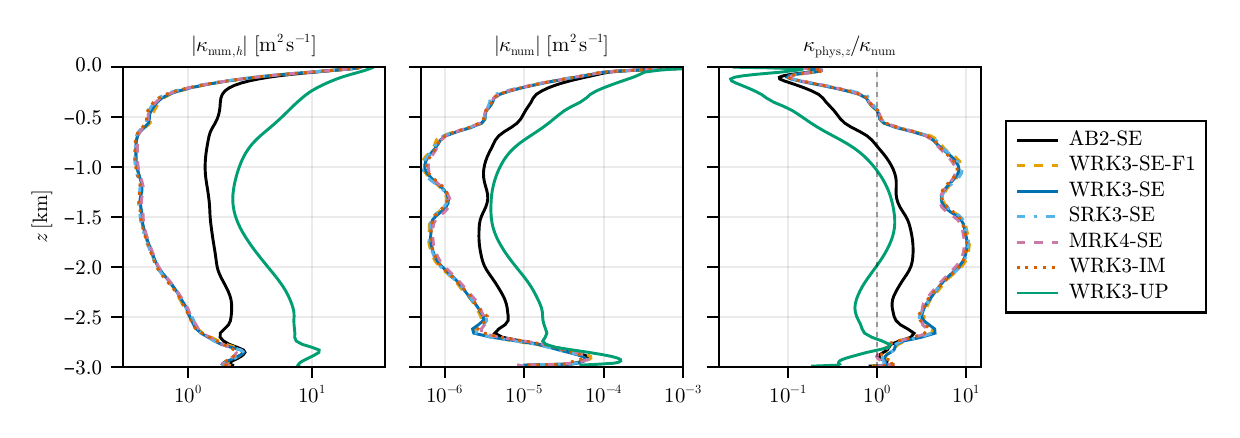}
    \caption{Numerical diffusivities for the channel configuration, averaged over years 35--40.
    From left to right: horizontal numerical diffusivity $- 0.5(\overline{P}_x + \overline{P}_y)/\overline{|\bnabla b|^2}$, total numerical diffusivity~\eqref{eq:numerical-diffusivity}, and the ratio
    of physical vertical diffusivity to total numerical diffusivity. The diffusivities are computed
    by averaging both the production and the variance over horizontal directions and time. The dashed vertical line in the right panel marks equal physical and numerical diffusivities.}
    \label{fig:dissipation-channel}
\end{figure}

\subsection{Idealized channel flow}
Our second test case is an idealized re-entrant channel with wind and buoyancy forcing and explicit vertical mixing following \citeA{hill2012}.
Unlike the previous transient configuration, this case relaxes toward statistical equilibrium, allowing us to assess how each time-stepping scheme affects the equilibrated solution and how numerical mixing interacts with a parameterized physical diffusivity.
The configuration consists of a 1000-kilometer-long, 2000-kilometer-wide, 3-kilometer-deep channel periodic in the $x$ direction and bounded in the $y$-direction. 
We use a 5-kilometer horizontal resolution and discretize the vertical direction with 90 non-uniformly spaced vertical layers, 10 m thick near the surface and the bottom, and 42 m thick in the interior.
To account for vertical physics, we use the CATKE vertical mixing parameterization \cite{catke} with a constant background viscosity $\nu = 10^{-4}\;\text{m}^2\,\text{s}^{-1}$ and no background diffusivity; no explicit horizontal viscosity or diffusivity is applied.
The channel is forced by a sinusoidal zonal wind stress and a surface buoyancy flux,
\begin{equation}
    \tau_x = -\tau_0 \sin\left(\dfrac{\pi y}{L_y}\right) \ , \quad Q_b = \begin{cases} Q_0 \cos\left(\dfrac{3\pi y}{L_y}\right) & y < \dfrac{5}{6}L_y \ , \\ 0 & \text{otherwise} \ , \end{cases}
\end{equation}
with $\tau_0 = 0.1\;\text{N}\,\text{m}^{-2}$ and $Q_0$ corresponding to a maximum heat flux of $10\;\text{W}\,\text{m}^{-2}$.
To obtain a non-trivial equilibrated solution, buoyancy is restored at the northern boundary within a sponge layer of width $L_y/10$, whose strength increases linearly toward the northern boundary, with a relaxation timescale of 7 days, toward the exponential profile
\begin{equation}
    b_{\rm targ}(z) = \Delta B \frac{e^{z/h} - e^{-L_z/h}}{1 - e^{-L_z/h}} \ ,
\end{equation}
where $\Delta B = 8\alpha g$, with $\alpha = 2\cdot 10^{-4}\;\text{K}^{-1}$, $h = 1000\;\text{m}$ and $L_z = 3000\;\text{m}$.
The Coriolis parameter varies on a beta-plane with $f_0 = -10^{-4} \text{s}^{-1}$ and a meridional gradient of $10^{-11} \text{m}^{-1}\text{s}^{-1}$. We apply a quadratic bottom drag with coefficient $1.1 \cdot 10^{-3}$.
More details about this setup can be found in \citeA{hill2012}.
We use a $\Delta t$ of $4.6$ min for AB2-SE, $16$ min for the RK3 discretizations and $25.8$ min for MRK4.
We initialize the buoyancy with an equilibrated solution generated with the MITgcm \cite{marshall1997} in the same configuration and, after a 100-day adjustment with $\Delta t = 1$ min, run for an additional 40 years. We then time-average the last 5 years (considered to be equilibrated to the new vertical mixing parameterization and numerical methods) and use it for our analyses below.

Figure~\ref{fig:dissipation-channel} shows the diagnosed numerical diffusivity profiles for the equilibrated state, with panels from left to right showing the horizontal diffusivity (only the $x+y$ directions), the total numerical diffusivity, and the ratio between the physical explicitly prescribed diffusivity and the total numerical diffusivity. For all the cases, the numerical diffusivity is $\sim1$–$10\;\text{m}^2\,\text{s}^{-1}$ in the horizontal direction, while the total numerical diffusivity varies from $10^{-6}$ to $10^{-4}\;\text{m}^2\,\text{s}^{-1}$. In terms of absolute values, the largest numerical diffusivity is observed in the WRK3-UP case, introduced by the comparatively lower-order advection scheme. 
The AB2-SE scheme has lower numerical mixing than WRK3-UP, but still higher than that of WRK3-SE (and the other Runge--Kutta configurations).
While the physical mixing is roughly the same across all cases (not shown), the ratio of physical to numerical mixing is different.
For WRK3-UP, $\kappa_{\rm phys} / \kappa_{\rm num} \lesssim 1$, indicating that the numerics exert a large control over the effective mixing in the whole water column.
The other cases show values larger than one, reaching local peaks $> 10$ in the interior of the domain, with the Runge--Kutta variants showing the largest values across the whole domain.

\begin{figure}[H]
    \centering
    \includegraphics[width=\linewidth,height=0.78\textheight,keepaspectratio]{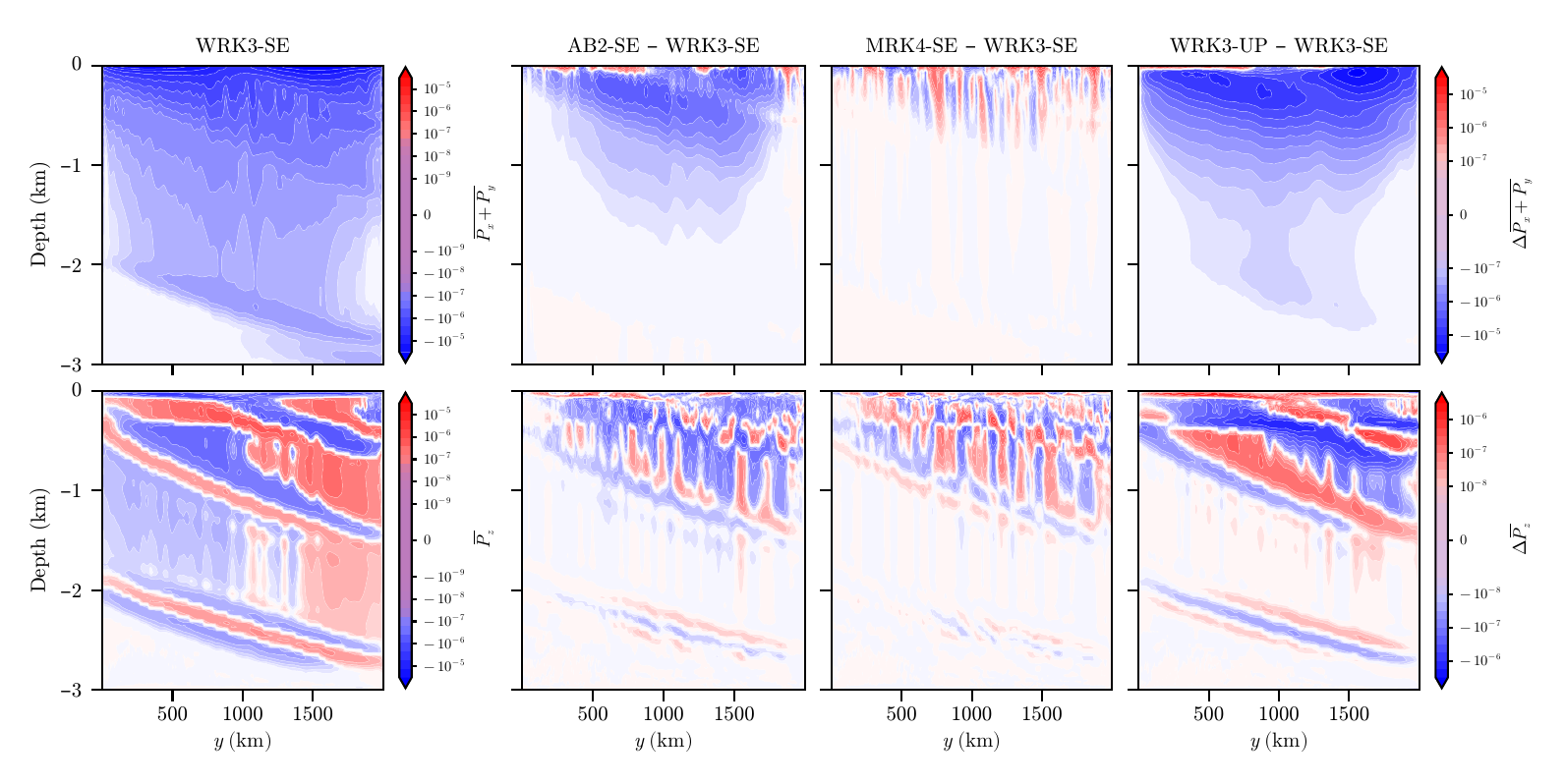}
    \caption{Zonal and time-mean numerical variance-budget terms for the channel configuration, averaged over years 35--40. The top row shows the horizontal contribution $P_x+P_y$ and the bottom row the vertical contribution $P_z$, in $\mathrm{m^5\,s^{-5}}$. The first column shows
    WRK3-SE; the remaining columns show each configuration minus WRK3-SE. Negative differences
    (blue) indicate a shift towards variance destruction, and positive differences (red) a shift
    towards variance production. Each row uses a common symmetric-logarithmic color scale for
    the differences, separate from the reference scale. No exponential depth weighting is applied.}
    \label{fig:dissipation-channel-contour}
\end{figure}

Figure~\ref{fig:dissipation-channel-contour} shows numerical dissipation averaged in the zonal direction split into horizontal and vertical dissipation for selected configurations. Positive values indicate regions of variance production, while negative values indicate variance destruction.
For all discretizations, unlike the internal tide test case, the horizontal contribution dominates the variance destruction, while the vertical contribution shows some anti-diffusive patches of variance production.
WRK3-UP shows by far the largest variance dissipation, confirming the large diffusivity values presented in Figure~\ref{fig:dissipation-channel}. 
AB2-SE also shows a larger dissipation in the horizontal direction when compared to WRK3-SE, explaining the much larger diffusivity in the left-most plot of Figure~\ref{fig:dissipation-channel}.
MRK4-SE shows a slightly lower variance destruction when compared to WRK3-SE.

\begin{figure}[H]
\centering
\includegraphics[width=\textwidth,height=0.78\textheight,keepaspectratio]{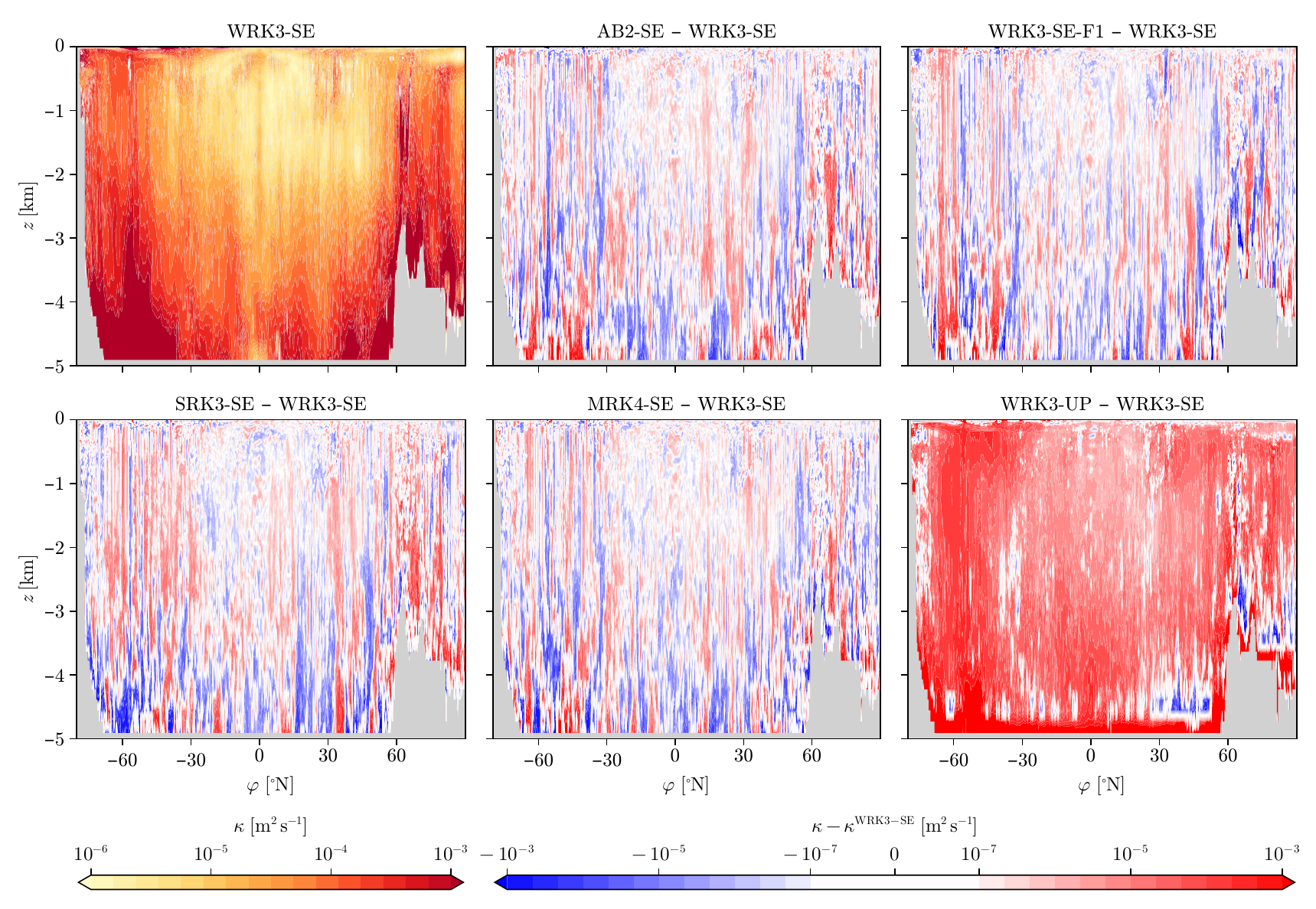}
\caption{Latitude--depth numerical diffusivity for WRK3-SE (upper left), and differences relative to WRK3-SE for the other configurations, over days 1080--1800. The reference uses a yellow--orange--red logarithmic scale from $10^{-6}$ to $10^{-3}\,\mathrm{m^2\,s^{-1}}$. Differences use a separate symmetric scale: blue is negative, white is near zero, and red is positive. The white band spans $\pm10^{-7}\,\mathrm{m^2\,s^{-1}}$; the endpoints are $\pm10^{-3}\,\mathrm{m^2\,s^{-1}}$. Values outside the displayed ranges use the endpoint colors. Gray denotes locations without valid diagnostics.}
\label{fig:global-section}
\end{figure}

\begin{figure}[H]
\centering
\includegraphics[width=\textwidth,height=0.78\textheight,keepaspectratio]{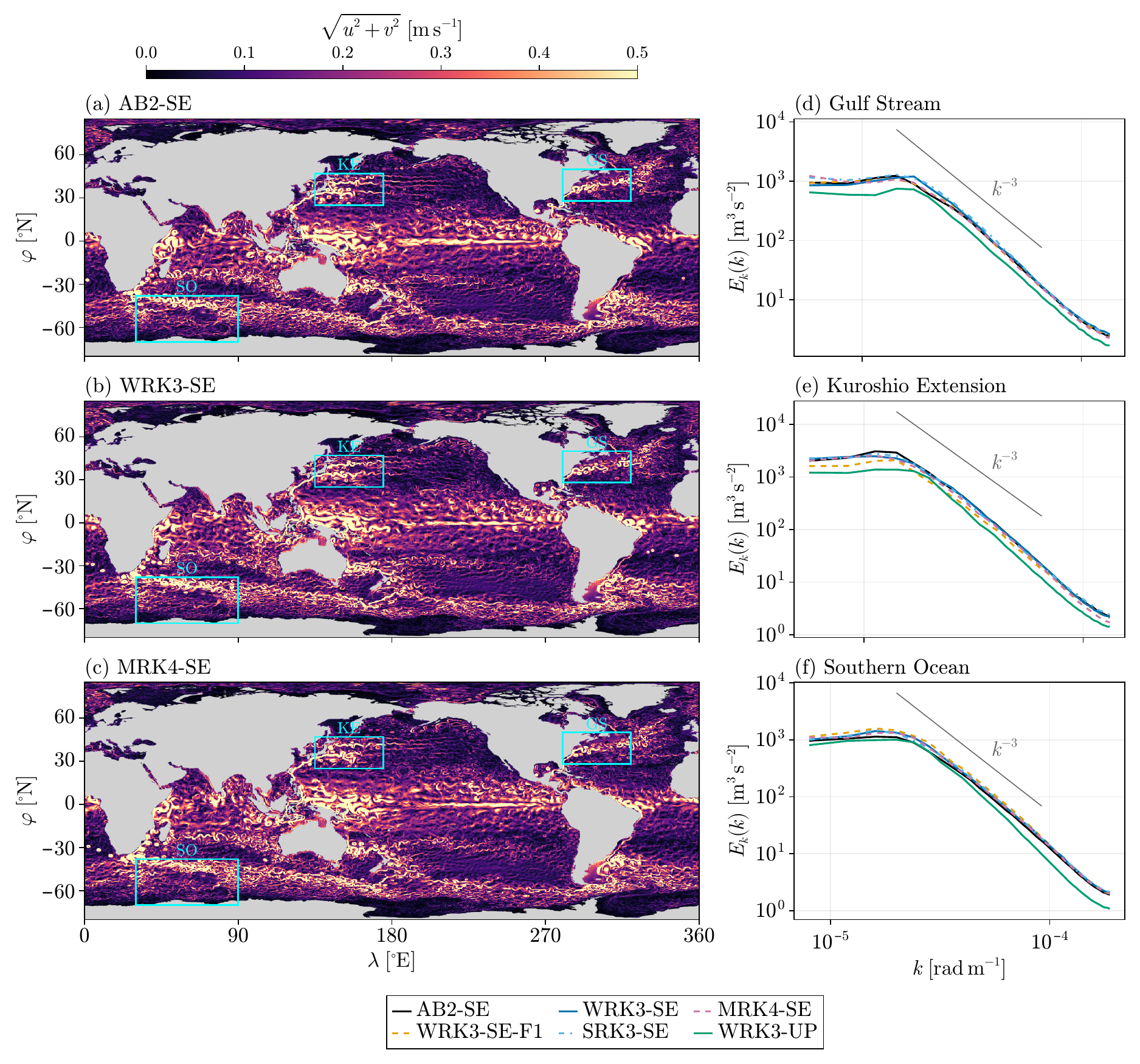}
\caption{Surface speed at day 1830 for AB2-SE, WRK3-SE, and MRK4-SE (left), and regional surface kinetic-energy spectra for all six configurations (right). Rectangles mark the Gulf Stream (GS), the Kuroshio Extension (KE) and the Southern Ocean (SO) regions used to compute the spectra.}
\label{fig:global-surface}
\end{figure}

\begin{table}[!t]
\centering
\caption{Measured cost of the quarter-degree global configurations on one Nvidia H100 GPU. SYPD denotes simulated years per wall-clock day, computed as $\Delta t/(365\,t_{\rm step})$ with both times in seconds. The relative throughput uses AB2-SE as the reference.}
\label{tab:global-cost}
\begin{tabular}{lrrrr}
\hline
Configuration & $\Delta t$ (min) & $t_{\rm step}$ (s) & SYPD & Relative throughput \\
\hline
AB2-SE      & 11 & 0.2534 & 7.136  & 1.00 \\
WRK3-SE     & 45 & 0.9097 & 8.132  & 1.14 \\
SRK3-SE     & 45 & 0.8873 & 8.337  & 1.17 \\
MRK4-SE     & 72 & 1.3859 & 8.540  & 1.20 \\
WRK3-SE-F1  & 45 & 0.8284 & 8.930  & 1.25 \\
WRK3-IM     & 45 & 2.1238 & 3.483  & 0.49 \\
WRK3-UP     & 45 & 0.7386 & 10.015 & 1.40 \\
\hline
\end{tabular}
\end{table}

\subsection{Quarter-degree configuration}
\label{sec:near-global}

We finally evaluate the proposed time discretizations in a quarter-degree global ocean configuration integrated for five years. The model uses the eORCA025 tripolar grid \cite{storkey2018}, with 100 vertical levels extending to 5000 m depth and a surface spacing of 5 m. Temperature and salinity are initialized from ECCO2 \cite{menemenlis2008}, while the surface forcing consists of an idealized latitude-dependent zonal wind stress and restoring of temperature and salinity towards monthly World Ocean Atlas climatologies \cite{locarnini2024,reagan2024} with a piston velocity of $1/6$ m day$^{-1}$. Salinity is restored through a virtual surface salt flux (no freshwater fluxes) and the configuration has no sea ice.
We use a linear equation of state to facilitate the computation of diapycnal mixing, CATKE vertical mixing, a quadratic bottom drag with coefficient $3\cdot 10^{-3}$, the WENO vector-invariant scheme for momentum, and ninth-order WENO reconstruction for tracer advection, with a separate third-order upwind comparison. Vertical advection of momentum and tracers is treated with an adaptive implicit scheme: the velocity exceeding a vertical Courant number of 1 (0.5 for AB2) is treated with implicit first-order upwinding.
The complete setup is described in the documentation available at \citeA{timesteppertestcases}.

The reduction in numerical mixing found in the idealized experiments does not occur in the global configuration. Over the complete budget window between days 1080 and 1800, the numerical diffusivities diagnosed from the time- and volume-integrated budgets remain close to AB2-SE for all the Runge--Kutta configurations except for WRK3-UP, which produces a widespread increase in numerical diffusivity (Figure~\ref{fig:global-section}), with a global value approximately two times larger than AB2-SE. We suspect this is because spatial discretization errors, in this relatively coarse, eddy-permitting setup, could dominate numerical mixing. At finer resolutions, the temporal contribution may become more apparent, as suggested by the idealized experiments. All tested configurations are robust and complete the five-year integrations.

Figure~\ref{fig:global-surface} shows surface-speed maps with major currents and eddy fields in the completed configurations as well as surface KE spectra, computed along grid rows and averaged in the last year of evolution over land-free sections in the regions shown in the surface speed maps. We use a Hann window to reduce boundary-adjacent artifacts. All cases have similar regional velocity spectra, with differences whose sign and magnitude depend on the region and wavelength. WRK3-UP has lower spectral energy across the displayed wavenumbers in all three regions, consistent with its larger numerical diffusivity. Together with the diffusivity diagnostics, these results show that the time discretizations maintain comparable surface variability in this quarter-degree configuration.

Table~\ref{tab:global-cost} reports the cost of the different time discretizations in this global setup. The Runge--Kutta timesteps follow the ratios of Table~\ref{tab:discretization} with $\Delta t = 45$~min for the three-stage schemes, and $\Delta t = 72$~min for MRK4, while AB2-SE runs with 11~min, below the 13~min implied by the ratio in Table~\ref{tab:discretization}. In this configuration, we empirically found AB2-SE to be unstable at 13~min, and 11~min is the largest stable timestep. The longer timesteps of the Runge--Kutta configurations more than compensate for their additional stage evaluations, and WRK3-SE, SRK3-SE, and MRK4-SE increase throughput by approximately 14\%, 17\% and 20\% relative to AB2-SE while costing 3.6, 3.5 and 5.5 times as much per step.
The barotropic treatment has a smaller effect: WRK3-SE-F1, which does roughly 4 times less barotropic work than WRK3-SE, is only about 10\% faster.
WRK3-IM, on the other hand, (not integrated for the full five years) reaches only half the throughput of AB2-SE, since the global elliptic solve of the free surface is expensive on the GPU.
We did not search for the largest stable timestep for the Runge--Kutta configurations as we did with AB2-SE, and we suspect we might be running at about 80\% of the stability limit, which could extend the computational savings even further. 

\section{Conclusions}
\label{sec:conclusions}

Baroclinic time-stepping in ocean models has long relied on low-order methods. A common choice is the over-extrapolated Adams--Bashforth scheme, which, while computationally efficient, damps waves across typical Courant numbers, even when those waves are well resolved. This dissipation contributes to spurious mixing that can degrade the solution, especially over long integrations, yet the contribution from temporal discretization has remained largely unquantified. 

As a lower-dissipation alternative, we developed here a low-storage Runge--Kutta framework that satisfies three essential constraints for free-surface ocean models: \emph{(i)}~consistency between free-surface evolution and grid geometry, \emph{(ii)}~tracer conservation and constancy preservation, and \emph{(iii)}~stability across coupled barotropic--baroclinic modes in a linearized setting. Both implicit and split--explicit free-surface treatments achieve these goals. We analyzed the coupling between barotropic and baroclinic modes through an eigenvalue perturbation analysis, separating direct contributions within each mode from indirect contributions through other modes. The analysis confirms the known stabilizing effect of an implicit free-surface treatment on the outer time discretization. For split--explicit schemes, on the other hand, stability depends on the barotropic solver, the averaging filter, and the reconstruction of the slow forcing. In the configurations considered, we find that \emph{(i)}~the dominant indirect contribution arises from coupling between barotropic and baroclinic modes, while contributions from coupling between baroclinic modes are negligible; \emph{(ii)}~damping of the barotropic mode is needed to suppress resonances in the indirect correction that arise when the uncoupled barotropic and baroclinic eigenvalues approach each other; \emph{(iii)}~reconstructing the slow forcing from previously computed stages reduces the coupling error associated with frozen forcing; and \emph{(iv)}~forward--backward barotropic solvers require second order in time damping to suppress the instability that arises at small timesteps, whereas the Runge--Kutta barotropic solvers considered here do not exhibit this failure mode. The extension of this analysis to an $f$-plane shows that the three-stage schemes remain stable up to $f\Delta t \simeq 0.7$, showing that rotation, rather than internal waves, may be the limiting factor in very coarse-resolution simulations or simulations of fast-spinning planets.

To demonstrate the benefits of this time discretization, we evaluated the proposed framework in two idealized test cases and a quarter-degree global configuration. In a nearly linear internal tide problem, the AB2 time discretization damps the flow to one third of the Runge--Kutta kinetic energy and, in the interior, mixes as much as a low-order spatial scheme. In an equilibrated re-entrant channel, Runge--Kutta schemes produce less variance destruction than AB2 throughout the water column, leading to lower numerical diffusivities and higher ratios of physical to numerical diffusivity. The reduction in numerical mixing does not occur in a five-year quarter-degree global configuration, where the numerical diffusivity remains close to the AB2 benchmark. We suspect that in this configuration, spatial discretization dominates mixing. Finally, we measured computational cost in the same quarter-degree configuration on an Nvidia H100 GPU. These measurements indicate that \emph{(i)}~the split--explicit barotropic configurations change throughput by at most 10\%; and \emph{(ii)}~the longer timesteps permitted by the outer Runge--Kutta schemes compensate for their additional stages, with a throughput increase of 14\% to 25\%.

The Runge--Kutta framework is implemented in Oceananigans \cite{oceananigans}, where it provides an alternative to Adams--Bashforth time-stepping. In the idealized configurations considered here, the split--explicit formulations reduce numerical mixing, while in the global configuration they leave it unchanged and increase computational throughput.
We recommend the SRK3-SE configuration, since it is stable in the linear analysis, remains stable on an $f$-plane over a wider range of $f\Delta t$ than WRK3, and preserves tracer bounds and positivity whenever the forward-Euler step does at the same timestep, which benefits monotone advection schemes and positive tracers such as biogeochemical ones. In our experiments, it matches the numerical mixing of WRK3-SE at a slightly higher throughput. Its frozen forcing gives a larger error than the progressive forcing of WRK3-SE and MRK4-SE in the linearized test, at the same order of accuracy, but this does not appear to affect the results of the more complex configurations of section~\ref{sec:5} which might be dominated by non-linear effects.
The eigenvalue perturbation analysis and exact tracer-variance budgets also provide complementary tools for assessing time discretization in other ocean models. Future work could extend the framework to non-hydrostatic dynamics and investigate the use of variance budgets to guide timestep selection. Our global experiments cover five years at eddy-permitting resolution. Longer integrations at mesoscale-resolving resolutions, including configurations coupled to biogeochemistry, are needed to assess how these improvements affect ocean circulation, tracer distributions, and climate simulations.

\section*{Acknowledgments}
We thank Milan Kl{\"o}wer for helpful discussions on the numerical precision of the split RK3 scheme and comments on an early draft of this paper.
This project is supported by Schmidt Sciences under the Climate Modeling Alliance project.
In addition, S.S.~is supported by the European Union's Horizon 2020 research and innovation program under the Marie Sk\l{}odowska-Curie grant agreement No.~101209201 and N.C.C.~is supported by the Australian Research Council under the Center of Excellence for the Weather of the 21st Century~CE230100012 and the Discovery Project~DP240101274.

\section*{Open Research Section}

The simulation outputs analyzed in this study, from the internal tide, the re-entrant channel and the quarter-degree global configuration, are archived on Zenodo under a CC BY 4.0 license \cite{timesteppertestcases-data}, with checksums, the final timings of Table~\ref{tab:global-cost} and a snapshot of the analysis source. The simulation setups, the variance-budget diagnostics, the linear stability analysis of Sections~2--4 and the scripts that produce every figure from the archived outputs are available in the TimestepperTestCases.jl repository \cite{timesteppertestcases} documented in the live documentation. AI tools (Claude Opus 5.5) were used to assist with the mathematical derivations, with launching the simulations, and with proofreading and correcting the manuscript. All AI-generated suggestions were thoroughly reviewed and verified by the authors who take full responsibility of the content of the manuscript, the mathematical derivations, and the configured experiments.

\appendix

\section{Averaging filters}
\label{sec:appendix-filters}

The first- and second-order filters belong to the low-dispersion family of~\citeA{shchepetkin}, which, with the exponents
used here, reads
\begin{equation}
    w(\tau) = \left(\frac{\tau}{\tau_0}\right)^{2} \left[1 - \left(\frac{\tau}{\tau_0}\right)^{4}\right] - r_0 \frac{\tau}{\tau_0} \ , \qquad \tau_0 = \frac{32}{21} \ ,
    \label{eq:sm05-filter}
\end{equation}
with $r_0 = 0.18927$ for F1, the default in Oceananigans, and $r_0 = 0.285$ for F2, the value that sets $M_2 = 0$. The third-order filter F3 is defined by
\begin{equation}
    w(\tau) = \begin{cases}
        \displaystyle (1-\hat \tau^2)^{\nu} \sum_{j=0}^{4} b_j \hat \tau^j \ , & |\hat \tau| \le 1 \ , \\[8pt]
        0 \ , & |\hat \tau| > 1 \ ,
    \end{cases}
    \qquad \hat \tau = \frac{\tau - 0.909165}{0.864972} \ , \qquad \nu = 0.595115 \ ,
\end{equation}
whose coefficients are the Euclidean minimum-norm solution of the conditions of normalization, centering and $M_2 = M_3 = 0$,
\begin{equation}
    (b_0, b_1, b_2, b_3, b_4) = (1.381694,\ 1.598408,\ -1.885242,\ -2.605757,\ -1.462849) \ .
\end{equation}
For $N_S$ substeps the shape is sampled at the fractional times $2s/N_S$ and truncated, after the last positive weight for
F1 and F2 and at the end of its support for F3, then normalized to unit sum, and the fractional substep $\Delta\tau$ is
rescaled from $2/N_S$ so that the centering of~\eqref{eq:filter-moments} holds exactly. The
coefficients of F3 and the coefficient $r_0$ of F2 are recomputed on the sampled times, so that the discrete moments satisfy
their conditions to round-off; the values above are the limit of the coefficients of F3 for large $N_S$, and $r_0 = 0.285$
that of F2.
The barotropic solver spans about $1.4\gamma_m$ for F1 and $1.76\gamma_m$ for F3.

\section{Assessing numerical mixing}
\label{sec:appendix-B}
Various methods to quantify numerical mixing have been developed in the ocean modeling literature \cite{klingbeil2014, hennel2023, banerjee2024, noemi}.
These methods were developed for single-step, single-stage methods and require adaptation to be used with the time discretization approaches we analyze in this paper.
In this section, we adapt the exact methodology described by \citeA{banerjee2024} to the AB2 scheme currently implemented in Oceananigans, and to the RK schemes developed in this paper.
For the restarted formulations (WRK3 and MRK4), only the last stage tendency contributes to the final tracer increment. For the blended formulation (SRK3), all three stage tendencies contribute with the weights of~\eqref{eq:ssp-weights}. The variance budget must use the effective tracer fluxes and volume transports of the complete timestep.

Therefore, the methodology described by \citeA{banerjee2024} reads
\begin{align}
    P_x &= 2 F^x_\star \delta_i \left( \frac{T^{n+1} + T^n}{2} \right) - U_\star \delta_i\left(T^{n+1}T^n\right) \ ,
    \label{eq:px-dissipation} \\ 
    P_y &= 2 F^y_\star \delta_j \left( \frac{T^{n+1} + T^n}{2} \right) - V_\star \delta_j\left(T^{n+1}T^n\right) \ ,
    \label{eq:py-dissipation} \\ 
    P_z &= 2 F^z_\star \delta_k \left( \frac{T^{n+1} + T^n}{2} \right) - W_\star \delta_k\left(T^{n+1}T^n\right)
    \label{eq:pz-dissipation} \ ,
\end{align}
where $\delta_d$ denotes a difference in the $d$ direction. 
For the restarted RK discretizations (WRK3 and MRK4):
\begin{alignat}{4}
    F^x_\star &= \mathcal{A}^x_{M-1} u_{M-1} \widetilde{T}_i^{M-1} \ \ \ && \text{and} \ \ \ U_\star &= && \mathcal{A}^x_{M-1} u_{M-1} \ , \\
    F^y_\star &= \mathcal{A}^y_{M-1} v_{M-1} \widetilde{T}_j^{M-1} \ \ \ && \text{and} \ \ \ V_\star &= && \mathcal{A}^y_{M-1} v_{M-1} \ , \\
    F^z_\star &= \mathcal{A}^z_{M-1} w_{M-1} \widetilde{T}_k^{M-1} \ \ \ && \text{and} \ \ \ W_\star &= && \mathcal{A}^z_{M-1} w_{M-1} \ , 
\end{alignat}
where the subscript indicates the stage ($M-1$), the velocities are the corrected velocities that advect the tracers, $\mathcal{A}$ are cell-areas, while $\widetilde{T}_i$, $\widetilde{T}_j$, and $\widetilde{T}_k$ indicate the finite volume reconstructions performed on $i$, $j$, and $k$ faces, respectively.

For SRK3, unrolling the tracer-content update~\eqref{eq:ssp-tracer} gives the effective fluxes
\begin{equation}
    F^d_\star = \sum_{m=1}^{3}\beta_m F^d_m \ , \qquad d \in \{x,y,z\} \ , \qquad
    (U_\star,V_\star,W_\star) = \sum_{m=1}^{3}\beta_m (U_m,V_m,W_m) \ .
    \label{eq:ssp-effective-fluxes}
\end{equation}
Here $F^d_m$ and $(U_m,V_m,W_m)$ are the tracer fluxes and volume transports used in the tracer update at stage $m$, including the stage-dependent geometry and the transport corrections of section~\ref{sec:3-ssp}.
These effective transports, and the corrected transports of the final stage for WRK3 and MRK4, are also those that model output should provide for offline tracer calculations.
For the AB2 time discretization, on the other hand, the effective fluxes and transports combine those of the two time-levels
\begin{equation}
F_\star^d = \left(\frac{3}{2} + \xi \right) F_n^d - \left(\frac{1}{2} + \xi \right) F_{n-1}^d \ , \quad (U_\star, V_\star, W_\star) = \left(\frac{3}{2} + \xi \right) (U_n, V_n, W_n) - \left(\frac{1}{2} + \xi \right) (U_{n-1}, V_{n-1}, W_{n-1}) \ , 
\end{equation}
where $F^d_n$ and $(U_n, V_n, W_n)$ are the fluxes and transports of the restarted case, with cell areas, velocities and reconstructions of time level $n$.
\begin{figure}[H]
    \centering
    \includegraphics[width=1\textwidth,height=0.78\textheight,keepaspectratio]{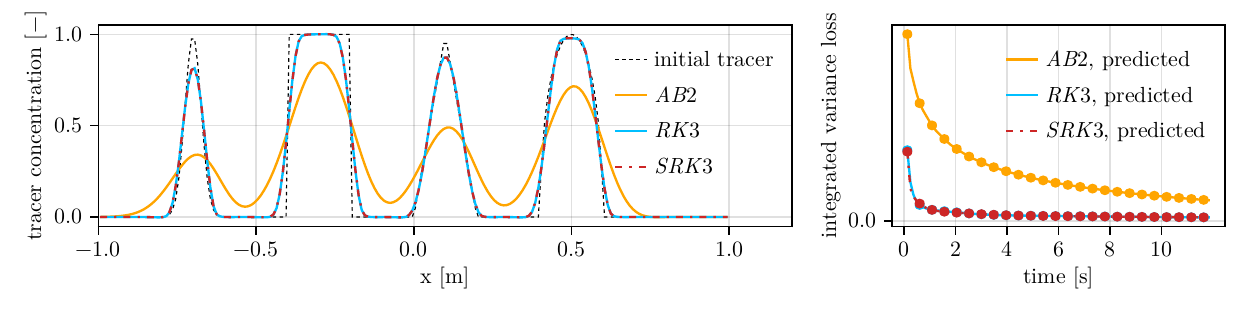}
    \caption{Validation of the tracer numerical mixing approach in a one-dimensional wave-advection test case. The left panel shows the advected tracer for AB2, WRK3 and SRK3 after 6 cycles compared to the initial distribution. The right panel shows the comparison between the integral variance destruction computed using the dissipation metric (lines) and equation~\eqref{eq:actual-dissipation} (markers).}
    \label{fig:dissipation-validation}
\end{figure}

Figure~\ref{fig:dissipation-validation} shows a validation of the above method applied to the WRK3, SRK3 and AB2 discretizations in the simple context of a one-dimensional tracer advection configuration using the WENO 7th order scheme for spatial reconstruction.
We advect an initial tracer profile consisting of various shapes by a constant velocity field $u = 1$.
We use a Courant number of 0.2 for AB2 and 0.5 for RK3.
The right panel compares the integrated tracer variance destruction, computed as
\begin{equation}
    \int_{-1}^1 \frac{\left( T^{n+1}\right)^2 - \left( T^{n}\right)^2} {\Delta t} \, \mathrm{d}x \ ,
    \label{eq:actual-dissipation}
\end{equation}
with the same quantity computed using the technique described in this section. The two approaches closely agree, verifying the implementation of the discrete variance budget.
We use this metric to assess the numerical mixing of buoyancy.
Since we use a linear equation of state that depends both on temperature~$T$ and salinity~$S$, we can compute diapycnal mixing, or mixing of buoyancy, by projecting the tracer fluxes in buoyancy space:
\begin{equation}
F_\star(b) = g \alpha_T  F_\star(T)  - g \beta_S F_\star(S) \ ,
\end{equation}
where $\alpha_T$ and $\beta_S$ are the (constant) thermal expansion and haline contraction coefficients.

\section{Coriolis treatment and stability}
\label{sec:rotation-results}

The Coriolis force acting on the depth-integrated flow contributes to the slow forcing of the barotropic solve, so that its treatment at the Runge--Kutta stages affects the stability of the coupled step. We extend the modal system of section~\ref{sec:4} to rotation on an $f$-plane and compare three forcing treatments for WRK3 and MRK4: frozen pressure and Coriolis, progressive pressure with frozen Coriolis, and progressive pressure and Coriolis. We also include SRK3 with frozen forcing. Supporting Information, section~\ref{si-sec:appendix-rotation}, gives the matrix construction.

\begin{figure}[H]
 \centering
 \includegraphics[width=\linewidth,height=0.78\textheight,keepaspectratio]{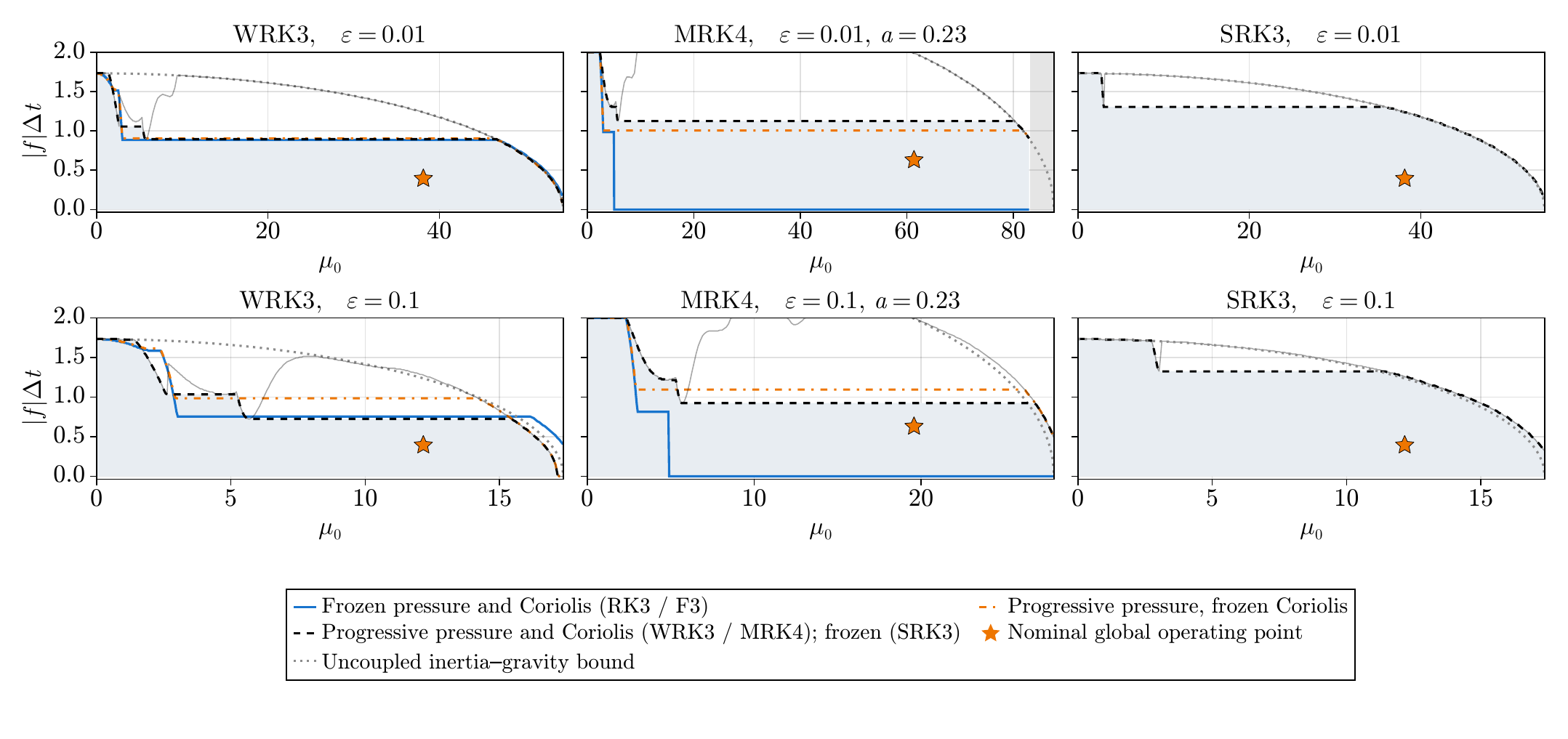}
\caption{Coupled stability boundaries for WRK3, MRK4, and SRK3. Shading marks stability for all $\mu \le \mu_0$ while stars indicate the global operating points. Thin gray curves show individual-wavenumber thresholds, and dotted curves show the uncoupled inertia--gravity bound. All cases retain eight baroclinic modes and use RK3 substeps, F3 filtering, and $N_S=96$. The gray band exceeds the fast-substep Courant limit. Frozen MRK4 is unstable near $\mu_0=5$ even without rotation.}
 \label{fig:rotation-coupling}
\end{figure}

With frozen forcing, the three-stage schemes advance a uniform inertial oscillation with the stability polynomial $R_3$ of section~\ref{sec:2}, which damps it up to $|f|\Delta t=\sqrt3$, for the depth-integrated and the baroclinic flow alike. At finite wavenumber, the uncoupled baroclinic modes are inertia--gravity waves, giving the joint condition
\begin{equation}
    \sqrt{(\alpha_1\mu_{0,\max})^2+(f\Delta t)^2}\leq\theta_{\max},
    \label{eq:rotation-joint-bound}
\end{equation}
where $\mu_{0,\max}$ is the largest resolved barotropic Courant number and $\theta_{\max}$ is the imaginary-axis stability limit of the outer scheme. Increasing the gravity-wave Courant number therefore reduces the available margin for rotation.

Figure~\ref{fig:rotation-coupling} includes this high-wavenumber restriction and the additional restriction from coupling. A model contains all resolved wavelengths, so the admissible timestep must satisfy the coupled boundary at every smaller $\mu_0$. This explains the horizontal portions of the black dashed curves: once an intermediate wavelength becomes unstable, extending the resolved wavenumber range cannot raise the limit. At $\varepsilon=0.1$, reconstructing pressure while holding Coriolis frozen gives a larger margin than reconstructing both for WRK3 and MRK4. Frozen MRK4 already has the nonrotating coupling instability discussed in section~\ref{sec:4}; progressive reconstruction removes it. At sufficiently large $\mu_0$, the inertia--gravity restriction becomes important for all schemes.
In the simulations of section~\ref{sec:5}, pressure, Coriolis and advection receive the same slow-forcing treatment: the progressive cases reconstruct them together, while SRK3-SE holds them at the preceding stage.
\section{Linear analysis of the coupled step}
\label{sec:short-linear-analysis}
This appendix describes how we construct the coupled evolution matrix and obtain the direct and indirect corrections used in section~\ref{sec:4}. The complete derivation is given in Supporting Information, section~\ref{si-sec:appendix-linear-analysis}. 

\subsection{From stage updates to the evolution matrix}
The modal equations are
\begin{equation}
 \partial_t\vect u=-ik\sqrt{gH}\,\vect h,\qquad
 \partial_t\vect h=-ik\sqrt{gH}\,\mathbf D\vect u,
 \qquad \mathbf D=\operatorname{diag}(\alpha_q^2).
\end{equation}
The modal weights satisfy $\chi_q=\psi_q\alpha_q^2$ and $\vect\Pi=\vect\chi-\vect\psi$. Thus $\bar u=\vect\chi^T\vect u$, $\eta=\vect\psi^T\vect h$, and the stage forcing are $\vect\Pi^T\vect h^m$. Supporting Information, section~\ref{si-sec:appendix-modes}, derives these identities and their small-$\varepsilon$ expansions. The analytical reconciliation uses $\vect\chi^T\vect\chi=1$, the complete-basis identity; its truncation error decreases with the number of retained modes.

For stage $m$ of duration $\gamma_m$, we use $\mu_m=k\sqrt{gH}\,\gamma_m$ and $\mu_0=k\sqrt{gH}\,\Delta t$. The continuity transport $\bar u^{\dagger,m}$ satisfies $\eta^m-\eta^n=-i\mu_m\bar u^{\dagger,m}$. For restarted stages, expanding the velocity reconciliation and the pressure update gives
\begin{equation}
\begin{aligned}
 \vect u^m&=\boldsymbol{\mathcal P}\vect u^n-i\mu_m\boldsymbol{\mathcal P}\vect h^{m-1}+\vect\chi\bar u^m,\\
 \vect h^m&=\vect h^n-i\mu_m\mathbf D\boldsymbol{\mathcal P}\vect u^{m-1}-i\mu_m\mathbf D\vect\chi\bar u^{\dagger,m},
 \qquad \boldsymbol{\mathcal P}=\mathbf I-\vect\chi\vect\chi^T.
\end{aligned}
\label{eq:short-stage}
\end{equation}
We write the forcing during stage $m$ as
\begin{equation}
 F^m(t')=\sum_{p=0}^{2}F_p^m t'^p,
 \label{eq:short-forcing-polynomial}
\end{equation}
where $t'=t-t^n$ and $F_p^m$ is the coefficient of $t'^p$; unused coefficients are zero. The coefficients follow from the stage values $\vect\Pi^T\vect h^j$ using the reconstruction in section~\ref{sec:3-mrk4}.
The barotropic outputs contain the prescribed substep solver, filter, and forcing sampling. The rows $(W_{11},W_{12})$ and $(\widetilde W_{11},\widetilde W_{12})$ map the unforced initial state $(\bar u^n,\eta^n)^T$ to the velocity and continuity transport, respectively. Linearity gives
\begin{equation}
\begin{aligned}
 \bar u^m&=W_{11}(\mu_m)\bar u^n+W_{12}(\mu_m)\eta^n+\sum_p F_p^m\bar u_p^m,\\
 \bar u^{\dagger,m}&=\widetilde W_{11}(\mu_m)\bar u^n+\widetilde W_{12}(\mu_m)\eta^n
                +\sum_p F_p^m\bar u_p^{\dagger,m}.
\end{aligned}
\label{eq:short-outputs}
\end{equation}
Here $\bar u_p^m$ and $\bar u_p^{\dagger,m}$ are the velocity and continuity-transport outputs with zero initial conditions and forcing $t'^p$. Frozen forcing has only $F_0^m=\vect\Pi^T\vect h^{m-1}$, with $\bar u_0^m=W_{12}(\mu_m)$ and $\bar u_0^{\dagger,m}=\widetilde W_{12}(\mu_m)$.

With $\vect y=(\vect u,\vect h)^T$, substituting~\eqref{eq:short-outputs} into~\eqref{eq:short-stage} and collecting the initial and preceding states gives $\vect y^m=\mathbf A_m\vect y^n+\mathbf B_m\vect y^{m-1}$. Writing $\vect y^m=\mathbf M_m\vect y^n$ and denoting the final stage by $M$, we obtain
\begin{equation}
 \mathbf M_0=\mathbf I,\qquad \mathbf M_m=\mathbf A_m+\mathbf B_m\mathbf M_{m-1},\qquad \mathbf M=\mathbf M_M.
\end{equation}
For reconstructed forcing, $F_p^m=\vect{\mathcal F}_p^m\vect y^n$, where $\vect{\mathcal F}_p^m$ is the row multiplying the initial state after substituting $\vect y^j=\mathbf M_j\vect y^n$ into the coefficient formulas. Its contribution therefore enters $\mathbf A_m$. For SRK3, each solve spans $\Delta t$ and the increment blend replaces the frozen-forcing matrices by
\begin{equation}
 \mathbf A_m^{\mathrm b}=(1-b_m)\mathbf I+b_m(\mathbf A_m-\mathbf I),\qquad
 \mathbf B_m^{\mathrm b}=b_m(\mathbf I+\mathbf B_m).
\end{equation}
Here $(b_1,b_2,b_3)=(1,1/4,2/3)$. Supporting Information, section~\ref{si-sec:app-stage-matrices}, gives the explicit blocks. The eigenvalues of $\mathbf M$ provide the full coupled stability test.

\subsection{From the evolution matrix to modal transfers}
The uncoupled matrix $\mathbf M_{\mathrm u}$ advances each mode independently: its barotropic block is $\mathbf W(\mu_0)$ and its baroclinic blocks apply the outer stability polynomial to the modal equations. Let $\lambda_j$ be its eigenvalues. We normalize its left and right eigenvectors by $\vect\ell_j\vect r_l=\delta_{jl}$ and define
\begin{equation}
 \mathbf V=\mathbf M-\mathbf M_{\mathrm u},\qquad
 V_{jl}=\vect\ell_j\mathbf V\vect r_l,\qquad
 \mathbf M\vect r_j=\lambda_j\vect r_j+\sum_l V_{lj}\vect r_l.
\label{eq:short-transfers}
\end{equation}
Taking $\vect y^n=\vect r_j$ selects one vertical mode and propagation direction. After one timestep, $V_{jj}$ changes its coefficient and $V_{bj}$ measures the generated barotropic mode. Taking $\vect y^n=\vect r_b$ instead gives the transfer $V_{jb}$ into mode $j$.

To evaluate these coefficients, subtract the uncoupled update from the coupled update with the same initial state. Denoting the resulting differences by $\delta u_q^M[\vect r_k]$ and $\delta h_q^M[\vect r_k]$ corresponding to an initial state at mode $k$, the projections are
\begin{equation}
 V_{jj}=\frac{1}{\sqrt2}
 	   \begin{bmatrix} 1 & \dfrac{\varsigma_j}{\alpha_q} \end{bmatrix}
	   \begin{bmatrix} \delta u_q^M[\vect r_j] \\ \delta h_q^M[\vect r_j] \end{bmatrix},  \quad
 V_{jb}=\frac{1}{\sqrt2}
 	   \begin{bmatrix} 1 & \dfrac{\varsigma_j}{\alpha_q} \end{bmatrix}
	   \begin{bmatrix} \delta u_q^M[\vect r_b] \\ \delta h_q^M[\vect r_b] \end{bmatrix},  \quad
 V_{bj}=\vect\ell_b\begin{bmatrix}\delta u_0^M[\vect r_j]\\\delta h_0^M[\vect r_j]\end{bmatrix} \ .
\end{equation}
Here $\varsigma_j=\pm1$ specifies the propagation direction. For input $\vect r_j$, the uncoupled update has zero barotropic amplitudes. The generated amplitudes therefore follow directly from the definitions of $\bar u^M$ and $\eta^M$:
\begin{equation}
 \delta u_0^M[\vect r_j]=\frac{\bar u^M-\sum_{r>0}\chi_ru_r^M}{\chi_0},\qquad
 \delta h_0^M[\vect r_j]=\frac{\eta^M-\sum_{r>0}\psi_rh_r^M}{\psi_0}.
\end{equation}
This specifies both the initial state and the projection for each transfer. Supporting Information, section~\ref{si-sec:app-transfers}, carries out the stage substitutions, first with general barotropic outputs and then for frozen forcing and the SRK3 blend.

\subsection{From modal transfers to growth}
For an isolated eigenvalue, first-order perturbation gives $\delta\lambda_j^{(1)}=V_{jj}$ and an eigenvector contribution $V_{lj}/(\lambda_j-\lambda_l)$ along $\vect r_l$. Applying $\mathbf V$ to that contribution and projecting back onto $\vect\ell_j$ gives the second-order correction:
\begin{equation}
 \delta\lambda_j^{(1)}=V_{jj},\qquad
 \delta\lambda_j^{(2)}=\sum_{l\ne j}\frac{V_{jl}V_{lj}}{\lambda_j-\lambda_l}.
\end{equation}
The first term is the direct correction. Each indirect term contains a transfer from $j$ to $l$ and back; a small eigenvalue separation increases its contribution. Section~\ref{sec:4} compares the full sum with the two barotropic contributions, for which $l=b$. Writing their combined eigenvalue correction as $\delta\lambda_j$, its contribution to the amplification magnitude is
\begin{equation}
 |\lambda_j+\delta\lambda_j|-|\lambda_j|
 =\frac{\operatorname{Re}(\overline{\lambda_j}\delta\lambda_j)}{|\lambda_j|}
 +\mathcal O\!\left(\frac{|\delta\lambda_j|^2}{|\lambda_j|}\right).
\end{equation}
Positive corrections oppose the uncoupled damping $1-|\lambda_j|$. The expansion requires separated eigenvalues; the full matrix remains the stability test near a coincidence. Supporting Information, section~\ref{si-sec:appendix-perturbation}, gives the perturbation derivation.

\subsection{Coriolis forcing in the split--explicit stages}
\label{sec:rotation-matrix-summary}
On an $f$-plane, taking the horizontal Fourier wavevector along $x$ gives the rotating modal system
\begin{equation}
    \frac{d\vect u}{dt}=f\vect v-ik\sqrt{gH}\vect h,\qquad
    \frac{d\vect v}{dt}=-f\vect u,\qquad
    \frac{d\vect h}{dt}=-ik\sqrt{gH}\mathbf D\vect u.
    \label{eq:rotation-summary-system}
\end{equation}
Here $f$ is the constant Coriolis parameter, and $\vect u$ and $\vect v$ contain the velocity amplitudes parallel and perpendicular to the wavevector. Coriolis couples the two velocities, while the pressure gradient acts along the wavevector.
We use $\vect y=(\vect u,\vect v,\vect h)^T$, with $\mathbf E_u=[\mathbf I\ \mathbf0\ \mathbf0]$, $\mathbf E_v=[\mathbf0\ \mathbf I\ \mathbf0]$, and $\mathbf E_h=[\mathbf0\ \mathbf0\ \mathbf I]$ selecting its three blocks. The filtered barotropic outputs are linear functions of the initial state:
\begin{equation}
    \bar u^m=\mathbf X_{m,1}\vect y^n,\qquad
    \bar v^m=\mathbf X_{m,2}\vect y^n,\qquad
    \eta^m=\mathbf X_{m,3}\vect y^n.
    \label{eq:rotation-summary-outputs}
\end{equation}
The rows $\mathbf X_{m,1}$, $\mathbf X_{m,2}$, and $\mathbf X_{m,3}$ include the chosen barotropic solver, averaging filter, and forcing reconstruction. Supporting Information, section~\ref{si-sec:appendix-rotation}, gives their construction and the complete derivation of the stage matrix. Using continuity to express the transport through $\eta^m-\eta^n$, substitution into the velocity reconciliation and pressure update gives
\begin{equation}
\mathbf M_m=\begin{bmatrix}
 \boldsymbol{\mathcal P}(\mathbf E_u+f\gamma_m\mathbf E_v\mathbf M_{m-1}-i\mu_m\mathbf E_h\mathbf M_{m-1})+\vect\chi\mathbf X_{m,1}\\[3pt]
 \boldsymbol{\mathcal P}(\mathbf E_v-f\gamma_m\mathbf E_u\mathbf M_{m-1})+\vect\chi\mathbf X_{m,2}\\[3pt]
 \mathbf E_h-i\mu_m\mathbf D\boldsymbol{\mathcal P}\mathbf E_u\mathbf M_{m-1}
       +\mathbf D\vect\chi(\mathbf X_{m,3}-\vect\psi^T\mathbf E_h)
\end{bmatrix}.
\label{eq:rotation-summary-stage}
\end{equation}
Here $\boldsymbol{\mathcal P}=\mathbf I-\vect\chi\vect\chi^T$ as above, and $\mathbf M_0=\mathbf I$. The forcing is either frozen, reconstructed for pressure only, or reconstructed for both pressure and Coriolis. The stage durations are those of the outer scheme.

For SRK3, each solve spans $\Delta t$ and uses the preceding stage forcing. After evaluating~\eqref{eq:rotation-summary-stage}, we apply the increment blend:
\begin{equation}
    \mathbf M_m \longleftarrow
    (1-b_m)\mathbf I
    +b_m\left(\mathbf M_{m-1}+\mathbf M_m-\mathbf I\right),
    \qquad (b_1,b_2,b_3)=(1,1/4,2/3).
    \label{eq:rotation-summary-blend}
\end{equation}
Applying these stages successively gives the full-step matrix used for the stability boundaries in Appendix~\ref{sec:rotation-results}.

\clearpage
\section*{Supporting Information}
\setcounter{section}{0}
\setcounter{equation}{0}
\renewcommand{\thesection}{S\arabic{section}}
\renewcommand{\theequation}{S\arabic{equation}}
\renewcommand{\theHsection}{supplement.\arabic{section}}
\renewcommand{\theHequation}{supplement.\arabic{equation}}

\section{Linear analysis of the coupled step}
\label{si-sec:appendix-linear-analysis}

This supporting information gives the derivations used in section~\ref{sec:4} to assess the stability of the proposed discretizations in a linear setting. We follow the linear-analysis framework of \citeA{demange} and expand the resting, stratified system in vertical modes, applying the numerical stages of section~\ref{sec:3} to their amplitudes.
We first derive the modal equations and assemble the evolution matrices, then we separate the uncoupled amplification factors from the corrections introduced by the coupling, and derive how the stage formulation, the time filters, the substep solvers and the forcing reconstructions affect damping. 
\subsection{Vertical modes}
\label{si-sec:appendix-modes}

Eliminating $u$, $w$, and $b$ from the linear system in section~\ref{sec:4} for a single wavenumber $k$, gives the associated Sturm--Liouville problem.
Combining the hydrostatic and buoyancy equations, the linear system reads
\begin{equation}
\hat{w} = - \frac{1}{N^2} \frac{\partial}{\partial z}\frac{\partial}{\partial t}(g\hat{\eta} +\hat{p}) \ , \qquad   \frac{\partial \hat{u}}{\partial t} = - ik (g\hat{\eta} +\hat{p}) \ ,  \qquad ik \hat{u} = - \frac{\partial\hat{w}}{\partial z} \ .
\label{si-eq:app-pressure-elimination}
\end{equation}
The boundary conditions are $\hat{p}(0) = 0$ and $\hat{w}(0) = \partial_t\hat{\eta}$ at the surface and $\hat{w}(-H) = 0$ at the bottom. Separating $g\hat\eta+\hat p$ into time-dependent amplitudes and vertical modes $\phi_q(z)$, as done in section~\ref{sec:4}, and using the boundary conditions, yields the eigenvalue problem
\begin{equation}
    -\frac{\mathrm{d}}{\mathrm{d}z}\Big(\frac{1}{N^{2}}\frac{\mathrm{d}\phi_q}{\mathrm{d}z}\Big) = \frac{\phi_q}{c_q^{2}} \ , \qquad
    \frac{\mathrm{d}\phi_q}{\mathrm{d}z}\Big|_{0} = -\frac{N^{2}}{g}\,\psi_q \ , \qquad
    \frac{\mathrm{d}\phi_q}{\mathrm{d}z}\Big|_{-H} = 0 \ .
    \label{si-eq:app-sl}
\end{equation}
For constant $N^2$, the solutions are the ones in equations~\eqref{eq:modes}. Normalization gives
\begin{equation}
\frac{1}{A_q^2} = \frac{1}{H} \int_{-H}^0 \cos^2\left[x_q\left(1+\frac{z}{H}\right)\right] dz = \frac{1}{2} + \frac{\sin{2x_q}}{4x_q} \ , \quad A_q > 0 \ .
\label{si-eq:app-modes}
\end{equation}
Integrating~\eqref{si-eq:app-sl} over the depth and using the boundary conditions gives
\begin{equation}
   - \left[\frac{1}{N^2}\frac{d\phi_q}{dz}\right]^0_{-H} =  \frac{\psi_q}{g} = \frac{1}{c_q^2}\int_{-H}^0\phi_q dz = \frac{H\,\chi_q}{c_q^{2}} \ ,
 \qquad
    \chi_q = \psi_q\,\alpha_q^{2} \ ,
    \label{si-eq:app-chi}
\end{equation}
where the last equality uses $\alpha_q^2 = c_q^2/(gH)$. For the barotropic root, expanding $x_0\tan x_0=x_0^2 + x_0^4/3 + 2x_0^6/15...$ 
gives $x_0^2 = \varepsilon - \varepsilon^2/3 + 4\varepsilon^3 /45 +...$, while for $q\ge1$ expansion about $q\pi$ gives $x_q = q\pi + \varepsilon / (q\pi) + \mathcal{O}(\varepsilon^2)$.
Inserting these roots into $\alpha_q = \sqrt{\varepsilon}/x_q$, $\psi_q = A_q\cos{x_q}$ and $\chi_q = A_q \sin{x_q}/x_q$ with~\eqref{si-eq:app-modes} gives
\begin{equation}
\begin{aligned}
    \alpha_0^2 & = 1 + \frac{\varepsilon}{3} + \mathcal{O}(\varepsilon^2) \ ,
    \qquad
    && \psi_0 = 1 - \frac{\varepsilon}{3} + \mathcal{O}(\varepsilon^2) \ , && 
    \chi_0 = 1 - \frac{\varepsilon^2}{90} + \mathcal{O}(\varepsilon^3) \ , \\ 
    \qquad
    \alpha_{q\ge1} & = \frac{\sqrt{\varepsilon}}{q\pi} + \mathcal{O}(\varepsilon^{3/2}) \ ,
    \qquad
    && \psi_{q\ge1} = \sqrt{2}(-1)^q + \mathcal{O}(\varepsilon)  
    \qquad &&\chi_{q\ge 1} = (-1)^{q}\,\frac{\sqrt{2}\,\varepsilon}{q^{2}\pi^{2}} + \mathcal{O}(\varepsilon^{2}) \ .
    \end{aligned}
    \label{si-eq:mode-expansions}
\end{equation}
Therefore, the barotropic mode travels slightly faster than $\sqrt{gH}$ and has a surface value slightly below one.
By orthonormality, the coefficient of the constant function along $\phi_q$ is $\chi_q$. Expanding that function and taking its squared norm gives
\begin{equation}
    1 = \sum_q\chi_q\,\phi_q \ , \qquad
    \sum_q\chi_q^{2} = 1 \ ,
    \label{si-eq:app-reads}
\end{equation}
The identities~\eqref{si-eq:app-reads} allow the velocity reconciliation to impose $\vect\chi^T\vect u^m=\bar u^m$. We use this substitution in the analytical derivation although truncating the modal sums introduces an error that decreases as more modes are retained.

Here $u_q$ and $h_q$ are the dimensionless modal amplitudes of horizontal velocity and total pressure, respectively, as defined in section~\ref{sec:4}.
Their depth averages are obtained by projecting the corresponding modal amplitude vectors onto $\vect{\chi}$.
The momentum equation directly gives an equation for $u_q$. To obtain an equation for $h_q$, insert the modal expansion into~\eqref{si-eq:app-pressure-elimination} and use~\eqref{si-eq:app-sl}
\begin{equation}
ik \sqrt{gH} \sum_q u_q \phi_q = - gH \sum_q \frac{\partial_t h_q}{c_q^2}\phi_q \ .
\label{si-eq:app-modal-continuity}
\end{equation}
Equating the coefficient for each $\phi_q$ yields
\begin{equation}
    \pd{u_q}{t} = -ik\sqrt{gH}\,h_q \ , \qquad
    \pd{h_q}{t} = -ik\sqrt{gH}\,\alpha_q^2 u_q \ ,
    \label{si-eq:app-rotation}
\end{equation}
These equations define the uncoupled modal evolution used in section~\ref{si-sec:appendix-perturbation}.

\subsection{Stage matrices}
\label{si-sec:app-stage-matrices}

Using the solution vector $\vect{y} = (\vect{u}^T, \vect{h}^T)^T$, all the proposed discretizations can be expressed as $\vect{y}^m = \mathbf{A}_m \vect{y}^n + \mathbf{B}_m \vect{y}^{m-1}$, so the evolution matrix satisfies
\begin{equation}
\begin{cases}
\mathbf{M}_0 & = \mathbf{I} \ , \\ 
\mathbf{M}_m & = \mathbf{A}_m + \mathbf{B}_m \mathbf{M}_{m-1} \ , \\
\mathbf{M} & = \mathbf{M}_M \ .
\end{cases}
\end{equation}
Linearizing the restarted stage of section~\ref{sec:3-stage} about rest, with unit specific thickness ($\sigma=1$), gives
\begin{equation}
    \begin{aligned}
    \vect{u}^{*} &= \vect{u}^{n} - i\mu_m\vect{h}^{m-1} \ , \qquad &
    & \vect{u}^{m} = \vect{u}^{*} + \vect{\chi}\big(\bar u^{m} - \vect{\chi}^{T}\vect{u}^{*}\big) \ , \\
    \vect{h}^{m} &= \vect{h}^{n} - i\mu_m\mathbf{D}\,\vect{u}^{\dagger} \ , \qquad &
    & \vect{u}^{\dagger} = \vect{u}^{m-1} + \vect{\chi}\big(\bar u^{\dagger, m} - \vect{\chi}^{T}\vect{u}^{m-1}\big) \ ,
    \end{aligned}
    \qquad \mathbf{D} = \mathrm{diag}\big(\alpha_q^{2}\big) \ ,
    \label{si-eq:app-stage}
\end{equation}
where the reconciliation~\eqref{eq:stage-reconcile} replaces the depth average of the predictor with the barotropic velocity and the transport~\eqref{eq:transport-evolution} replaces the average of the tracer velocity in the $\vect{h}$ equations.

Since the barotropic equations (in this framework) are linear, we compute separately the outputs of the barotropic solve due to the initial state and those due to the forcing. Following section~\ref{sec:3-mrk4}, we write the forcing during stage $m$ as
\begin{equation}
    F^m(t^\prime) = \sum_p F^m_p t^{\prime p} \ ,
\end{equation}
where $F^m_p$ denotes the coefficient of $t^{\prime p}$ corresponding to the $\mathbf{G}_p^m$ coefficients introduced in section~\ref{sec:3-mrk4}. For a frozen forcing, $F^m_p = 0$ for $p > 0$ and $F_0^m = \vect{\Pi}^T\vect{h}^{m-1}$. We denote by $\bar u^m_p$ and $\bar u^{\dagger, m}_p$ the velocity and continuity-transport outputs for forcing $t^{\prime p}$ with zero initial conditions, using the prescribed solver, filter, and forcing sampling over stage $m$.
The forcing contributions are then
\begin{equation}
    \bar u^m_{\text{forced}} = \sum_p F^m_p \bar u^m_p \ ,  \qquad \bar u^{\dagger, m}_{\text{forced}} = \sum_p F^m_p \bar u^{\dagger, m}_p \ .
\end{equation}
Adding the initial-state and forcing contributions gives the barotropic outputs, with $W_{ij}$ and $\widetilde W_{ij}$ evaluated at $\mu_m$:
\begin{equation}
    \begin{aligned}
        \bar u ^m & = W_{11} \vect{\chi}^T \vect{u}^n + W_{12} \vect{\psi}^T\vect{h}^n + \sum_p F^m_p \bar u^m_p \\
        \bar u ^{\dagger, m} & = \widetilde W_{11} \vect{\chi}^T \vect{u}^n +  \widetilde W_{12} \vect{\psi}^T\vect{h}^n + \sum_p F^m_p \bar u^{\dagger, m}_p  
    \end{aligned}
    \label{si-eq:app-barotropic-complete}
\end{equation}
For frozen forcing $dF^m_0/dt^\prime = 0$, so $(\bar u, \eta + F_0^m)$ satisfy the unforced barotropic equations of section~\ref{sec:4}. With zero initial state and unit forcing, the initial state is $(0, 1)^T$. Its velocity and transport outputs are the second-column coefficients
\begin{equation}
\bar u^m_0 = W_{12}(\mu_m) \ ,  \qquad \bar u_0^{\dagger, m} = \widetilde W_{12}(\mu_m)
\label{si-eq:app-unit-constant-forcing}
\end{equation}
Multiplying by $F^m_0 = \vect{\Pi}^T\vect{h}^{m-1}$ in~\eqref{si-eq:app-barotropic-complete} gives
\begin{equation}
    \begin{aligned}
        \bar u ^m & = W_{11} \vect{\chi}^T \vect{u}^n + W_{12} (\vect{\psi}^T\vect{h}^n + \vect{\Pi}^T \vect{h}^{m-1}) \\
        \bar u ^{\dagger, m} & = \widetilde W_{11} \vect{\chi}^T \vect{u}^n +  \widetilde W_{12} (\vect{\psi}^T\vect{h}^n + \vect{\Pi}^T \vect{h}^{m-1}) 
    \end{aligned}
    \label{si-eq:app-barotropic}
\end{equation}
Here $W_{ij}$ and $\widetilde W_{ij}$ are evaluated at $\mu_m$. We introduce $\boldsymbol{\mathcal{P}} = \mathbf{I} - \vect{\chi}\vect{\chi}^T$, the operator that subtracts the depth average from a variable (i.e. $\boldsymbol{\mathcal{P}}\vect{\chi} = 0$ in the complete modal basis), and expand the reconciliation and pressure update in~\eqref{si-eq:app-stage} yielding
\begin{equation}
\begin{aligned}
\vect{u}^m & = \boldsymbol{\mathcal{P}} \vect{u}^n - i \mu_m \boldsymbol{\mathcal{P}} \vect{h}^{m-1} + \vect{\chi}\bar u^m \ , \\
\vect{h}^m & = \vect{h}^n - i\mu_m \mathbf{D} \boldsymbol{\mathcal{P}} \vect{u}^{m-1} - i\mu_m \mathbf{D} \vect{\chi} \bar{u}^{\dagger, m} \ .
\label{si-eq:app-stage-collected}
\end{aligned}
\end{equation}
Substituting~\eqref{si-eq:app-barotropic} into the last term of each row and collecting the coefficients of $(\vect{u}^n, \vect{h}^n)$ into $\mathbf{A}_m$ and the ones of $(\vect{u}^{m-1}, \vect{h}^{m-1})$ into $\mathbf{B}_m$ gives the stage matrices
\begin{equation}
    \mathbf{A}_m = \begin{bmatrix}
        \vect{\mathcal{P}} + W_{11}\vect{\chi}\vect{\chi}^{T} & W_{12}\,\vect{\chi}\vect{\psi}^{T} \\[2pt]
        -i\mu_m\widetilde W_{11}\,\mathbf{D}\vect{\chi}\vect{\chi}^{T} & \mathbf{I} - i\mu_m\widetilde W_{12}\,\mathbf{D}\vect{\chi}\vect{\psi}^{T}
    \end{bmatrix} , \quad
    \mathbf{B}_m = \begin{bmatrix}
        \mathbf{0} & -i\mu_m\vect{\mathcal{P}} + W_{12}\,\vect{\chi}\vect{\Pi}^{T} \\[2pt]
        -i\mu_m\mathbf{D}\vect{\mathcal{P}} & -i\mu_m\widetilde W_{12}\,\mathbf{D}\vect{\chi}\vect{\Pi}^{T}
    \end{bmatrix} .
    \label{si-eq:app-AB}
\end{equation}
Note that these matrices (for frozen forcing), once we consider the relation between the transport weights $\widetilde{W}_{ij}$ and the barotropic weights $W_{ij}$ reduce to the ones derived in \citeA{demange}.

For a polynomial forcing, the coefficients $F^m_p$ are obtained using the reconstructions in section~\ref{sec:3-mrk4}, replacing each $\overline{\mathbf{G}}{}^j$ with $\vect{\Pi}^T\vect{h}^j$. To construct the evolution matrices, we use $\vect{y}^m = \mathbf{M}_m \vect{y}^n$. For example, the linear coefficient at the second stage $F^2_1$ is
\begin{equation}
    F^2_1 = \frac{\vect{\Pi}^T(\vect{h}^1 - \vect{h}^n)}{c_1 \Delta t} = \underbrace{\frac{1}{c_1 \Delta t} \begin{bmatrix}
        \vect{0}^T & \vect{\Pi}^T 
    \end{bmatrix}(\mathbf{M}_1 - \mathbf{I})}_{\vect{\mathcal{F}}^2_1} \vect{y}^n \ .
\end{equation}
We can express all the remaining coefficients $F^m_p = \vect{\mathcal{F}}^m_p\vect{y}^n$ where $\vect{\mathcal{F}}^m_p$ contains the dependence on the previously computed stage matrices. 
Since each $\vect{\mathcal F}_p^m$ multiplies $\vect y^n$, its forcing contribution enters $\mathbf A_m$. Removing the frozen-forcing terms from $\mathbf B_m$ gives
\begin{equation}
    \mathbf{A}_m \longrightarrow
    \mathbf{A}_m +
    \begin{bmatrix}
        \vect{\chi}\displaystyle\sum_p \bar{u}^m_p\vect{\mathcal F}_p^{\,m} \\[3pt]
        -i\mu_m\mathbf{D}\vect{\chi}
        \sum_p  \bar{u}^{\dagger, m}_p\vect{\mathcal F}_p^{\,m}
    \end{bmatrix},
    \qquad
    \mathbf{B}_m \longrightarrow
    \begin{bmatrix}
        \mathbf{0} & -i\mu_m\boldsymbol{\mathcal P} \\
        -i\mu_m\mathbf{D}\boldsymbol{\mathcal P} & \mathbf{0}
    \end{bmatrix}.
    \label{si-eq:app-AB-reconstructed}
\end{equation}

For SRK3, we start from the frozen forcing matrices~\eqref{si-eq:app-AB}. Each barotropic solve spans the full timestep and the $b_m$ weights are applied to the state and the barotropic increment $\Delta \bar u ^m$ giving
\begin{equation}
    \bar u^{m} = \big(1 - b_m\big)\bar u^{n} + b_m\big(\bar u^{m-1} + \Delta\bar u^{m}\big)  \ .
\end{equation}
With $\vect{\chi}^T\vect{\chi} = 1$, expanding the blended velocity reconciliation and pressure update gives
\begin{equation}
 \vect{y}^m = (1-b_m) \vect{y}^n + b_m [\vect{y}^{m-1} + (\mathbf{A}_m - \mathbf{I})\vect{y}^n + \mathbf{B}_m\vect{y}^{m-1}]
 \label{si-eq:app-blended-matrix-substitution}
\end{equation}
Collecting the two states gives
\begin{equation}
    \begin{gathered}
    \mathbf{A}^{\mathrm{b}}_m = \big(1 - b_m\big)\mathbf{I} + b_m\big(\mathbf{A}_m - \mathbf{I}\big) \ , \qquad
    \mathbf{B}^{\mathrm{b}}_m = b_m\big(\mathbf{I} + \mathbf{B}_m\big) \ ,
    \end{gathered}
    \label{si-eq:app-AB-blended}
\end{equation}
where the superscript $\mathrm b$ denotes the blended stage. $\mathbf A_m$ and $\mathbf B_m$ are the stage matrices for restarted stages with frozen forcing. 

\subsection{Eigenvalue perturbation}
\label{si-sec:appendix-perturbation}
The stage recurrence gives $\mathbf M$ over one complete timestep. Following the decomposition introduced in section~\ref{sec:4-coupled}, we compare the coupled matrix $\mathbf M$ with an uncoupled reference $\mathbf M_{\mathrm u}$ and derive the eigenvalue corrections due to $\mathbf V=\mathbf M-\mathbf M_{\mathrm u}$. 
The uncoupled reference $\mathbf{M}_u$ is the discrete evolution matrix of each mode by itself, and it can be expressed as a block matrix for mode $q$ where
\begin{equation}
    \mathbf{M}_{\mathrm{u}, 0} = \mathbf{W} \ ,  \qquad \mathbf{M}_{\mathrm{u}, q} = R\left(-i \mu_0 \begin{bmatrix}
        0 & 1 \\ \alpha_q^2 & 0
    \end{bmatrix} \right) \ , \quad \text{for} \quad q \ge 1 \ ,
\end{equation}
where $\mathbf{W}$ is the $2\times 2$ barotropic evolution matrix and $R$ is the stability polynomial of the particular outer discretization shown in section~\ref{sec:2}.
$\mathbf M_u$ advances each mode independently. The off-diagonal modal blocks of $\mathbf M-\mathbf M_\mathrm{u}$ therefore measure the transfers generated by the coupled discretization.

We denote the eigenvalues of $\mathbf M_{\mathrm u}$ by $\lambda_j$ and its left and right eigenvectors by $\vect{\ell}_j$ and $\vect{r}_j$, respectively. The index $j$ identifies a vertical mode and a propagation direction. These satisfy
\begin{equation}
    \mathbf M_{\mathrm u}\vect r_l = \lambda_l\vect r_l \ , \qquad
    \vect{\ell}_l\mathbf M_{\mathrm u} = \lambda_l\vect{\ell}_l \ ,
\end{equation}
and we normalize them so that $\vect{\ell}_j\vect r_l=\delta_{jl}$.
We write the coupling corrections in this eigenvector basis as
\begin{equation}
    V_{jl} \defn \vect{\ell}_j\,\mathbf V\,\vect r_l \ .
\end{equation}
Then $V_{jj}$ is the correction projected back onto the same mode and propagation direction, while $V_{jl}$ describes the contribution from $l$ to $j$.
For a simple, isolated eigenvalue $\lambda_j$, we expand the corresponding eigenvalue $\lambda$ and eigenvector $\vect r$ of the coupled step as
\begin{equation}
    \big(\mathbf{M}_{\mathrm{u}} + \mathbf{V}\big)\,\vect{r} = \lambda\,\vect{r} \ , \qquad
    \lambda = \lambda_j + \delta\lambda^{(1)} + \delta\lambda^{(2)} \ , \qquad
    \vect{r} = \vect{r}_j + \delta\vect{r}^{(1)} +  \delta\vect{r}^{(2)} \ , \qquad
    \vect{\ell}_j\,\delta\vect{r}^{(1)} = 0 \ .
\end{equation}
Then, at first order
\begin{equation}
    \big(\mathbf{M}_{\mathrm{u}} - \lambda_j\mathbf{I}\big)\,\delta\vect{r}^{(1)} = \big(\delta\lambda^{(1)}\mathbf{I} - \mathbf{V}\big)\,\vect{r}_j \ .
\end{equation}
Multiplying by $\vect{\ell}_j$ gives $0 = \delta \lambda_j^{(1)} - V_{jj}$. Multiplying by $\vect{\ell}_l$ gives $(\lambda_l - \lambda_j)\vect{\ell}_l \delta \vect{r}^{(1)} = - V_{lj}$. Therefore,
\begin{equation}
    \delta\lambda^{(1)} = V_{jj} \ , \qquad
    \vect{\ell}_l\,\delta\vect{r}^{(1)} = \frac{V_{lj}}{\lambda_j - \lambda_l} \ ,
\end{equation}
The normalization $\vect\ell_j\delta\vect r^{(1)}=0$ removes the $\vect r_j$ term in the eigenvector expansion, so
\begin{equation}
 \delta\vect r^{(1)}=\sum_{l\ne j}\left(\vect\ell_l\delta\vect r^{(1)}\right)\vect r_l
 =\sum_{l\ne j}\frac{V_{lj}}{\lambda_j-\lambda_l}\vect r_l.
\end{equation}
At second order the eigenvalue equation reads
\begin{equation}
(\mathbf{M}_{\mathrm{u}} - \lambda_j\mathbf{I}\big)\,\delta\vect{r}^{(2)} = \delta \lambda^{(2)} \vect{r}_j + (\delta \lambda^{(1)}\mathbf{I} - \mathbf{V})\delta \vect{r}^{(1)} \ .
\label{si-eq:app-second-order-balance}
\end{equation}
Multiplication by $\vect{\ell}_j$ removes the left-hand side and the term $\delta \lambda^{(1)}\vect{\ell}_j \delta \vect{r}^{(1)}$. Substituting the expression for $\delta\vect{r}^{(1)}$ finally gives
\begin{equation}
    \delta\lambda^{(2)} = \vect{\ell}_j\,\mathbf{V}\,\delta\vect{r}^{(1)} = \sum_{l\neq j}\frac{V_{jl}V_{lj}}{\lambda_j - \lambda_l} \ ,
\end{equation}
which is~\eqref{eq:perturbation}. To obtain~\eqref{eq:linearized-eig-equation}, expand
\begin{equation}
 |\lambda_j+\delta\lambda_j|-|\lambda_j|
 =|\lambda_j|\left(\left|1+\frac{\delta\lambda_j}{\lambda_j}\right|-1\right)
 =\frac{\operatorname{Re}(\overline{\lambda_j}\delta\lambda_j)}{|\lambda_j|}
   +\mathcal O\left(\frac{|\delta\lambda_j|^2}{|\lambda_j|}\right).
\end{equation}
Following the comparison in Figure~\ref{fig:c4-balance}, we retain only the direct correction and the two barotropic indirect terms below.

\subsection{Direct and indirect coupling corrections}
\label{si-sec:app-transfers}
We evaluate the direct correction $V_{jj}$ and the indirect correction $\sum_b V_{jb} V_{bj} / (\lambda_j - \lambda_b)$ for a selected baroclinic mode and propagation direction $j$, where the sum includes both barotropic propagation directions.

We compute the coupled stage updates, subtract the uncoupled updates, and project the resulting corrections onto the left eigenvectors. Equation~\eqref{eq:linearized-eig-equation} then gives their contribution to growth. We first take $\varepsilon \rightarrow 0$ at fixed $\mu_0$; the small-timestep limit is considered separately below.

The matrix $[0 \ \ 1; \ \ \alpha_q^2 \ \ 0]$ has right eigenvectors $(1, \varsigma_j\alpha_q)^T$ and left eigenvectors $(1, \varsigma_j/\alpha_q)$ with eigenvalues $\varsigma_j\alpha_q$ and $\varsigma_j = \pm 1$. Applying the Runge--Kutta polynomial $R$ preserves these eigenvectors, such that, normalizing their product to one, gives the eigenvectors of $\mathbf{M}_{\mathrm{u}, q}$, 
\begin{equation}
\vect{r}_j = \frac{1}{\sqrt{2}}    \begin{bmatrix}
    1 \\ \varsigma_j \alpha_q 
\end{bmatrix} \ , \qquad \vect{\ell}_j = \frac{1}{\sqrt{2}}    \begin{bmatrix}
    1 & \dfrac{\varsigma_j}{\alpha_q}
\end{bmatrix} \ , \qquad \varsigma_j = \pm 1 \ .
\end{equation}
The displayed entries belong to mode $q$; all other entries are zero.
To compute $V_{bj}$ and $V_{jj}$, we set the initial condition to $\vect{y}^n = \vect{r}_j$ and apply one complete timestep.
To compute $V_{jb}$ we instead set $\vect{y}^n = \vect{r}_b$ using the barotropic eigenvectors
\begin{equation}
\vect{r}_b = \begin{bmatrix}
    r_{b, 1} \\ r_{b, 2}
\end{bmatrix} \ , \qquad \vect{\ell}_b = \begin{bmatrix}
    \ell_{b, 1} & \ell_{b, 2}
\end{bmatrix} \ , \qquad \mathbf{W}\vect{r}_b = \lambda_b\vect{r}_b \ , \qquad \vect{\ell}_b\mathbf{W}= \lambda_b \vect{\ell}_b \ , \qquad  \vect{\ell}_b \vect{r}_{b^\prime} = \delta_{bb^\prime} \ .
\end{equation}
Here the barotropic matrix $\mathbf{W}$ is evaluated at the barotropic Courant number $\mu_0$ and each barotropic eigenvector has non-zero components only in the barotropic block.
For either initial state $\vect{r}_k$, with $k$ equal to $j$ or $b$, $\delta u_r^M[\vect{r}_k]$ and $\delta h_r^M[\vect{r}_k]$ denote the coupled final update minus the uncoupled one from that same initial state.
\begin{equation}
    (\mathbf{V}\vect{r}_k)_r = [(\mathbf{M} - \mathbf{M}_u)\vect{r}_k]_r = \begin{bmatrix}
        \delta u_r^M[\vect{r}_k] \\ 
        \delta h_r^M[\vect{r}_k] 
    \end{bmatrix} \ , \quad \text{with} \quad r = 0, q \ , 
\end{equation}
where the subscript $r$ identifies the receiving modal block ($0$ in this case is the barotropic block) and $[\vect{r}_k]$ the initial state. Therefore
\begin{equation}
    V_{jb} = \vect{\ell}_j (\mathbf{V}\vect{r}_b)_q \ , \qquad V_{bj} = \vect{\ell}_b (\mathbf{V}\vect{r}_j)_0 \ , \qquad V_{jj} = \vect{\ell}_j (\mathbf{V}\vect{r}_j)_q \ .
\end{equation}

For the initial state $\vect{r}_j$, $\delta u_0^M[\vect{r}_j] = u_0^M$ and $\delta h_0^M[\vect{r}_j] = h_0^M$, because the uncoupled update has no mode-0 component (an initially baroclinic perturbation stays baroclinic). From the definitions of $\bar u$ and $\eta$, we can then express these update differences as
\begin{equation}
    \delta u_0^M[\vect{r}_j] = u_0^M = \frac{\bar u^M - \sum_{r>0} \chi_r u_r^M}{\chi_0} \ ,  \quad \delta h_0^M[\vect{r}_j] = h_0^M = \frac{\eta^M - \sum_{r>0} \psi_r h_r^M}{\psi_0}
    \label{si-eq:app-barotropic-differences}
\end{equation}
For $\vect y^n=\vect r_j$, the eigenvector and the modal identities in section~\ref{si-sec:appendix-modes} give
\begin{equation}
 \bar u^n=\chi_qu_q^n=\mathcal O(\varepsilon),\qquad
 \eta^n+\vect\Pi^T\vect h^n=(\psi_q+\Pi_q)h_q^n
 =\chi_qh_q^n=\mathcal O(\varepsilon^{3/2}).
\end{equation}

\subsubsection{The restarted stage}

For WRK3 and MRK4, the final stage is $M=3$ and $M=4$, respectively, and $\mu_{M-1} = \mu_0 / 2$ in both cases. We keep the barotropic outputs from~\eqref{si-eq:app-barotropic-complete} explicit, so the following transfers apply to either frozen or reconstructed forcing. 

For $\vect{y}^n = \vect{r}_b$, the mode-q entries of~\eqref{si-eq:app-stage} are
\begin{equation}
\begin{aligned}
u_q^m & = -i \mu_m h_q^{m-1} + \chi_q (\bar u^m - \bar u^n + i\mu_m \vect{\chi}^T \vect{h}^{m-1}) \ , \\
h_q^m & = -i \mu_m \alpha_q^{2}\left[u^{m-1}_q + \chi_q\left(\bar u^{\dagger, m} - \vect{\chi}^T \vect{u}^{m-1}\right)\right] \ . \\
\end{aligned}
\label{si-eq:forward-mode-updates}
\end{equation}
Starting from $u_q^n = h_q^n = 0$, the velocity update multiplies by $\chi_q$ and the pressure has the additional factor $\alpha_q^2$. Therefore, these equations give $u_q^m = \mO(\e)$ and $h_q^m = \mO(\e^2)$ at each stage.
As a consequence
\begin{equation}
\vect{\Pi}^T\vect{h}^{m} = \Pi_0 h_0^m + \sum_{r>0} \Pi_r h^m_r = \mO(\e) + \mO(\e^2) = \mO(\e)
\label{si-eq:forward-forcing-order}
\end{equation}
Frozen forcing and progressive forcing that reconstructs from these stage values therefore contribute at most $\mO(\e)$ to the barotropic outputs. In~\eqref{si-eq:app-barotropic-complete}, substituting $\bar u^n = \chi_0r_{b, 1}$ and $\eta^n = \psi_0r_{b, 2}$ gives
\begin{equation}
\begin{aligned}
\bar u^M & = W_{11}(\mu_0)\chi_0r_{b,1} + W_{12}(\mu_0)\psi_0 r_{b, 2} + \mO(\e) \\
	      & = W_{11}(\mu_0)r_{b,1} + W_{12}(\mu_0)r_{b, 2} + \mO(\e) \\
	      & = \lambda_b r_{b, 1} + \mO(\e) \ ,  \\
\eta^{M-1} & = W_{21}(\mu_0/2)\chi_0r_{b,1} + W_{22}(\mu_0/2)\psi_0 r_{b, 2} + \mO(\e) \\
	      & = W_{21}(\mu_0/2)r_{b,1} + W_{22}(\mu_0/2)r_{b, 2} + \mO(\e) \ ,
\end{aligned}
\label{si-eq:forward-barotropic-outputs}
\end{equation}
where we have omitted the higher order terms in $\chi_0$ and $\psi_0$. The final equality for $\bar u^M$ is the first row of $\mathbf{W}(\mu_0)\vect{r}_b = \lambda_b \vect{r}_b$ while the formula for $\eta^{M-1}$ uses the second row of the half step matrix.
Given that the forcing is $\mO(\e)$ in these initial conditions~\eqref{si-eq:forward-forcing-order} we can substitute $\vect{\chi}^T\vect{h}^{M-1} \approx \eta^{M-1}$. 
The differences with the uncoupled evolution, showing the leading non-zero terms in $\e$ are
\begin{equation}
 \delta u^M_q[\vect{r}_b] = \chi_q (\bar u^M - \bar u^n + i\mu_0 \eta^{M-1}) + \mO(\e^2) \ ,  \qquad  \delta h^M_q[\vect{r}_b] = \mO(\e^2) \ .
\end{equation}
The final projection on $\vect{\ell}_j$ to compute $V_{jb}$ is
\begin{equation}
V_{jb}  = \frac{1}{\sqrt{2}} \left(  \delta u^M_q[\vect{r}_b]  + \frac{\varsigma_j}{\alpha_q} \delta h^M_q[\vect{r}_b]  \right) = \frac{\delta u^M_q[\vect{r}_b]}{\sqrt{2}} + \mO(\e^{3/2})
\label{si-eq:forward-projection-step}
\end{equation}
Substituting~\eqref{si-eq:forward-barotropic-outputs} into this velocity correction gives
\begin{equation}
V_{jb} = \frac{\chi_q}{\sqrt{2}} \left[(\lambda_b - 1)r_{b, 1} + i\mu_0 \left\{ W_{21}(\mu_0/2)r_{b, 1} + W_{22}(\mu_0/2)r_{b, 2}\right\}\right] + \mO(\e^{3/2})
\label{si-eq:general-forward-transfer}
\end{equation}
Given these barotropic only initial conditions, the forcing only contributes at high order for this forward transfer. For this reason, the above expression is general and valid for both frozen and progressive forcing.

Starting from $\vect{y}^n = \vect{r}_j$ recovers the remaining transfer coefficients $V_{bj}$ and $V_{jj}$. The mode-$q$ pressure update in~\eqref{si-eq:app-stage} is
\begin{equation}
h_q^m = h_q^n - i \mu_m \alpha^2_q u^{m-1}_q - i \mu_m \alpha_q^2\chi_q(\bar u^{\dagger , m} - \vect{\chi}^T\vect{u}^{m-1}) \ .
\label{si-eq:reverse-pressure-step}
\end{equation}
For these purely baroclinic initial conditions, $\bar u^{\dagger,m}$ and $\vect{\chi}^T\vect u^{m-1}$ are both $\mO(\e)$. Since $\alpha_q^2\chi_q=\mO(\e^2)$, the last term in~\eqref{si-eq:reverse-pressure-step} is $\mO(\e^3)$. Substituting $u_q^{m-1}=u_q^n+\mO(\e^{1/2})$ into the remaining velocity term introduces an $\mO(\e^{3/2})$ remainder, giving
\begin{equation}
u_q^m = u_q^n + \mO(\e^{1/2}) \ ,  \qquad h_q^m=h_q^n-i\mu_m\alpha_q^2u_q^n+\mO(\e^{3/2}) \ .
\label {si-eq:general-stage-pressures}
\end{equation}
For $r > 0, r \ne q$, the same equations start from $u_r^n = h_r^n = 0$ and give $u^m_r = \mO(\e^2)$ and $h_r^m = \mO(\e^3)$. Therefore, their contributions to each sum in~\eqref{si-eq:app-barotropic-differences} are $\mO(\e^3)$. Substituting the amplitudes of mode $q$ explicitly gives
\begin{equation}
\begin{aligned}
\chi_0 \delta u^M_0[\vect{r}_j] & = \bar u^M - \chi_q[u_q^n + \mO(\e^{1/2})] + \mO(\e^3) = \bar u^M - \bar u^n + \mO(\e^{3/2})  \ , \\
\psi_0 \delta h^M_0[\vect{r}_j] & = \eta^n - i\mu_0 \bar u^{\dagger, M} - \psi_q[h_q^n - i\mu_0 \alpha_q^2 u_q^n] + \mO(\e^{3/2})   \\
& = -i\mu_0 (\bar u^{\dagger, M} - \bar u^n) + \mO(\e^{3/2}) \ .
\end{aligned}
\label{si-eq:reverse-modal-substitution}
\end{equation}
The last line uses $\eta^n = \psi_qh^n_q$ and $\psi_q\alpha_q^2u_q^n = \chi_q u^n_q = \bar u^n$. We can then approximate $\chi_0, \psi_0 \approx 1$ at this leading order, and projecting onto $\vect{\ell}_b$ gives
\begin{equation}
V_{bj} = \vect{\ell_b} \begin{bmatrix}
\bar u^M - \bar u^n \\ 
 -i\mu_0 (\bar u^{\dagger, M} - \bar u^n)
\end{bmatrix} + \mO(\e^{3/2}) \ .
\label{si-eq:general-reverse-transfer}
\end{equation}

For the same initial state $\vect{r}_j$, the uncoupled restarted update is 
\begin{equation}
 u_{q, \mathrm{u}}^m = u_q^n - i\mu_m h_{q, \mathrm{u}}^{m-1} \ , \qquad  h_{q, \mathrm{u}}^m = h_q^n - i\mu_m\alpha_q^2 u_{q, \mathrm{u}}^{m-1} \ , 
 \label {si-eq:uncoupled-restarted-stages}
\end{equation}
where the subscript $\mathrm{u}$ denotes the uncoupled amplitudes, so $\delta u_q^m  = u_q^m - u_{q, \mathrm{u}}^m$ and $\delta h_q^m  = h_q^m - h_{q, \mathrm{u}}^m$. Subtracting~\eqref{si-eq:uncoupled-restarted-stages} from~\eqref{si-eq:app-stage} gives
\begin{equation}
\begin{aligned}
\delta u^m_q & = - i \mu_m \delta h_q^{m-1} + \chi_q (\bar u^m - \bar u^n + i\mu_m \vect{\chi}^T\vect{h}^{m-1})  \ ,  \\ 
\delta h^m_q & = - i \mu_m \alpha_q^2 \left[\delta u_q^{m-1} + \chi_q (\bar u^{\dagger, m} - \vect{\chi}^T\vect{u}^{m-1})\right]  \  .
\end{aligned}
\label {si-eq:general-stage-corrections}
\end{equation}
In the first equation, the reconciliation is $\chi_q \mO(\e) = \mO(\e^2)$. In the second, the transport counterpart has the additional factor $\alpha_q^2$ and is $\mO(\e^3)$. Since the coupled and uncoupled evolutions start from the same state, $\delta u_q^0=\delta h_q^0=0$. For successive stages $\delta u_q^m = \mO(\e^2)$ and $\delta h_q^m = \mO(\e^3)$. To evaluate $\vect{\chi}^T\vect{h}^m$, we apply the modal-sum identity of~\eqref{si-eq:reverse-modal-substitution} at stage $m$, 
\begin{equation}
\psi_0h_0^m
=-i\mu_m\left(\bar u^{\dagger,m}-\bar u^n\right)
+\mathcal O(\varepsilon^{3/2}).
\label {si-eq:weighted-pressure-elimination}
\end{equation}
Substituting into the modal sum gives
\begin{equation}
\begin{aligned}
\boldsymbol\chi^T\mathbf h^m
&=\frac{\chi_0}{\psi_0}\,\psi_0h_0^m
+\chi_qh_q^m+\mathcal O(\varepsilon^3)\\
&=-i\mu_m\left(\bar u^{\dagger,m}-\bar u^n\right)
+\mathcal O(\varepsilon^{3/2}),
\end{aligned}
\end{equation}
where we used $\chi_0/\psi_0=1+\mathcal O(\varepsilon)$, $\bar u^{\dagger,m}-\bar u^n=\mathcal O(\varepsilon)$, and $\chi_qh_q^m=\mathcal O(\varepsilon^{3/2})$. At the final stage, the equations in~\eqref{si-eq:general-stage-corrections} become
\begin{equation}
\delta u_q^M = - i \mu_0 \delta h_q^{M-1} + \chi_q (\bar u^M - \bar u^n) + i\mu_0 \chi_q \left[ -i\mu_{M-1} (\bar u^{\dagger, M-1} - \bar u^n) + \mO(\e^{3/2})\right] \ , \qquad \delta h^M_q = \mO(\e^3) \ .
\label {si-eq:direct-final-substitution}
\end{equation}
Since $\mu_{M-1} = \mu_0 / 2$, the equation for the velocity difference gives
\begin{equation}
\delta u_q^M = \chi_q \left[ \bar u^M - \bar u^n + \frac{\mu_0^2}{2} (\bar u^{\dagger, M-1} - \bar u^n) \right]  + \mO(\e^{5/2}) \ .
\end{equation}
The pressure correction contributes $\delta h^M_q / \alpha_q = \mO(\varepsilon^{5/2})$ to the projection. Therefore, 
\begin{equation}
V_{jj} = \frac{1}{\sqrt{2}}\left(\delta u_q^M + \frac{\varsigma_j}{\alpha_q}\delta h^M_q\right) = \frac{\chi_q}{\sqrt{2}}\left[ \bar u^M - \bar u^n + \frac{\mu_0^2}{2} (\bar u^{\dagger, M-1} - \bar u^n) \right]  + \mO(\e^{5/2}) \ .
\label {si-eq:general-direct-projection}
\end{equation}
Finally, \cref{si-eq:general-forward-transfer,si-eq:general-reverse-transfer,si-eq:general-direct-projection} give the leading transfers for frozen or reconstructed stage forcing through its barotropic
outputs. Their substitution for progressive forcing and the resulting leading cancellation are discussed in section~\ref{sec:4}.

We now compute the eigenvalue shift $\delta \lambda$.  We start with computing the indirect correction $\delta\lambda^{(2)}$. At fixed $\mu_0$, away from eigenvalue coincidences, $\lambda_j = 1 + \mO(\e^{1/2})$, so we approximate $\lambda_j\approx 1$ in the denominator. Substituting $V_{bj}$ in the expression for $\delta\lambda^{(2)}$ we get
\begin{equation}
\delta\lambda_j^{(2)}
=
\left[\sum_b\frac{V_{jb}\boldsymbol\ell_b}{1-\lambda_b}\right]
\begin{bmatrix}
\bar u^M-\bar u^n\\
-i\mu_0(\bar u^{\dagger,M}-\bar u^n)
\end{bmatrix}
+\mathcal O(\varepsilon^{5/2}).
\end{equation}
Using $\mathbf{W} = \sum_b \lambda_b\vect{r}_b\vect{\ell}_b$ and $\vect{\ell}_b\vect{r}_{b^\prime} = \delta_{bb^\prime}$, we have that
\begin{equation}
\sum_b \vect{r}_b\vect{\ell}_b = \mathbf{I} \ , \qquad \sum_b \frac{\vect{r}_b\vect{\ell}_{b}}{1 - \lambda_b} = \left[\mathbf{I} - \mathbf{W}(\mu_0)\right]^{-1} \ .
\end{equation}
The two terms in $V_{jb}$ then give
\begin{equation}
\sum_b
\frac{(\lambda_b-1)r_{b,1}\boldsymbol\ell_b}{1-\lambda_b}
=-\sum_b r_{b,1}\boldsymbol\ell_b
=-\begin{bmatrix}1&0\end{bmatrix},
\end{equation}
and
\begin{equation}
\sum_b
\frac{\left[W_{21}(\mu_0/2)r_{b,1}
+W_{22}(\mu_0/2)r_{b,2}\right]\boldsymbol\ell_b}
{1-\lambda_b}
=
\begin{bmatrix}W_{21}(\mu_0/2)&W_{22}(\mu_0/2)\end{bmatrix}
[\mathbf I-\mathbf W(\mu_0)]^{-1} \ .
\label {si-eq:indirect-eigenvector-contractions}
\end{equation}
Substituting the two sums into the indirect correction gives
\begin{equation}
\begin{aligned}
\delta\lambda_j^{(2)}
=\frac{\chi_q}{\sqrt2}
\Bigg[
-\begin{bmatrix}1&0\end{bmatrix}
+i\mu_0
\begin{bmatrix}W_{21}(\mu_0/2)&W_{22}(\mu_0/2)\end{bmatrix}
[\mathbf I-\mathbf W(\mu_0)]^{-1}
\Bigg]
\begin{bmatrix}
\bar u^M-\bar u^n\\
-i\mu_0(\bar u^{\dagger,M}-\bar u^n)
\end{bmatrix}
+\mathcal O(\varepsilon^{5/2}) \ .
\end{aligned}
\end{equation}
Which we can rewrite as
\begin{equation}
\begin{aligned}
\delta\lambda_j^{(2)}
={}&-\frac{\chi_q}{\sqrt2}(\bar u^M-\bar u^n)\\
&+\frac{i\chi_q\mu_0}{\sqrt2}
\begin{bmatrix}W_{21}(\mu_0/2)&W_{22}(\mu_0/2)\end{bmatrix}
[\mathbf I-\mathbf W(\mu_0)]^{-1}
\begin{bmatrix}
\bar u^M-\bar u^n\\
-i\mu_0(\bar u^{\dagger,M}-\bar u^n)
\end{bmatrix}
+\mathcal O(\varepsilon^{5/2}) \ .
\end{aligned}
\end{equation}
Adding $V_{jj}$ cancels the first term and adds a term proportional to the barotropic transport
\begin{equation}
\begin{aligned}
\delta\lambda_j
&=\frac{\chi_q\mu_0^2}{2\sqrt{2}}(\bar u^{\dagger, M-1}-\bar u^n)\\
&+\frac{i\chi_q\mu_0}{\sqrt2}
\begin{bmatrix}W_{21}(\mu_0/2)&W_{22}(\mu_0/2)\end{bmatrix}
[\mathbf I-\mathbf W(\mu_0)]^{-1}
\begin{bmatrix}
\bar u^M-\bar u^n\\
-i\mu_0(\bar u^{\dagger,M}-\bar u^n)
\end{bmatrix}
+\mathcal O(\varepsilon^{5/2}) \ .
\end{aligned}
\label {si-eq:general-leading-coupling}
\end{equation}

\subsubsection{Frozen forcing}

The leading expression for $V_{jb}$ does not depend on the forcing, so it is already reduced and suitable for frozen forcing. We now simplify $V_{bj}$ and $V_{jj}$.
For frozen forcing, $F^m_0 = \vect{\Pi}^T\vect{h}^{m-1}$ and $F^m_p = 0$ for $p>0$. The constant-forcing barotropic outputs are $\bar u_0^m = W_{12}(\mu_m)$ and $\bar u_0^{\dagger, m} = \widetilde{W}_{12}(\mu_m)$ as derived in the appendix above. For an $\vect{r}_j$ initial state, equation~\eqref{si-eq:general-stage-pressures} and $\Pi_q \alpha_q^2 = \chi_q(\alpha_q^2 -1)$ give
\begin{equation}
F_0^m = \Pi_qh_q^n - i\mu_{m-1} \Pi_q\alpha_q^2 u_q^n + \mO(\e^{3/2}) = F_0^1 + i\mu_{m-1} \bar u^n + \mO(\e^{3/2})
\label{si-eq:frozen-stage-values}
\end{equation}
Substituting the above into~\eqref{si-eq:app-barotropic} and using the fact that, for this initial baroclinic state, $\eta^n + F_0^1 = \mO(\e^{3/2})$, we can derive
\begin{equation}
\begin{aligned}
\bar u^M - \bar u^n & = \left[ W_{11}(\mu_0) - 1 + \frac{i\mu_0}{2}W_{12}(\mu_0) \right]\bar u^n + \mO(\e^{3/2}) \ , \\
- i\mu_0(\bar u^{\dagger, M} - \bar u^n) & = \left[ W_{21}(\mu_0) + \frac{i\mu_0}{2}(1+ W_{22}(\mu_0))\right] \bar u^n + \mO(\e^{3/2}) \ , \\
\bar u^{\dagger, M-1} - \bar u^n & = \left[ \widetilde{W}_{11}(\mu_0/2) - 1 + \frac{i\mu_0}{3}\widetilde{W}_{12}(\mu_0/2) \right] \bar u^n  + \mO(\e^{3/2}) \ .
\end{aligned}
\label{si-eq:frozen-output-differences}
\end{equation}
Inserting the first two rows into~\eqref{si-eq:general-reverse-transfer} and using
\begin{equation}
\begin{aligned}
\ell_{b,1}(W_{11} − 1) + \ell_{b,2}W_{21} & = (\lambda_b − 1)\ell_{b,1} \ ,  \\
\ell_{b,1}W_{12} + \ell_{b,2}(1 +W_{22})&  = (1 + \lambda_b)\ell_{b,2} \ ,
\end{aligned}
\label {si-eq:frozen-left-projection}
\end{equation}
with $\bar u^n = \chi_q / \sqrt2$ gives a closed expression for the frozen forcing reverse transfer coefficient
\begin{equation}
V_{bj} = \frac{\chi_q}{\sqrt2}\left[ (\lambda_b - 1)\ell_{b, 1} + \frac{i\mu_0}{2} (1+\lambda_b)\ell_{b, 2}\right] + \mO(\e^{3/2}) \ .
\label {si-eq:frozen-reverse-transfer}
\end{equation}
For $V_{jj}$, the term involving $\bar u^{\dagger,M-1}-\bar u^n$ in~\eqref{si-eq:general-direct-projection} becomes
\begin{equation}
\frac{\mu_0^2}{2}\left[ \widetilde{W}_{11}(\mu_0/2) - 1 + \frac{i\mu_0}{3}\widetilde{W}_{12}(\mu_0/2) \right]  = i\mu_0 W_{21}(\mu_0/2) - \frac{\mu_0^2}{2} - \frac{\mu_0^2}{3}[W_{22}(\mu_0/2) - 1] \ , 
\label {si-eq:frozen-half-step-transport}
\end{equation}
and adding the final velocity from~\eqref{si-eq:frozen-output-differences} gives
\begin{equation}
V_{jj} = \frac{\chi_q^2}{2}\left[W_{11}(\mu_0) - 1 + \frac{i\mu_0}{2} W_{12}(\mu_0) + i\mu_0 W_{21}(\mu_0/2) - \frac{\mu_0^2}{3}(W_{22}(\mu_0/2)-1) - \frac{\mu_0^2}{2} \right] + \mO(\e^{5/2}) \ .
\label {si-eq:direct-finite-step}
\end{equation}
For the combined correction, the first two rows of~\eqref{si-eq:frozen-output-differences} give
\begin{equation}
\begin{bmatrix}
\bar u^M - \bar u^n \\ 
- i\mu_0(\bar u^{\dagger, M} - \bar u^n) 
\end{bmatrix} = (\mathbf W(\mu_0) - \mathbf I) 
\begin{bmatrix} 
1 \\
i\mu_0 / 2
\end{bmatrix} \bar u^n + 
\begin{bmatrix}
0 \\ i\mu_0
\end{bmatrix}\bar u^n + \mO(\e^{3/2})
\end{equation}
Multiplication by $(\mathbf I - \mathbf W(\mu_0))^{-1}$ changes the first term to $- (1, i\mu_0/2)^T\bar u^n$. Substituting the vector in~\eqref{si-eq:general-leading-coupling} with the explicit expressions in~\eqref{si-eq:frozen-output-differences} and using~\eqref{si-eq:frozen-half-step-transport}, while factoring out $\chi_q \bar u^n / \sqrt2$, gives
\begin{equation}
\begin{aligned}
    \delta\lambda_j
    = \frac{\chi_q\bar u^n}{\sqrt{2}}
    \Bigg[
        &i\mu_0 W_{21}(\mu_0/2)
        - \frac{\mu_0^2}{2}
        - \frac{\mu_0^2}{3}\big[W_{22}(\mu_0/2)-1\big] \\
        -& i\mu_0 W_{21}(\mu_0/2)
        + \frac{\mu_0^2}{2}W_{22}(\mu_0/2) \\
        - &\mu_0^2
        \left\{
            \mathbf W(\mu_0/2)
            \big[\mathbf I-\mathbf W(\mu_0)\big]^{-1}
        \right\}_{22}
    \Bigg]
    + \mathcal O(\varepsilon^{5/2}).
\end{aligned}
\label{si-eq:frozen-combined-substitution}
\end{equation}
The terms in $W_{21}$ cancel. Here $\bar u^n$ corresponds to the input $\vect y^n=\vect r_j$, so \(\bar u^n=\chi_q/\sqrt2\), giving 
\begin{equation}
    \delta\lambda_j
    = \frac{\chi_q^2\mu_0^2}{12}
    \left[
        W_{22}(\mu_0/2)-1
        -6\left\{
            \mathbf W(\mu_0/2)
            \big[\mathbf I-\mathbf W(\mu_0)\big]^{-1}
        \right\}_{22}
    \right]
    + \mathcal O(\varepsilon^{5/2}).
    \label{si-eq:generic-total-coupling}
\end{equation}
Together with~\eqref{si-eq:general-forward-transfer}, these are the frozen-forcing expressions used in section~\ref{sec:4}.

\subsubsection{The blended stage}

Subtracting the uncoupled Shu–Osher update from the coupled predictor and reconciliation \cref{eq:ssp-velocity,eq:ssp-reconcile} gives
\begin{equation}
\delta u^m_q = b_m \delta u_q^{m-1} - ib_m \mu_0 \delta h^{m-1}_q + \chi_q \left[ \bar u^m - (1 - b_m)\bar u^n - b_m \bar u^{m-1} + ib_m \mu_0 \vect{\chi}^T\vect{h}^{m-1}\right] \ .
\label {si-eq:ssp-velocity-difference}
\end{equation}
We move $\chi_q(\bar u^m - \bar u^n)$ to the left hand side. The remaining velocity difference is multiplied by $b_m$ at each stage, therefore, substitution of stages one and two into stage three yields
\begin{equation}
\delta u^3_q - \chi_q (\bar u^3 - \bar u^n) = i\mu_0 \left[ b_3b_2b_1 (\chi_q\vect{\chi}^T\vect{h}^n - \delta h^0_q)+ b_3b_2(\chi_q\vect{\chi}^T\vect{h}^1 - \delta h^1_q) +b_3(\chi_q\vect{\chi}^T\vect{h}^2 - \delta h^2_q)\right] 
\label {si-eq:ssp-unrolled-velocity}
\end{equation}
where $(b_3b_2b_1, b_3b_2, b_3) = (1/6, 1/6, 2/3) = \beta$.

Once again, for $V_{jb}$, we take $\vect{y}^n = \vect{r}_b$. As for the restarted stage, the forcing is $\mO(\e)$, so we can consider the unforced barotropic evolution. Using $b_1 = 1$ and $b_2 = 1/4$ in~\eqref{eq:ssp-etablend} gives
\begin{equation}
\begin{aligned}
 \eta^1&=W_{21}(\mu_0)\bar u^n+W_{22}(\mu_0)\eta^n+\mathcal O(\varepsilon),\\
 \eta^2&=\frac34\eta^n+\frac14\left[\eta^1+(\eta^1-\eta^n)\right]+\mathcal O(\varepsilon)
        =\frac{\eta^n+\eta^1}{2}+\mathcal O(\varepsilon).
\end{aligned}
\end{equation}
The weighted free surface value entering the velocity update is
\begin{equation}
\begin{aligned}
 \sum_{m=1}^3\beta_m\eta^{m-1}
 &=\frac16\eta^n+\frac16\eta^1+\frac23\eta^2
 =\frac{\eta^n+\eta^1}{2}+\mathcal O(\varepsilon)\\
 &=\frac12W_{21}(\mu_0)\bar u^n+\frac12[1+W_{22}(\mu_0)]\eta^n+\mathcal O(\varepsilon).
\end{aligned}
\label{si-eq:ssp-weighted-surface}
\end{equation}
For this input $\vect{\chi}^T\vect{h}^{m-1} = \eta^{m-1} + \mO(\e)$ by~\eqref{si-eq:forward-forcing-order}, while $\delta h^{m-1}_q = \mO(\varepsilon^2)$. Equation~\eqref{si-eq:ssp-unrolled-velocity} consequently gives
\begin{equation}
\delta u^3_q[\vect{r}_b] = \chi_q \left[ \bar u^3 - \bar u^n + i\mu_0 \sum_{m=1}^3 \beta_m \eta^{m-1}\right] + \mO(\e^2) \ .
\label {si-eq:ssp-forward-velocity}
\end{equation}
The unforced velocity increments are identical in each solve. Their blend returns
\begin{equation}
\begin{aligned}
    \bar u^3-\bar u^n
    &= \left(\sum_{m=1}^{3}\beta_m\right)
       \left[(W_{11}-1)\bar u^n+W_{12}\eta^n\right]
       + \mathcal O(\varepsilon) \\
    &= (\lambda_b-1)r_{b,1}+\mathcal O(\varepsilon).
\end{aligned}
\label{si-eq:ssp-final-barotropic-velocity}
\end{equation}
where $\sum\beta = 1$ and the substitution of the initial state is the same as in~\eqref{si-eq:forward-barotropic-outputs}. The second row of $\mathbf{W}(\mu_0) \vect{r}_b = \lambda_b \vect{r}_b$ reduces~\eqref{si-eq:ssp-weighted-surface} to 
\begin{equation}
 W_{21}(\mu_0)r_{b,1}+W_{22}(\mu_0)r_{b,2}=\lambda_b r_{b,2},
 \qquad
 \sum_{m=1}^3\beta_m\eta^{m-1}=\frac{1+\lambda_b}{2}r_{b,2}+\mathcal O(\varepsilon).
\end{equation}
The pressure update contributes $\delta h_q^3[\vect r_b]=\mathcal O(\alpha_q^2\chi_q)=\mathcal O(\varepsilon^2)$. Its projection is therefore order $\varepsilon^{3/2}$, and projection onto $\vect\ell_j$ gives
\begin{equation}
\begin{aligned}
 V_{jb}&=\frac1{\sqrt2}\left(\delta u_q^3[\vect r_b]
                  +\frac{\varsigma_j}{\alpha_q}\delta h_q^3[\vect r_b]\right)\\
 &=\frac{\chi_q}{\sqrt2}\left[(\lambda_b-1)r_{b,1}
         +\frac{i\mu_0}{2}(1+\lambda_b)r_{b,2}\right]
         +\mathcal O(\varepsilon^{3/2}).
\end{aligned}
\label{si-eq:ssp-forward-transfer}
\end{equation}

We now derive the coefficient with initial state $\vect{r}_j$. For this initial condition, the blended pressure state gives
\begin{equation}
 h_q^1 = h_q^n - i\mu_0 \alpha_q^2 u_q^n + \mO(\e^{3/2}) \ , \qquad h^2_q = h_q^n  -  \frac{i\mu_0}{2} \alpha_q^2 u_q^n + \mO(\e^{3/2}) \ , 
\end{equation}
where the expression for $h^2_q$ follows from $h_q^2 = 3h_q^n/4 + (h_q^1 - i\mu_0\alpha^2_q u_q^1)/4$. Multiplying by $\Pi_q$ gives
\begin{equation}
\begin{aligned}
    \vect{\Pi}^T\vect h^1
    &= \vect{\Pi}^T\vect h^n+i\mu_0\bar u^n
       +\mathcal O(\varepsilon^{3/2}), \\
    \vect{\Pi}^T\vect h^2
    &= \vect{\Pi}^T\vect h^n+\frac{i\mu_0}{2}\bar u^n
       +\mathcal O(\varepsilon^{3/2}).
\end{aligned}
\label{si-eq:ssp-stage-forcing-values}
\end{equation}
Unrolling the increment blend weights the three barotropic solves by $\beta=(1/6,1/6,2/3)$. The solves use frozen forcing, so the effective constant forcing is
\begin{equation}
\begin{aligned}
    &\frac16\vect{\Pi}^T\vect h^n
    +\frac16\vect{\Pi}^T\vect h^1
    +\frac23\vect{\Pi}^T\vect h^2
 \\
    &\qquad=
    \vect{\Pi}^T\vect h^n + i\mu_0 \left(\frac{1}{6} + \frac{2}{3}\frac{1}{2}\right) \bar u^n + \mO(\e^{3/2})
    =\vect{\Pi}^T\vect h^n+\frac{i\mu_0}{2}\bar u^n
    +\mathcal O(\varepsilon^{3/2}).
\end{aligned}
\end{equation}
All three solves use the same $\mathbf W(\mu_0)$ and initial state. Summing their increments with these weights, and using
$\eta^n+\vect{\Pi}^T\vect h^n=\mathcal O(\varepsilon^{3/2})$, gives
\begin{equation}
\begin{aligned}
    \bar u^3-\bar u^n
    &=\left[W_{11}-1+\frac{i\mu_0}{2}W_{12}\right]\bar u^n
      +\mathcal O(\varepsilon^{3/2}), \\
    \eta^3-\eta^n
    &=\left[W_{21}+\frac{i\mu_0}{2}(W_{22}-1)\right]\bar u^n
      +\mathcal O(\varepsilon^{3/2}).
\end{aligned}
\label{si-eq:ssp-final-output-differences}
\end{equation}
The final pressure blend gives
\begin{equation}
\begin{aligned}
    h_q^3
    &=\frac13h_q^n+\frac23\left(h_q^2-i\mu_0\alpha_q^2u_q^2\right)
      +\mathcal O(\varepsilon^3) \\
    &=h_q^n-i\mu_0\alpha_q^2u_q^n
      +\mathcal O(\varepsilon^{3/2}).
\end{aligned}
\label{si-eq:ssp-final-mode-pressure}
\end{equation}
where $h_q^2$ is the half-step value above and $u_q^2=u_q^n+\mathcal O(\varepsilon^{1/2})$. Each solve spans the full step and satisfies the transport identity, so the increment blend gives
\begin{equation}
    \eta^3-\eta^n=\sum_m\beta_m\Delta\eta^m=-i\mu_0\bar u^{\dagger,3} \ , \qquad
    \bar u^{\dagger,3}\defn\sum_m\beta_m\bar u^{\dagger,m} \ .
    \label{si-eq:ssp-blended-transport}
\end{equation}
In~\eqref{si-eq:app-barotropic-differences}, these amplitudes give
\begin{equation}
\begin{aligned}
    \chi_0\delta u_0^3[\vect r_j]
    &=\bar u^3-\bar u^n+\mathcal O(\varepsilon^{3/2}), \\
    \psi_0\delta h_0^3[\vect r_j]
    &=\eta^3-\eta^n+i\mu_0\bar u^n
      +\mathcal O(\varepsilon^{3/2})
     =-i\mu_0\left(\bar u^{\dagger,3}-\bar u^n\right)
      +\mathcal O(\varepsilon^{3/2}) \\
    &=\left[W_{21}+\frac{i\mu_0}{2}(1+W_{22})\right]\bar u^n
      +\mathcal O(\varepsilon^{3/2}) \ ,
\end{aligned}
\label{si-eq:ssp-reverse-modal-substitution}
\end{equation}
where the last line uses the frozen outputs~\eqref{si-eq:ssp-final-output-differences}. The first two lines hold for any forcing and, replacing $\chi_0$ and $\psi_0$ by one at this order, give the reverse transfer in the same form as for the restarted stage~\eqref{si-eq:general-reverse-transfer},
\begin{equation}
    V_{bj} = \vect{\ell}_b
    \begin{bmatrix}
        \bar u^3-\bar u^n \\
        -i\mu_0(\bar u^{\dagger,3}-\bar u^n)
    \end{bmatrix}
    +\mathcal O(\varepsilon^{3/2}) \ .
    \label{si-eq:ssp-general-reverse-transfer}
\end{equation}
For frozen forcing, the last line and the left-eigenvector identities in~\eqref{si-eq:frozen-left-projection} give
\begin{equation}
    V_{bj}
    =\frac{\chi_q}{\sqrt2}
    \left[
        (\lambda_b-1)\ell_{b,1}
        +\frac{i\mu_0}{2}(1+\lambda_b)\ell_{b,2}
    \right]
    +\mathcal O(\varepsilon^{3/2}).
\end{equation}
which is equivalent to the explicit expression derived for frozen forcing in the restarted stage~\eqref{si-eq:frozen-reverse-transfer}.

For the same input, $h_0^n=0$. In~\eqref{si-eq:ssp-unrolled-velocity}, we use $\vect{\chi}^T\vect h^{m-1}=h_0^{m-1}+\mathcal O(\varepsilon^{3/2})$ and $\delta h_q^{m-1}=\mathcal O(\varepsilon^3)$. The retained terms are
\begin{equation}
    \delta u_q^3[\vect r_j]
    =\chi_q\left[
        \bar u^3-\bar u^n
        +i\mu_0\left(\frac16h_0^1+\frac23h_0^2\right)
    \right]
    +\mathcal O(\varepsilon^{5/2}).
    \label{si-eq:ssp-direct-velocity}
\end{equation}
The discarded pressure difference enters $V_{jj}$ divided by $\alpha_q$, at order $\varepsilon^{5/2}$. To calculate $h_0^1$ and $h_0^2$, we combine the definition of the barotropic pressure, the transport identity of each solve with the blend $b_2=1/4$, and the two pressure stages above,
\begin{equation}
\begin{aligned}
    \psi_0h_0^m &= \eta^m-\eta^n-\psi_q\left(h_q^m-h_q^n\right) \ , \\
    \eta^1-\eta^n &= -i\mu_0\bar u^{\dagger,1} \ , &
    \psi_q\left(h_q^1-h_q^n\right) &= -i\mu_0\bar u^n \ , \\
    \eta^2-\eta^n &= -\frac{i\mu_0}{4}\left(\bar u^{\dagger,1}+\bar u^{\dagger,2}\right) \ , &
    \psi_q\left(h_q^2-h_q^n\right) &= -\frac{i\mu_0}{2}\bar u^n \ .
\end{aligned}
\end{equation}
At leading order, this gives
\begin{equation}
    h_0^1 = -i\mu_0\left(\bar u^{\dagger,1}-\bar u^n\right) \ , \qquad
    h_0^2 = -\frac{i\mu_0}{4}\left(\bar u^{\dagger,1}+\bar u^{\dagger,2}-2\bar u^n\right) \ ,
    \label{si-eq:ssp-stage-surfaces}
\end{equation}
and~\eqref{si-eq:ssp-direct-velocity} gives, for any forcing,
\begin{equation}
    V_{jj} = \frac{\chi_q}{\sqrt2}\left[\bar u^3-\bar u^n+\frac{\mu_0^2}{6}\left(2\bar u^{\dagger,1}+\bar u^{\dagger,2}-3\bar u^n\right)\right]+\mathcal O(\varepsilon^{5/2}) \ .
    \label{si-eq:ssp-general-direct-projection}
\end{equation}
Summing over $b$ with the forward transfer~\eqref{si-eq:ssp-forward-transfer} as for the restarted stage gives
\begin{equation}
    \delta\lambda_j = \frac{\chi_q}{\sqrt2}\left[\frac{\mu_0^2}{6}\left(2\bar u^{\dagger,1}+\bar u^{\dagger,2}-3\bar u^n\right)
    +\frac{i\mu_0}{2}\left\{\left[\mathbf I+\mathbf W(\mu_0)\right]\left[\mathbf I-\mathbf W(\mu_0)\right]^{-1}
    \begin{bmatrix}
        \bar u^3-\bar u^n \\
        -i\mu_0(\bar u^{\dagger,3}-\bar u^n)
    \end{bmatrix}\right\}_2\right]+\mathcal O(\varepsilon^{5/2}) \ .
    \label{si-eq:ssp-general-leading-coupling}
\end{equation}
For frozen forcing, the first constant-forcing solve gives
$\eta^1-\eta^n=W_{21}(\mu_0)\bar u^n+\mathcal O(\varepsilon^{3/2})$. In the second solve, the extra forcing $i\mu_0\bar u^n$ from~\eqref{si-eq:ssp-stage-forcing-values} contributes
$i\mu_0(W_{22}-1)\bar u^n$ to the increment. Applying $b_2=1/4$ gives
\begin{equation}
    \eta^2-\eta^n
    =\frac{\bar u^n}{4}
    \left[2W_{21}(\mu_0)+i\mu_0\big(W_{22}(\mu_0)-1\big)\right]
    +\mathcal O(\varepsilon^{3/2}).
\end{equation}
The same steps as above yield
\begin{equation}
\begin{aligned}
    h_0^1
    &=\bar u^n\left[W_{21}(\mu_0)+i\mu_0\right]
      +\mathcal O(\varepsilon^{3/2}), \\
    h_0^2
    &=\frac{\bar u^n}{4}
      \left[2W_{21}(\mu_0)+i\mu_0\big(1+W_{22}(\mu_0)\big)\right]
      +\mathcal O(\varepsilon^{3/2}).
\end{aligned}
\end{equation}
The combination required in~\eqref{si-eq:ssp-direct-velocity} is therefore
\begin{equation}
    \frac16h_0^1+\frac23h_0^2
    =\bar u^n\left[
        \frac12W_{21}(\mu_0)
        +\frac{i\mu_0}{6}\big(2+W_{22}(\mu_0)\big)
    \right]
    +\mathcal O(\varepsilon^{3/2}).
\end{equation}
Substitution into~\eqref{si-eq:ssp-direct-velocity}, using the final velocity from the reverse-transfer calculation and $\bar u^n=\chi_q/\sqrt2$, gives
\begin{equation}
    V_{jj}
    =\frac{\chi_q^2}{2}
    \left[
        W_{11}(\mu_0)-1
        +\frac{i\mu_0}{2}\big(W_{12}(\mu_0)+W_{21}(\mu_0)\big)
        -\frac{\mu_0^2}{6}\big(2+W_{22}(\mu_0)\big)
    \right]
    +\mathcal O(\varepsilon^{5/2}).
    \label{si-eq:ssp-direct-transfer}
\end{equation}

Finally, to compute the total eigenvalue shift, we multiply the forward and reverse transfers and sum over $b$. Using completeness, $\sum_b\lambda_b\vect r_b\vect\ell_b=\mathbf W$, and
\begin{equation}
    \frac{(1+\lambda_b)^2}{1-\lambda_b}
    =\frac{4}{1-\lambda_b}-3-\lambda_b,
\end{equation}
the velocity, cross, and pressure products give, respectively,
\begin{equation}
\begin{aligned}
    \sum_b\frac{(\lambda_b-1)^2r_{b,1}\ell_{b,1}}{1-\lambda_b}
    &=1-W_{11}, \\
    \sum_b\frac{(\lambda_b-1)(1+\lambda_b)
        (r_{b,1}\ell_{b,2}+r_{b,2}\ell_{b,1})}{1-\lambda_b}
    &=-(W_{12}+W_{21}), \\
    \sum_b\frac{(1+\lambda_b)^2r_{b,2}\ell_{b,2}}{1-\lambda_b}
    &=4\left\{(\mathbf I-\mathbf W)^{-1}\right\}_{22}-3-W_{22}.
\end{aligned}
\label{si-eq:ssp-product-contractions}
\end{equation}
Multiplying the second and third rows by $i\mu_0/2$ and $-\mu_0^2/4$ gives
\begin{equation}
\begin{aligned}
    \sum_b\frac{V_{jb}V_{bj}}{\lambda_j-\lambda_b}
    =\frac{\chi_q^2}{2}\Bigg[
        &1-W_{11}-\frac{i\mu_0}{2}(W_{12}+W_{21}) \\
        &-\mu_0^2\left\{(\mathbf I-\mathbf W)^{-1}\right\}_{22}
        +\frac{\mu_0^2}{4}(3+W_{22})
    \Bigg]
    +\mathcal O(\varepsilon^{5/2}),
\end{aligned}
\end{equation}
where $\mathbf W=\mathbf W(\mu_0)$ and $\lambda_j$ has been replaced by one at leading order. Adding~\eqref{si-eq:ssp-direct-transfer} cancels $1-W_{11}$ and the cross terms. The remaining polynomial coefficient is
\begin{equation}
    -\frac{2+W_{22}}{6}+\frac{3+W_{22}}{4}
    =\frac{5+W_{22}}{12}.
\end{equation}
This gives
\begin{equation}
 \delta\lambda_j=\frac{\chi_q^2\mu_0^2}{24}\left[5+W_{22}(\mu_0)
       -12\left\{[\mathbf I-\mathbf W(\mu_0)]^{-1}\right\}_{22}\right]
       +\mathcal O(\varepsilon^{5/2}).
       \label{si-eq:ssp-total-coupling}
\end{equation}
\subsubsection{Implicit free surface}

Substituting the backward-Euler barotropic matrix~\eqref{eq:W-im} into~\eqref{si-eq:generic-total-coupling} and~\eqref{si-eq:ssp-total-coupling} gives
\begin{equation}
    \delta\lambda_j^{IM} = - \chi_q^2\mu_0^2\begin{cases}
        \dfrac{\mu_0^2+12}{12(\mu_0^2+4)}  \ \ & \text{(restarted)} \ \\[8pt]
        \dfrac{7\mu_0^2+6}{24(\mu_0^2+1)}  \ \ & \text{(blended)} \
    \end{cases} + \mathcal{O}(\varepsilon^{5/2}) \ .
    \label{si-eq:implicit-coupling}
\end{equation}
Both expressions are real and strictly negative, so coupling with an implicit free surface adds damping at every $\mu_0$.

\subsubsection{Sign of the indirect correction}

The indirect correction is the difference between the totals~\eqref{si-eq:generic-total-coupling} and~\eqref{si-eq:ssp-total-coupling} and the direct corrections~\eqref{si-eq:direct-finite-step} and~\eqref{si-eq:ssp-direct-transfer}. For exact barotropic evolution, $\mathbf W(\mu)=e^{-i\mu\mathbf J}$, it reduces to
\begin{equation}
    \delta\lambda_{j,\mathrm{ind}}^{\mathrm{exact}} \simeq \frac{\chi_q^2}{4}\left[\mu_0\cos\frac{\mu_0}{2}-2\sin\frac{\mu_0}{2}\right]
    \begin{cases}
        \mu_0-2\sin\dfrac{\mu_0}{2} \ , & \text{(restarted)} \ , \\[8pt]
        \mu_0\cos\dfrac{\mu_0}{2}-2\sin\dfrac{\mu_0}{2} \ , & \text{(blended)} \ .
    \end{cases}
    \label{si-eq:exact-indirect}
\end{equation}
Since $\mu_0 > 2\sin(\mu_0/2)$, the blended (SRK3) correction is never negative, while the restarted (WRK3) correction is negative until the first positive root of $\mu_0\cos(\mu_0/2) = 2\sin(\mu_0/2)$, at $\mu_0 \simeq 8.99$. The ratio of the two corrections approaches $-1$ at $\mu_0 = 2\pi$.
Averaging the exact evolution uniformly over $0 \leq \tau \leq 2$, with weight $w_s = 1/2$, gives the barotropic eigenvalue
\begin{equation}
    \lambda_b^{\mathrm{filt}}(\mu_0) = \frac12\int_0^2 e^{-i\mu_0\tau}\,\mathrm{d}\tau = e^{-i\mu_0}\frac{\sin\mu_0}{\mu_0} \ .
\end{equation}
Evaluating the indirect correction at $\mu_0 = 2\pi$ for the averaged and for the exact evolution gives
\begin{equation}
    \delta\lambda_{j,\mathrm{ind}}^{\mathrm{filt}}(2\pi)
    \simeq \frac{\chi_q^2}{2}
    \begin{cases}
        1>0, & \text{WRK3} \ ,\\
        1-\pi^2<0, & \text{SRK3} \ ,
    \end{cases}
    \qquad
    \delta\lambda_{j,\mathrm{ind}}^{\mathrm{exact}}(2\pi)
    \simeq \chi_q^2
    \begin{cases}
        -\pi^2<0, & \text{WRK3} \ ,\\
        \pi^2>0, & \text{SRK3} \ ,
    \end{cases}
    \label{si-eq:averaged-indirect}
\end{equation}
so the averaging reverses the sign of both corrections, while the two schemes keep opposite signs.

\subsection{The small timestep limit}
\label{si-sec:appendix-small-step}
We recover the leading small-timestep direct correction caused by the barotropic substep solver. In contrast to the previous sections, we take $\mu_0$ to zero at fixed $\varepsilon$. For a FB substep scheme, the barotropic substep matrix in the barotropic variables is
\begin{equation}
\mathbf W_B
=\mathbf I
-i\frac{\mu_m\gamma_B}{\gamma_m}\mathbf J
-\left(\frac{\mu_m\gamma_B}{\gamma_m}\right)^2
\begin{bmatrix}0&0\\0&1\end{bmatrix}.
\end{equation}
We raise the above at the $s$-th power to second order in $\mu_m$:
\begin{equation}
\mathbf W_B^s
=\mathbf I+s(\mathbf W_B-\mathbf I)
+\frac{s(s-1)}2(\mathbf W_B-\mathbf I)^2+\mO(\mu_m^3) \ .
\end{equation}
Then, since $\mathbf J^2= \mathbf I$ we have that
\begin{equation}
\mathbf W_B^s
=\mathbf I
-is\frac{\mu_m\gamma_B}{\gamma_m}\mathbf J - \left(\frac{\mu_m\gamma_B}{\gamma_m}\right)^2 \left[  s \left(\begin{matrix} 0 & 0 \\ 0 & 1 \end{matrix} \right)+ \frac{s(s-1)}{2} \mathbf{I}\right] + \mO(\mu_m^3) \ .
\end{equation}
Using $\tau_s = s\gamma_B/\gamma_m$ gives
\begin{equation}
    \mathbf{W_B}^s=\mathbf{I} - i\mu_m \tau_s \mathbf{J} - \frac{\mu_m^2\tau_s^2}{2} \mathbf{I} + \frac{\mu_m^2 \tau_s}{2}\frac{\gamma_B}{\gamma_m}\begin{bmatrix}
        1 & 0 \\ 0 & -1
    \end{bmatrix} + \mathcal{O}(\mu_m^3)
\end{equation}
From the filter's moment conditions, performing $\sum_s w_s \mathbf{W}_B^s$ with $\sum_s w_s \tau_s = 1$ and $\sum_s w_s\tau_s^2 = 1 + M_2$ gives
\begin{equation}
    \mathbf{W}_{FB} = \mathbf{I} - i\mu_m  \mathbf{J} - \frac{\mu_m^2}{2}(1+M_2) \mathbf{I} + \frac{\mu_m^2}{2}\frac{\gamma_B}{\gamma_m}\begin{bmatrix}
        1 & 0 \\ 0 & -1
    \end{bmatrix} + \mathcal{O}(\mu_m^3)
\end{equation}

For RK3 a single substep uses $R(-i \mu_m (\gamma_B/\gamma_m)\mathbf J)$. Using the same procedure we used with the FB filtered matrix yields
\begin{equation}
    \mathbf{W}_{RK3} = \mathbf{I} - i\mu_m \mathbf{J} - \frac{\mu_m^2}{2}(1+M_2) \mathbf{I} + \mathcal{O}(\mu_m^3) \ , 
\end{equation}
while an unfiltered exact evolution has
\begin{equation}
    \mathbf W_{\mathrm{exact}}
=e^{-i\mu_m\mathbf J}
=\mathbf I-i\mu_m\mathbf J-\frac{\mu_m^2}{2}\mathbf I +O(\mu_m^3) \ , 
\end{equation}
showing that an RK3 substep procedure with $M_2 = 0$ agrees with exact evolution through quadratic order in $\mu_m$.

For the constant forcing $F$ used in this calculation, each matrix acts on $[\bar u^n,\eta^n+F]^T$. Subtracting the exact, unfiltered matrix therefore gives the following stage-output differences, where $\delta$ in this paragraph denotes numerical minus exact barotropic evolution:
\begin{equation}
    \begin{aligned}
        & \delta \bar{u}^m = \frac{\mu_m^2}{2} \left(\frac{\gamma_B}{\gamma_m} - M_2 \right)\bar{u}^n \ ,  \ \ && \delta \eta^m = - \frac{\mu_m^2}{2} \left(\frac{\gamma_B}{\gamma_m} + M_2 \right)(\eta^n + F)\ ,  \ \ &&  \text{for FB} \ , \\
        & \delta \bar{u}^m = - \frac{\mu_m^2M_2}{2}\bar{u}^n \ ,  \ \ && \delta \eta^m = - \frac{\mu_m^2M_2}{2}(\eta^n + F)\ ,  \ \ && \text{for RK3} \ , \\
    \end{aligned}
\end{equation}
For $\vect y^n=\vect r_j$ in vertical mode $q$, $\bar{u}^n = \chi_q u_q^n$ and $\eta^n + F = \chi_qh^n_q + \mathcal{O}(\mu_0)$; then the barotropic continuity equation gives $\delta \bar{u}^\dagger = \delta \eta / (-i\mu_m)$ and the velocity reconciliation and pressure update give 
\begin{equation}
    \delta u_q = \chi_q \delta \bar{u} \ , \qquad \delta h_q = - i \mu_m \alpha_q^2 \chi_q \delta \bar{u}^\dagger = \alpha_q^2\chi_q \delta \eta \ .
\end{equation}
To pass from a stage-output difference to $\mathbf V$, first consider the restarted final stage. With exact, unfiltered evolution and $\mu_{M-1}=\mu_0/2$, the preceding-stage values satisfy
\begin{equation}
 \vect u^{M-1}=\vect u^n-\frac{i\mu_0}{2}\vect h^n+\mathcal O(\mu_0^2),\qquad
 \vect h^{M-1}=\vect h^n-\frac{i\mu_0}{2}\mathbf D\vect u^n+\mathcal O(\mu_0^2).
\end{equation}
For the frozen stage forcing, expanding the exact barotropic outputs and the velocity predictor gives
\begin{equation}
\begin{aligned}
 \bar u^M&=\bar u^n-i\mu_0\vect\chi^T\vect h^n
 -\frac{\mu_0^2}{2}(\bar u^n+\vect\Pi^T\mathbf D\vect u^n)+\mathcal O(\mu_0^3),\\
 \vect\chi^T\vect u^*&=\bar u^n-i\mu_0\vect\chi^T\vect h^n
 -\frac{\mu_0^2}{2}\vect\chi^T\mathbf D\vect u^n+\mathcal O(\mu_0^3),\\
 \bar u^{\dagger,M}&=\bar u^n-\frac{i\mu_0}{2}\vect\chi^T\vect h^n+\mathcal O(\mu_0^2).
\end{aligned}
\end{equation}
Since $\vect\chi=\vect\psi+\vect\Pi$ and $\vect\psi^T\mathbf D=\vect\chi^T$, the two velocity expressions agree through order $\mu_0^2$. The transport difference $\bar u^{\dagger,M}-\vect\chi^T\vect u^{M-1}=\mathcal O(\mu_0^2)$ enters the pressure update multiplied by $\mu_0$, contributing only $\mathcal O(\mu_0^3)$. Thus the exact, unfiltered coupled update agrees with the uncoupled update through order $\mu_0^2$.
The numerical/filter differences above supply the order-$\mu_0^2$ contribution to $\mathbf V$. Errors from earlier stages are order $\mu_0^2$ and enter a subsequent restarted update multiplied by $\mu_0$, so they do not enter this coefficient.

At the final stage $\mu_m=\mu_0$ and $\gamma_m=\Delta t$. The two output differences above are diagonal in the initial velocity and pressure; off-diagonal terms start at order $\mu_0^3$. Hence the retained modal block is
\begin{equation}
    \left.\mathbf V_{qq}\right|_{\text{diagonal}}
    =
    \dfrac{\chi_q^2\mu_0^2}{2} \begin{cases}
        
        \operatorname{diag}\left(
            \dfrac{\gamma_B}{\Delta t}-M_2,\,
            -\alpha_q^2\left[\dfrac{\gamma_B}{\Delta t}+M_2\right]
        \right),
        & \text{(FB)},\\[5pt]
        -
        M_2\operatorname{diag}\left(1,\alpha_q^2\right),
        & \text{(RK3)}
    \end{cases}
    +\mathcal O(\mu_0^3).
\end{equation}
The left and right eigenvectors give equal weights $1/2$ to the two diagonal entries. Projecting this block therefore gives
\begin{equation}
    V_{jj}=\vect{\ell}_j\mathbf V_{qq}\vect r_j
    =
    \dfrac{\chi_q^2\mu_0^2}{4} \begin{cases}
            (1-\alpha_q^2)\dfrac{\gamma_B}{\Delta t}
                        -(1+\alpha_q^2)M_2 \ , 
        & \text{(FB)},\\[5pt]
        -(1+\alpha_q^2)M_2 \ ,
        & \text{(RK3)}
    \end{cases}
    +\mathcal O(\mu_0^3).
\end{equation}
This is the direct correction for the restarted schemes with frozen stage forcing. Its real quadratic term enters the magnitude correction through~\eqref{eq:linearized-eig-equation}; the distinction between FB and RK3 comes from the term proportional to $\gamma_B$. For FB, nonpositive quadratic growth requires $(1+\alpha_q^2)M_2\ge(1-\alpha_q^2){\gamma_B}/{\Delta t}$. For weak stratification this reduces to $M_2\gtrsim{\gamma_B}/{\Delta t}$, the filter condition discussed in section~\ref{sec:4}.
\section{Coriolis forcing in the split--explicit stages}
\label{si-sec:appendix-rotation}

We derive the rotating modal matrices compared in Appendix~\ref{sec:rotation-results}. We consider three WRK3 (and MRK4) cases which differ in the reconstruction of pressure and Coriolis forcing while SRK3 holds both frozen at the preceding stage value.
Taking the horizontal Fourier wavevector along $x$ gives
\begin{equation}
    \frac{d\vect u}{dt}=f\vect v-ik\sqrt{gH}\vect h,\qquad
    \frac{d\vect v}{dt}=-f\vect u,\qquad
    \frac{d\vect h}{dt}=-ik\sqrt{gH}\mathbf D\vect u.
    \label{si-eq:rot-modal}
\end{equation}
Here $\vect v$ contains the transverse velocity amplitudes. Each vertical mode has a stationary geostrophic solution and two inertia--gravity waves. Their exact amplification factors are $1$ and $\exp[\pm i\sqrt{\mu_0^2\alpha_q^2+(f\Delta t)^2}]$. The wavevector choice loses no horizontal direction in this isotropic Fourier problem: the longitudinal velocity feels the pressure gradient, and Coriolis couples it to the transverse velocity.

Using $\vect\chi^T\vect h=\eta+\vect\Pi^T\vect h$ and $\vect\psi^T\mathbf D=\vect\chi^T$, the depth-averaged equations are
\begin{equation}
\begin{aligned}
    \frac{d\bar u}{dt'}&=f\widehat{\bar v}^{\,m}(t')-i\frac{\mu_0}{\Delta t}[\eta+F^m(t')] \ ,\\
    \frac{d\bar v}{dt'}&=-f\widehat{\bar u}^{\,m}(t') \ ,\\
    \frac{d\eta}{dt'} & =-i\frac{\mu_0}{\Delta t}\bar u.
\end{aligned}
\label{si-eq:rot-barotropic}
\end{equation}
The hats specify the velocities supplied to Coriolis from the stage values.

We use the forcing reconstructions defined in section~\ref{sec:3-mrk4}. The three choices differ in which stage values are reconstructed:
\begin{itemize}
    \item \emph{Frozen forcing}: $F^m=\vect\Pi^T\vect h^{m-1}$, $\widehat{\bar u}^{\,m}=\vect\chi^T\vect u^{m-1}$, and $\widehat{\bar v}^{\,m}=\vect\chi^T\vect v^{m-1}$.
    \item \emph{Progressive pressure}: reconstruct $F^m$ from the available values of $\vect\Pi^T\vect h^j$, while keeping both Coriolis velocities frozen at stage $m-1$.
    \item \emph{Progressive pressure and Coriolis}: apply the same reconstruction to $\vect\Pi^T\vect h^j$, $\vect\chi^T\vect u^j$, and $\vect\chi^T\vect v^j$.
\end{itemize}
To express these quantities in terms of $\vect y^n$, we use $\vect y^j=\mathbf M_j\vect y^n$, with $\mathbf M_0=\mathbf I$ and $\vect y=(\vect u,\vect v,\vect h)^T$. We define selectors as
\begin{equation}
    \mathbf E_u=[\,\mathbf I\ \mathbf0\ \mathbf0\,],\qquad
    \mathbf E_v=[\,\mathbf0\ \mathbf I\ \mathbf0\,],\qquad
    \mathbf E_h=[\,\mathbf0\ \mathbf0\ \mathbf I\,]
    \label{si-eq:rot-selectors}
\end{equation}
that allows us to express the forcings as
\begin{equation}
\begin{aligned}
    \vect\Pi^T\vect h^j&=\vect\Pi^T\mathbf E_h\mathbf M_j\vect y^n,\\
    \vect\chi^T\vect u^j&=\vect\chi^T\mathbf E_u\mathbf M_j\vect y^n,\\
    \vect\chi^T\vect v^j&=\vect\chi^T\mathbf E_v\mathbf M_j\vect y^n.
\end{aligned}
\label{si-eq:rot-stage-values}
\end{equation}
Substituting these expressions into the reconstruction formulas and factoring out $\vect y^n$ gives
\begin{equation}
\begin{aligned}
    F^m(t')&=\sum_p t'^p\vect{\mathcal F}_p^{\,m}\vect y^n,\\
    \widehat{\bar u}^{\,m}(t')&=\sum_p t'^p\vect{\mathcal U}_p^{\,m}\vect y^n,\\
    \widehat{\bar v}^{\,m}(t')&=\sum_p t'^p\vect{\mathcal V}_p^{\,m}\vect y^n.
\end{aligned}
\label{si-eq:rot-reconstruction}
\end{equation}
Here $\vect{\mathcal F}_p^{\,m}$ is the forcing row defined in section~\ref{si-sec:app-stage-matrices}; $\vect{\mathcal U}_p^{\,m}$ and $\vect{\mathcal V}_p^{\,m}$ are obtained in the same way from the two velocity samples. For frozen forcing, only $p=0$ is nonzero:
\begin{equation}
\begin{aligned}
    \vect{\mathcal F}_0^{\,m}&=\vect\Pi^T\mathbf E_h\mathbf M_{m-1},\\
    \vect{\mathcal U}_0^{\,m}&=\vect\chi^T\mathbf E_u\mathbf M_{m-1},\\
    \vect{\mathcal V}_0^{\,m}&=\vect\chi^T\mathbf E_v\mathbf M_{m-1}.
\end{aligned}
\label{si-eq:rot-cases}
\end{equation}
We write the barotropic solution as
\begin{equation}
    \begin{bmatrix}\bar u(t')\\\bar v(t')\\\eta(t')\end{bmatrix}
    =\mathbf X_m(t')\vect y^n.
    \label{si-eq:rot-output-definition}
\end{equation}
Substitution into~\eqref{si-eq:rot-barotropic}, followed by factoring out $\vect y^n$, gives
\begin{equation}
\begin{aligned}
 \frac{d\mathbf X_m}{dt'}&=
 \begin{bmatrix}0&0&-i\mu_0/\Delta t\\0&0&0\\-i\mu_0/\Delta t&0&0\end{bmatrix}\mathbf X_m
 +\sum_p t'^p\begin{bmatrix}
 f\vect{\mathcal V}_p^{\,m}-i(\mu_0/\Delta t)\vect{\mathcal F}_p^{\,m}\\
 -f\vect{\mathcal U}_p^{\,m}\\\vect0^T
 \end{bmatrix},\\
 \mathbf X_m(0)&=\begin{bmatrix}\vect\chi^T\mathbf E_u\\\vect\chi^T\mathbf E_v\\\vect\psi^T\mathbf E_h\end{bmatrix}.
\end{aligned}
\label{si-eq:rot-output-matrix}
\end{equation}
Each restarted solve uses the same initial state; its stage dependence enters through the forcing coefficients and duration $\gamma_m$.

We integrate~\eqref{si-eq:rot-output-matrix} with RK3 barotropic substeps, evaluating the prescribed forcing at each substep midpoint, and apply the F3 averaging filter. Denoting the three filtered output rows by $\mathbf X_{m,1}$, $\mathbf X_{m,2}$, and $\mathbf X_{m,3}$ gives
\begin{equation}
    \bar u^m=\mathbf X_{m,1}\vect y^n,\qquad
    \bar v^m=\mathbf X_{m,2}\vect y^n,\qquad
    \eta^m=\mathbf X_{m,3}\vect y^n.
    \label{si-eq:rot-filtered-outputs}
\end{equation}
Using $\eta^m-\eta^n=-i\mu_m\bar u^{\dagger,m}$, continuity gives the corresponding transport:
\begin{equation}
    \bar u^{\dagger,m}=\frac{\mathbf X_{m,3}-\vect\psi^T\mathbf E_h}{-i\mu_m}\vect y^n.
    \label{si-eq:rot-transport-row}
\end{equation}
The velocity predictors are
\begin{equation}
\begin{aligned}
    \vect u^*&=\vect u^n+f\gamma_m\vect v^{m-1}-i\mu_m\vect h^{m-1},\\
    \vect v^*&=\vect v^n-f\gamma_m\vect u^{m-1}.
\end{aligned}
\label{si-eq:rot-predictors}
\end{equation}
Reconciliation and the pressure update, as in~\eqref{si-eq:app-stage-collected}, give
\begin{equation}
\begin{aligned}
    \vect u^m&=\boldsymbol{\mathcal P}\vect u^*+\vect\chi\bar u^m,\\
    \vect v^m&=\boldsymbol{\mathcal P}\vect v^*+\vect\chi\bar v^m,\\
    \vect h^m&=\vect h^n-i\mu_m\mathbf D\left(\boldsymbol{\mathcal P}\vect u^{m-1}+\vect\chi\bar u^{\dagger,m}\right).
\end{aligned}
\label{si-eq:rot-reconciled-updates}
\end{equation}
Substituting~\eqref{si-eq:rot-predictors}, \eqref{si-eq:rot-filtered-outputs}, and~\eqref{si-eq:rot-transport-row} into these updates and using $\vect y^{m-1}=\mathbf M_{m-1}\vect y^n$ gives
\begin{equation}
\mathbf M_m=\begin{bmatrix}
 \boldsymbol{\mathcal P}(\mathbf E_u+f\gamma_m\mathbf E_v\mathbf M_{m-1}-i\mu_m\mathbf E_h\mathbf M_{m-1})+\vect\chi\mathbf X_{m,1}\\[3pt]
 \boldsymbol{\mathcal P}(\mathbf E_v-f\gamma_m\mathbf E_u\mathbf M_{m-1})+\vect\chi\mathbf X_{m,2}\\[3pt]
 \mathbf E_h-i\mu_m\mathbf D\boldsymbol{\mathcal P}\mathbf E_u\mathbf M_{m-1}
       +\mathbf D\vect\chi(\mathbf X_{m,3}-\vect\psi^T\mathbf E_h)
\end{bmatrix}.
\label{si-eq:rot-stage}
\end{equation}
For WRK3, $(\gamma_1,\gamma_2,\gamma_3)=\Delta t(1/3,1/2,1)$, while for MRK4  $(\gamma_1,\gamma_2,\gamma_3,\gamma_4)=\Delta t(a,1/3,1/2,1)$. Applying this construction successively gives $\mathbf M_1$, $\mathbf M_2$, and $\mathbf M_3$ (and $\mathbf M_4$) for each forcing treatment. At $f=0$, each recovers its nonrotating $(\vect u,\vect h)$ matrix, and the transverse velocity remains constant.

For SRK3, we first evaluate the restarted-stage matrix in \eqref{si-eq:rot-stage} with $\gamma_m=\Delta t$ and the frozen stage forcing in \eqref{si-eq:rot-cases}. We then apply the increment blend:
\begin{equation}
    \mathbf M_m \longleftarrow
    (1-b_m)\mathbf I
    +b_m\left(\mathbf M_{m-1}+\mathbf M_m-\mathbf I\right),
    \qquad (b_1,b_2,b_3)=(1,1/4,2/3).
    \label{si-eq:rot-ssp-matrix}
\end{equation}

\bibliography{mybibfile}
\end{document}